\documentclass[11pt, reqno]{amsart}

\usepackage{amsmath, amssymb, comment, amsthm, amscd}
\usepackage{graphicx}
\usepackage{ragged2e}
\usepackage{epstopdf}
\usepackage{setspace}
\usepackage{abstract}
\usepackage{subfigure}
\usepackage{color}
\usepackage[foot]{amsaddr}

\usepackage{enumitem}
\usepackage[percent]{overpic}

\usepackage[pdftex, linktocpage=true]{hyperref}

\DeclareFontFamily{U}{mathx}{\hyphenchar\font45}
\DeclareFontShape{U}{mathx}{m}{n}{
      <5> <6> <7> <8> <9> <10>
      <10.95> <12> <14.4> <17.28> <20.74> <24.88>
      mathx10
      }{}
\DeclareSymbolFont{mathx}{U}{mathx}{m}{n}
\DeclareFontSubstitution{U}{mathx}{m}{n}
\DeclareMathAccent{\widecheck} {0}{mathx}{"71}

\usepackage[format=plain,labelfont=bf,up,width=.99\textwidth]{caption} 

\theoremstyle{plain}
\newtheorem{thm}{Theorem}[section]
\newtheorem{lemma}[thm]{Lemma}
\newtheorem{prop}[thm]{Proposition}
\newtheorem{cor}{Corollary}[thm]

\theoremstyle{definition}
\newtheorem{defn}[thm]{Definition}

\newtheorem{remarka}[thm]{Remark}

\newcommand{\R}{\mathbb{R}}
\newcommand{\N}{\mathbb{N}}
\newcommand{\C}{\mathbb{C}}

\newcommand{\D}{\mathbb{D}}

\newcommand{\bigO}{\mathcal{O}}

\newcommand{\fq}{\mathbf{q}}
\newcommand{\s}{\mathbb{S}}
\newcommand{\Pet}{\mathcal{P}}
\newcommand{\ds}{\displaystyle}

\def \hvec #1#2{\left(#1, #2\right)}

\def \diag #1#2#3#4#5#6#7#8{\begin{CD} #1 @>#5>> #2\\ @V#6VV  @VV#7V\\ #3 @>#8>> #4 \end{CD}}
\def\He{{H\'enon }}

\def\limind{{\underrightarrow{\lim}}}

\def\proof{\par\medskip\noindent {\bf Proof.\ \ }}
\def\qed{\hfill $\square$\\ }

\title[Semi-parabolic tools for hyperbolic H\'enon maps]{Semi-parabolic tools for hyperbolic H\'enon maps and continuity of Julia sets in $\mathbb{C}^{2}$}
\author[Remus Radu and Raluca Tanase]{Remus Radu$^1$ and Raluca Tanase$^1$}
\address{$^1$Institute for Mathematical Sciences, Stony Brook University, Stony Brook, NY 11794-3660}
\email{remus.radu@stonybrook.edu, raluca.tanase@stonybrook.edu}

\date{\today}

\begin{document}
\maketitle

\begin{abstract}
\noindent {\sc abstract.} 
We prove some new continuity results for the Julia sets $J$ and $J^{+}$ of the complex H\'enon map $H_{c,a}(x,y)=(x^{2}+c+ay, ax)$, where $a$ and $c$ are complex parameters. We look at the parameter space of  dissipative H\'enon maps which have a fixed point with one eigenvalue $(1+t)\lambda$, where $\lambda$ is a root of unity and $t$ is real and small in absolute value. These maps have a semi-parabolic fixed point when $t$ is $0$, and we use the techniques that we have developed in \cite{RT} for the semi-parabolic case to describe nearby perturbations. 
We show that for small nonzero $|t|$, the H\'enon map is hyperbolic and has connected Julia set.  We prove that the Julia sets $J$ and $J^{+}$ depend continuously on the parameters as $t\rightarrow 0$, which is a two-dimensional analogue of radial convergence from one-dimensional dynamics. Moreover, we prove that this family of H\'enon maps is stable on $J$ and $J^{+}$ when $t$ is nonnegative.
\end{abstract}

\tableofcontents

\section{Introduction}\label{sec:intro}

Complex analytic maps with a parabolic fixed point or cycle have generated much interest in dynamics in one complex variable as they play a fundamental role in understanding the parameter space of rational maps. Moreover, they provide important models for understanding non-hyperbolic behavior. 

In \cite{RT} we studied the family of \He maps with a semi-parabolic fixed point or cycle, and showed that the family has nice stability properties. In this paper we want to unravel the mystery about how these semi-parabolic maps sit in the parameter space of \He maps and describe the Julia sets of nearby perturbations.  

Chapter \ref{ch:Continuity1} provides a useful digression to
dynamics in one complex variable. Consider a quadratic polynomial $p(x)=x^2+c$ with a parabolic fixed point and denote by
$J_p$ its Julia set. The parameter $c$ lies in the boundary of the
Mandelbrot set. It is well-known from the work of P. Lavaurs \cite{L} and  A. Douady
\cite{D} that on a neighborhood of the parameter $c$ in $\C$ the Julia set does not vary continuously in the Hausdorff topology. Parabolic implosion represents the source of discontinuity and of obtaining limit Julia sets with enriched dynamics. Using quasiconformal techniques, C. McMullen
\cite{Mc} (and also P. Ha\"issinsky \cite{Hai}) showed that $J_{p_n}$ converges to $J_p$ when $p_n$ converges to $p$ horocyclically or radially (i.e. non-tangentially with respect to the boundary of Mandelbrot set in the quadratic case). These tools are harder, if not impossible, to apply to several complex variables, where an analogue of the Uniformization Theorem does not exist. We first set out to give a topological proof of the continuity result for polynomial Julia sets under a stronger radial convergence assumption. The proof involves recovering the Julia set as the image of the unique
fixed point $f^*$ of a (weakly) contracting operator in an appropriate function
space. In Section \ref{sec:contraction} we prove that $f^*$
depends continuously on the parameter, and thus the corresponding
Julia sets converge to the Julia set of the parabolic polynomial, in the Hausdorff topology. 
After gaining some valuable insight from the study of the one dimensional problem, notably from Lemma 
~\ref{lemma: Browder-function-h}, we pursue the two-dimensional problem and prove 
some new  continuity results in Chapter \ref{ch:Henon}
for the Julia sets $J$ and $J^+$ of a complex \He map.

We consider the family of complex H\'enon maps $H_{c,a}(x,y) = (p(x)+ay, ax)$, where $p$ is a quadratic polynomial, $p(x)=x^2+c$, and $a$ is a complex parameter. When $a\neq 0$, this is a polynomial automorphism of $\C^2$. The dynamics of \He maps bears some resemblance to the dynamics of 1-D polynomials, however extending results from one to several variables requires envisioning new techniques and approaches, 
and in many cases the emerging picture is substantially different and contains new and thrilling phenomena not present in the one-dimensional world. In order to describe the dynamics of the \He map, one studies the sets $K^{+}$ and $K^{-}$ of points which do not escape to infinity under forward and respectively backward iterations. The topological complements of $K^{\pm}$ in $\C^2$ are denoted by $U^{\pm}$ and called the escaping sets. The most interesting dynamics occurs on the boundary of the sets $K^{\pm}$ and $U^{\pm}$ where chaotic behavior is present. The sets $J^{\pm}=\partial K^{\pm}=\partial U^{\pm}$ and $J=J^{+}\cap J^{-}$ are called the Julia sets of the \He map. The sets $J$ and $K= K^{+}\cap K^{-}$ are compact, while the sets $J^{\pm}$ are closed, connected and unbounded \cite{BS1}. 

A quadratic \He map is uniquely determined by the eigenvalues $\lambda$ and $\nu$ at a fixed point so we will sometimes write $H_{\lambda, \nu}$ in place of $H_{c,a}$ to mark this dependence. The precise formula for $H_{\lambda,\nu}$ is given in at the beginning of Chapter \ref{ch:Henon}. 
We say that a
\He map is semi-parabolic if it has a fixed point (or cycle) with one
eigenvalue $\lambda$, a  root of unity, and one eigenvalue
smaller than one in absolute value. 
Unlike hyperbolic \He maps,
semi-parabolic ones are not stable under
perturbations. E. Bedford, J. Smillie and T. Ueda have described some
semi-parabolic bifurcations in $\C^2$ for $\lambda=1$ in \cite{BSU}. 
In particular, they show that at a parameter value with a semi-parabolic fixed point with the eigenvalues $\lambda=1$ and $|\nu| < 1$, the sets $J$, $J^+$, $K$ and $K^+$ vary discontinuously with the parameters, while $J^{-}$ and $K^{-}$ vary continuously with the parameters. The phenomenon described in \cite{BSU} is a two-dimensional analogue of parabolic implosion that occurs in complex dimension one.

In order to state our main results, consider a primitive root of unity $\lambda=e^{2\pi i p/q}$ and let $\lambda_t=(1+t)\lambda$.  For $t$ real and 
small in absolute value, we look at the parameter space
$\mathcal{P}_{\lambda_t}$ of complex \He maps which have a fixed point with
one eigenvalue $\lambda_t$. The equation of the curve $\mathcal{P}_{\lambda_t}$ is given in Proposition \ref{prop:Plambda}. When $t=0$ these maps are semi-parabolic; when $t\neq 0$, we regard the maps corresponding to parameters from  $\mathcal{P}_{\lambda_t}$ as perturbations of the semi-parabolic ones. We show in Section \ref{subsec:JJ+} that there exists $\delta>0$ such that for $(c,a)\in \mathcal{P}_{\lambda_t}$ and $0<|a|<\delta$,  the Julia sets $J$ and $J^{+}$
depend continuously on the parameters as $t$ approaches $0$. An equivalent formulation is given in the theorem below.
These results can be regarded as a natural extension of the concept of radial convergence of Julia sets \cite{Mc} to higher dimensions, in the context of polynomial automorphisms of $\C^2$.

\begin{thm}[\textbf{Continuity}]\label{thm:continuity}
There exists $\delta>0$ such that if $|\nu_t|<\delta$ and $\nu_t\rightarrow \nu$ as $t\rightarrow 0$, then 
the Julia sets $J$ and $J^+$ depend continuously on the parameters, i.e. 
\[
    J^{+}_{(\lambda_{t},\nu_t)}\rightarrow J^{+}_{(\lambda,\nu)}\ \ \ \ \ \mbox{and}\ \ \ \     J_{(\lambda_t,\nu_t)}\rightarrow J_{(\lambda,\nu)}
\]
in the Hausdorff topology.
\end{thm}

For the set $J^{+}$ we are taking the Hausdorff topology on the one-point compactification of  $\C^2$. What we prove in Theorem \ref{thm:continuity} is the continuity of Julia sets $J$ and $J^+$ as we approach a semi-parabolic parameter from the interior of a hyperbolic component of the \He connectedness locus, similar to radial convergence from 1-D dynamics. Our next theorem describes the dynamical nature of the perturbed semi-parabolic maps.

\begin{thm}[\textbf{Hyperbolicity}]\label{thm:HypRegion}
There
exist $\delta,\delta'>0$ such that in the parametric region
\begin{eqnarray*}
\mathcal{HR}_{\delta, \delta'}=\left\{(c,a)\in \mathcal{P}_{\lambda_t}\ :\ 0<|a|<\delta\ \ \mbox{and}\
 -\delta'<t<\delta',\  t\neq 0\right\}
\end{eqnarray*}
the Julia set $J_{c,a}$ is connected and the \He map $H_{c,a}$ is hyperbolic.
\end{thm}

By definition, the connectedness locus for the \He family is the set of parameters $(c,a)\in\C^{2}$ such that the Julia set $J_{c,a}$ is connected. 
Theorem \ref{thm:HypRegion} shows that the parametric region $\{(c,a)\in P_{\lambda}: |a|<\delta\}$ of semi-parabolic \He maps lies in the boundary of a hyperbolic component of the \He connectedness locus. In fact, when $\lambda\neq 1$, it lies in the boundary of two such hyperbolic components.  A mechanism for loss of hyperbolicty at the boundary of the horseshoe region through the development of tangencies between stable and unstable manifolds is described by E. Bedford and J. Smillie in \cite{BS}, and more recently by Z. Arai and Y. Ishii  in \cite{AI}. In Theorem \ref{thm:HypRegion} we describe a different mechanism for loss of hyperbolicity, through the creation of a semi-parabolic fixed point. 
We first do a local analysis in Sections  \ref{sec:localdynamics}, \ref{subsec:horizontal} and \ref{section:parabolic2hyperbolic}, and show how to deform the local semi-parabolic structure into a hyperbolic structure; these sections are applicable to holomorphic germs of diffeomorphisms of $(\C^2,0)$ with a semi-parabolic fixed point, and their perturbations. 
We complete the proof of Theorem \ref{thm:HypRegion} in Section \ref{subsec:hyper-region}.

Theorem \ref{thm:HypRegion} proves the existence of a larger region of hyperbolicity for complex \He maps than what was previously known.  It is in general very hard to exhibit regions of hyperbolicity for \He maps. Z. Arai developed a computer program for detecting hyperbolicity, that relies on heavy numerical computations. Otherwise, the only  \He maps proven to be hyperbolic using only theoretical arguments correspond to the horseshoe region  and to perturbations of 1-D hyperbolic maps. It is also known from early work of J. Hubbard and R. Oberste-Vorth \cite{HOV1, HOV2} 
(see also J.E. Forn{\ae}ss and N. Sibony \cite{FS}) that \He maps that come from perturbations of hyperbolic polynomials with connected Julia sets inherit both of these properties. However, the proof gave no control on the admissible size of perturbations as we approach the boundary of the Mandelbrot set, i.e. it was not known that the size of the region $\mathcal{HR}_{\delta, \delta'}$ does not decrease to $0$ as $t\rightarrow 0$. Z. Arai \cite{A} gave  a computer-assisted proof for the existence of hyperbolic plateaus for the family of complex H\'enon maps $H_{c,a}$ with both parameters $c$ and $a$ real. In our language, these regions corresponds to strips on the right/left side of the real curves $\mathcal{P}^{n}_{\pm 1}\cap\R^{2}$, where $\mathcal{P}^{n}_{\pm 1}$ is the set of parameters $(c,a)\in \C^{2}$ for which the  \He map $H_{c,a}$ has a cycle of order $n$ with one eigenvalue $\pm 1$. The existence of these regions is established in Theorem \ref{thm:HypRegion}.

\begin{cor}\label{cor:realhyp} There exists an $\epsilon>0$ such that the real parametric region
\begin{equation*}\label{eq:Be}
\left\{(c,a)\in \R\times(-\epsilon,\epsilon)\ :\ a\neq 0,\ \frac{(1+a)^2}{4}-\epsilon<c<\frac{(1+a)^2}{4}\right\}
\end{equation*}
is a region of hyperbolicity for the \He family $H_{c,a}(x,y) = (x^2+c-ay,x)$. The \He map is written in the standard parametrization and has 
Jacobian $a$.  
\end{cor}

To compare our results with \cite{B} and \cite{BSU}, suppose $\lambda$ is $1$. Theorems  \ref{thm:continuity} and \ref{thm:HypRegion} answer Questions $3$ and $4$ of E. Bedford, from \cite{B}.
Corollary \ref{cor:realhyp} is formulated as a specific answer to Question 3 and the set $c=(a+1)^2/4$ from the corollary is simply the defining equation of parabola $\mathcal{P}_1$.  When $(c,a)\in \mathcal{P}_1$, the \He map has a double fixed point with one eigenvalue $1$. From the ``right" of the real parabola $\mathcal{P}_1\cap \R^2$ we have semi-parabolic implosion described in \cite{BSU}. More specifically, in \cite{BSU} it is shown that there exists a sequence $\epsilon_n\rightarrow 0$ which converges to $0$ tangentially to the positive real axis ($Re(\epsilon_n)>0$ and $Im(\epsilon_n)<{\rm const.}\, |\epsilon_n|^2$) such that the Julia set $J_{c_n,a}$ corresponding to the sequence $c_n=(a+1)^2/4+\epsilon_n$ does not converge to the Julia set $J_{c,a}$ in the Hausdorff topology. By comparison, Theorem \ref{thm:continuity} shows that we have continuity of $J$ and $J^{+}$ from the ``left" of the real parabola $\mathcal{P}_1\cap \R^2$. When $(c,a)\in \mathcal{P}_1\cap \R^2$ and $0<|a|<\delta$ we get that
$J_{c-\epsilon,a}\rightarrow J_{c,a}$ and $J^+_{c-\epsilon,a}\rightarrow J^+_{c,a}$ as $\epsilon\rightarrow 0^+$. Note that Theorem \ref{thm:continuity} gives no information on what happens to the ``right'' of parabola $\mathcal{P}_1$; indeed, when $\lambda=1$, both curves $\mathcal{P}_{1+t}$ and $\mathcal{P}_{1-t}$ are to the left of $\mathcal{P}_{1}$.  This can be seen from the fact that for $a=0$ and $t\neq 0$, the polynomial $p_t$ has two distinct fixed points, with multipliers $1\pm t$; therefore regardless of whether $t$ is positive or negative, $p_t$ has an attracting fixed point, and belongs to the interior of the main cardioid of the Mandelbrot set.

\begin{figure}[htb]
\begin{center}
\includegraphics[scale=0.95]{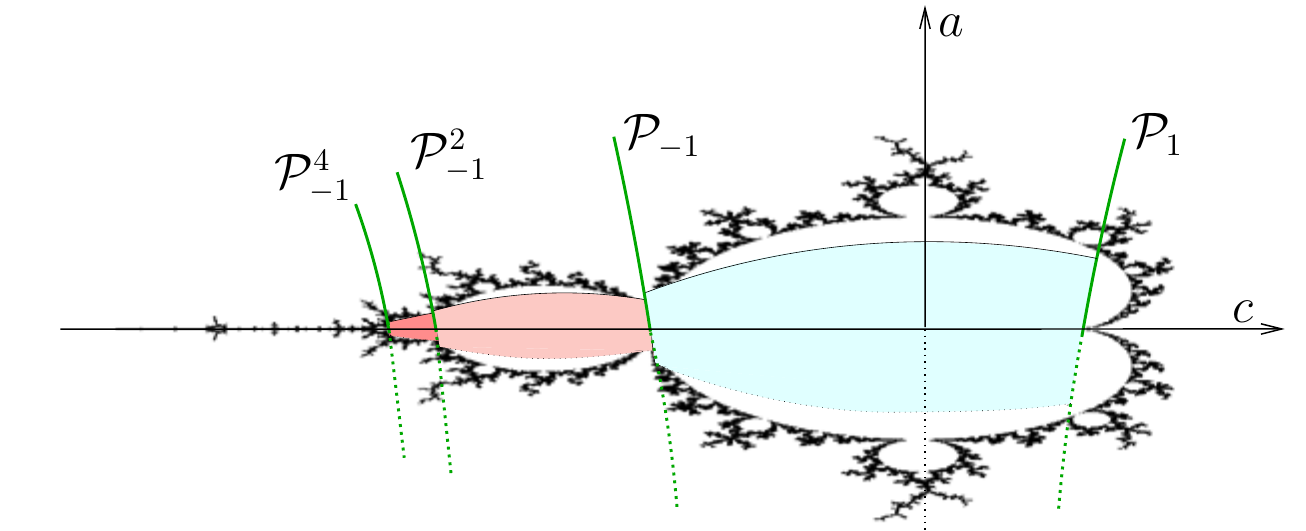}
\end{center}
\caption{The curves $\mathcal{P}_1$, $\mathcal{P}_{-1}$,  $\mathcal{P}_{-1}^{2}$, $\mathcal{P}_{-1}^{4}$ are examples of parametric curves containing semi-parabolic \He maps. There exist regions of hyperbolicity between $\mathcal{P}_1$ and $\mathcal{P}_{-1}$  and $\mathcal{P}_{-1}$ and  $\mathcal{P}_{-1}^{2}$  which belong to hyperbolic components of the \He connectedness locus. }
\label{fig:perspective}
\end{figure}

This paper is built on previous work done by the authors in \cite{R} and \cite{T}. We use the tools that we have developed for the study of semi-parabolic germs/\He maps in \cite{RT} to extend the results from \cite{RT} to nearby perturbations of semi-parabolic germs/\He maps. We can actually say more about the stability properties of our family of \He maps when the parameter $t$ is non-negative:

\begin{thm}[\textbf{Stability}]\label{thm:stability} The family of complex \He maps $\mathcal{P}_{\lambda_t} \ni (c,a) \to  H_{c,a}$   is a structurally stable family  on $J$ and $J^{+}$ for  $0<|a| < \delta$ and $0\leq t<\delta'$.
\end{thm}

We say that the family of \He maps $\mathcal{P}_{\lambda_t} \ni (c,a) \to  H_{c,a}$ is structurally stable on $J$ when $t\in [0,\delta')$ and $|a|<\delta$ if for any two pairs $(c_{i},a_{i})\in \mathcal{P}_{\lambda_{t_{i}}}$, with $|a_{i}|<\delta$ and $t_i\in[0,\delta')$ for $i=1,2$, we have $(H_{c_{1},a_{1}}, J_{c_{1},a_{1}})$ conjugate to $(H_{c_{2},a_{2}}, J_{c_{2},a_{2}})$. Consequently, the Julia sets $J_{c_{1},a_{1}}$ and $J_{c_{2},a_{2}}$ are homeomorphic. We explain structural stability on $J^{+}$ in a similar way, by replacing $J$ with $J^{+}$.  We complete the proof of Theorem \ref{thm:stability} in Section \ref{subsec:JJ+}.

Another notion of stability (called weak stability) was introduced by R. Dujardin and M. Lyubich \cite{DL} for  holomorphic families $(f_{z})_{z\in\Lambda}$ of moderately dissipative polynomial automorphisms of $\C^2$, where $\Lambda$ is a connected complex manifold. Weak stability is defined in terms of branched holomorphic motions of the set $J^*$ (the closure of the saddle periodic points), 
but an equivalent easier definition is the following: 
the family is weakly stable if periodic orbits do not bifurcate. The equivalence between weak stability and continuity of $J^*$ is discussed in \cite{DL}, and 
the relation between weak stability and uniform hyperbolicity on $J^*$ is analyzed by P. Berger and R. Dujardin in \cite{BD}. These results do not apply to our context, but they are of independent interest.

\medskip
\noindent {\it Acknowledgements.} We warmly thank John Hubbard and John Smillie for their guidance and for many useful discussions about \He maps and suggestions on this project.

\section{Continuity of polynomial Julia sets}\label{ch:Continuity1}

In this section we focus only on one-dimensional dynamics. We first discuss continuity of polynomial Julia sets, which will prove useful in understanding continuity of Julia sets for \He maps. This will be treated in Section \ref{ch:Henon}.

Assume that $p$ is a quadratic polynomial. The {\it filled Julia set} of the polynomial $p$ is
\[
    K_p = \{z\in \C\ :\ |p^{\circ n}(z)|\ \mbox{bounded as}\ n\rightarrow \infty \},
\]
and the {\it Julia set} of $p$ is $J_p=\partial K_p$. The filled Julia set $K_{p}$ is connected iff the orbit of the unique critical point is bounded.
If $K_p$ is connected (or equivalently $J_{p}$ is connected) then there exists a unique analytic isomorphism
\begin{equation}\label{eq:psi-p}
    \Psi_{p}:\C-\overline{\D} \rightarrow \C-K_p
\end{equation}
such that $\Psi_{p}(z^2)=p(\Psi_{p}(z))$ and $\Psi_{p}(z)/z\rightarrow 1$ as
$z\rightarrow \infty$. 
If $J_p$ is locally connected then $\Psi_{p}$ extends to the boundary $\s^{1}$ and defines a continuous surjection (see \cite{Mi})
\begin{equation}\label{eq:gamma-p}
\gamma_{p}:\s^{1} \rightarrow J_p.
\end{equation}
The Julia set of a hyperbolic or parabolic polynomial is connected and locally connected (see \cite{DH}). The map $\Psi_{p}^{-1}$ is the {\it B\"ottcher coordinate} of the polynomial $p$, while the map $\Psi_{p}$  is called the {\it inverse B\"ottcher isomorphism} (or the {\it B\"ottcher  chart} \cite{H}). The boundary map $\gamma_{p}$ is called the {\it Carath\'eodory loop} of $p$.  

The continuous map $G_p:\C\rightarrow \R$, defined by $G_{p}(z)=\log
|\Psi_{p}^{-1}(z)| $ for $z\in \C-K_p$ and $G_{p}(z)=0$ for $z\in K_p$, is called the Green function of the polynomial $p$. Each level set of the Green function $\{z : G_p(z) = \log(R) \}$ with $R>1$ is called an equipotential for the polynomial $p$. This is the image of the circle of radius $R$ under $\Psi_p$.

\subsection{Horocyclic and radial convergence}\label{sec:horocyclic} 
The topic of convergence of Julia sets in the Hausdorff topology (of compact sets in $\mathbb{P}^1$) is very vast and has been covered by many authors (A. Douady \cite{D}, P. Lavaurs \cite{L}, C. McMullen \cite{Mc}, P. Ha\"issinsky \cite{Hai}, etc.). We only recall here a theorem from \cite{Mc} about horocyclic and radial convergence of rational maps, and give the simplified form of the
theorem for quadratic polynomials.

\begin{defn}[Hausdorff topology] The compact sets $K_n$ converge to
the compact set $K$ in the Hausdorff topology if the following
conditions hold
\begin{itemize}
\item[(a)] Every neighborhood of a point $x\in K$ meets all but
finitely many $K_n$.
\item[(b)] If every neighborhood of $x$ meets infinitely many $K_n$,
then $x\in K$.
\end{itemize}
\end{defn}

\begin{thm}[\cite{Mc}]\label{thm: McMullen} Let $f$ be a geometrically finite rational map
and suppose that $f_n\rightarrow f$ horocyclically (or radially),
preserving critical relations. Then $J_{f_n}\rightarrow J_{f}$ in
the Hausdorff topology.
\end{thm}

Theorem \ref{thm: McMullen}
can be expressed in a simplified form when we restrict to the family
of quadratic polynomials. Let $p$ be a quadratic polynomial. The sequence of polynomials $p_n$ converges to $p$ algebraically if ${\rm deg}(p_n)={\rm deg}(p )$ for all $n$, and the coefficients of $p_n$ converge to the coefficients of $p$.

\begin{defn}[Horocyclic/radial convergence of multipliers]
Let $\lambda_n\rightarrow 1$ in $\C^*$, 
\[
\lambda_n = e^{L_n+i \theta_n}\ \mbox{and}\ \theta_n\rightarrow 0.
\]
The sequence $\lambda_n$  converges to $1$ horocyclically if
$\theta_n^2/L_n\rightarrow 0$. The sequence $\lambda_n$ converges to
$1$ radially if $\theta_n=\bigO(|L_n|)$, that is there exists $M>0$ such
that $|\theta_n|\leq M L_n$ for $n>0$.
\end{defn}

\begin{thm}
Let $p$ be a quadratic polynomial with a parabolic fixed point
$\alpha_0$ with multiplier $e^{2\pi i p/q}$. Let $p_n$ be a
sequence of quadratic polynomials, such that $p_n\rightarrow p$ algebraically.
Assume that each $p_n^q$ has a fixed point $\alpha_n$, such that
$\alpha_n\rightarrow \alpha_0$ and such that the sequence of
multipliers $\lambda_n = (p_n^q)'(\alpha_n)$ converges to $1$
horocyclically (or radially). Then
\[
J_{p_n}\rightarrow J_p
\]
in the Hausdorff topology.
\end{thm}

The proof of Theorem \ref{thm: McMullen} is quite involved and uses
quasiconformal theory and it is not very clear how one could extend it to higher dimensions. We would therefore like to first outline a more topological proof of continuity in one dimension.

\subsection{A topological proof of continuity in dimension one}\label{sec:Cont}

Let $\lambda= e^{2 \pi i p/q}$ be a primitive root
of unity of order $q$. Set
\[
\lambda_{t}=(1+t)\lambda,
\]
for $t$ real and sufficiently small. Consider the  family of quadratic polynomials 
\begin{equation}\label{eq:ptct}
p_{t}(x)=x^{2}+c_{t}, \ \mbox{where}\ c_{t}=\frac{\lambda_{t}}{2} -
\frac{\lambda_{t}^{2}}{4}.
\end{equation}
For $t>0$ the polynomial
$p_{t}$  is hyperbolic and has a repelling fixed
point $\alpha_t=\lambda_t/2$ of multiplier $\lambda_{t}$  and
a $q$-periodic attractive orbit. For $t=0$ the polynomial $p_0$
has a parabolic fixed point $\alpha_0$ of multiplier $\lambda$. The multiplicity of the
fixed point $\alpha_0$ as a solution of the equation $p_0^{\circ q}(x)=x$ is $q+1$. 
Finally, when $t<0$, $p_t$ has an attracting fixed point $\alpha_t$ 
of multiplier $\lambda_{t}$ and a $q$-periodic repelling orbit.

We have  $p_{t} \rightarrow
p_0$ uniformly as $t\rightarrow 0$. The continuity of the
corresponding Julia sets (that we state below as Theorem
\ref{thm:cont-pol}) is an easy consequence of McMullen's Theorem
\ref{thm: McMullen}. The sequence of multipliers
$
\lambda_t^q=(1+t)^q
$
has no imaginary part,
therefore it converges horocyclically and radially to $1$.

\begin{thm}\label{thm:cont-pol}
The Julia set $J_{p_{t}}$ of the polynomial $p_{t}$ depends
continuously on the parameter $t$, that is $J_{p_{t}}\rightarrow J_{p_0}$
in the Hausdorff topology.
\end{thm}

We give a new proof of the continuity result for the
family $p_t$ which does not use quasiconformal theory.
The proof relies on the techniques developed by
Douady and Hubbard in \cite{DH} for proving the local connectivity
of Julia sets of polynomials where all critical points are attracted to attracting or
parabolic cycles. An adaptation of this technique was also
used in \cite{KW} to build semi-conjugacies between Julia sets of
geometrically finite rational maps.

We build a continuous family of bounded
metrics $\mu_{t}$ on the neighborhood of the Julia set $J_{p_t}$,
with respect to which the polynomial $p_t$ is weakly expanding. Then we
will recover the Julia set $J_{p_t}$ as the image of the unique
fixed point $f_t^*$ of a weakly contracting operator $F_t$ in an
appropriate function space. We will show that the fixed point
$f_t^*$ is continuous with respect to $t$, and conclude that the
Julia sets $J_{p_t}$ converge in the Hausdorff topology as
$t\rightarrow 0$ to the Julia set of the parabolic polynomial.

We will illustrate the technique for $t\geq 0$. The case when $t$ is negative is almost identical, but it requires a small technical adjustment, which we will discuss in a more general setting in Chapter \ref{ch:Henon}, where we adapt the construction to a family of polynomial automorphisms of $\C^2$.

\subsection{Normalizing coordinates at a repelling fixed point}\label{sec:Normalizing}

When $t>0$, the polynomial $p_t$ is hyperbolic and expanding with respect to the Poincar\'e
metric on a suitable neighborhood of $J_{p_t}$. 
In order to make the choice
of metrics continuous with respect to the parameter $t$ when $t\rightarrow
0$, we need to correct the metric $\mu_t$ near the repelling fixed
point $\alpha_t$,  which becomes parabolic when $t=0$. 
A naive idea would be to take a small disk $D_t$
around the repelling point $\alpha_t$ on which the polynomial $p_t$ is analytically
conjugated to its linear part $z\rightarrow \lambda_t z$, $|\lambda_t|>1$, hence naturally expanding with respect to the Euclidean metric.
This is not very helpful however, because the radius of $D_t$ converges to $0$ as $t$ converges to $0$. 
This issue can be dealt with by constructing ``normalizing
coordinates" around $\alpha_t$, similar to the parabolic case $t=0$.
We will build a larger neighborhood $\D_{\rho}$, with $\rho$
independent of $t$, around $\alpha_t$, on which the polynomial is
not fully linearized, but rather conjugated to a ``normal form".

Let $\epsilon_{0}=\tan(2\pi/9)$ and  $\epsilon_{1} =
\epsilon_{0}/\sqrt{\epsilon_{0}^{2}+1}$. The meaning of these constants is fully explained by Equation
\eqref{eq:eps1} and the discussion following it.

\begin{prop}\label{prop:normal-polynom-t}
There exist $\delta'>0$ and $\rho>0$ such that for all $t$ with $|t|<\delta'$  there exists a coordinate
transformation  $\phi_t:\D_{\rho'}(\alpha_{t})\rightarrow \D_{\rho}(0)$
defined in a neighborhood of the repelling fixed point $\alpha_t$ such that in the new
coordinates the polynomial $p_{t}$ can be written as $\widetilde{p_{t}}(x)=\lambda_{t} (x+x^{q+1}+C_t x^{2
q+1}+\bigO(x^{2q+2}))$.
Suppose $t\in [0,\delta']$. In the region
\[
\Delta^{-}=\left\{ |x| \leq \rho :
\mbox{Re}(x^{q})>\epsilon_{0}|\mbox{Im}(x^{q})| \right\}
\]
the derivative $\widetilde{p_{t}}'$ is expanding, with a factor of
\[
|\lambda_{t}|\left(1 + (q+3/2)\epsilon_{1}|x|^{q}\right)>|\lambda_t|\geq 1.
\]
The compact region
\[
\Delta^{+}=\left\{ |x| \leq \rho : \mbox{Re}(x^{q})\leq
\epsilon_{0}|\mbox{Im}(x^{q})| \right\}
\]
satisfies
$\Delta^{+}\subset int(K_{\widetilde{p_t}})\cup \{0\}$ and
$\widetilde{p_{t}}(\Delta^{+})\subset int(\Delta^{+})\cup \{0\}$.
\end{prop}
\proof One performs for the family $p_t$ the same sequence of
coordinate transformations as the ones done in \cite{BH} or \cite{DH} in the
parabolic case. After a global coordinate change that brings the
fixed point $\alpha_t$ to the origin, we can assume that
$p_t(x)=\lambda_t x +x^2$. Suppose by induction that for $k\geq 2$ the maps $p_t$
have the form
\[
x_{1}= \lambda_t x+a_{t}x^{k}+ \bigO(x^{k+1})
\]
where $a_t\neq 0$ for $|t|<\delta'$.

Consider the coordinate transformation
\begin{equation*}
\begin{array}{l}
    X=x+b_t x^{k}
\end{array} \ \ \ \mbox{with inverse}\ \ \
\begin{array}{l}
    x=X-b_t X^{k}+\ldots
\end{array}
\end{equation*}
In the new coordinate system, we get
\begin{eqnarray*}
X_{1}&=&x_{1}+b_t x_{1}^{k} = \lambda_t x +(a_{t}+b_t\lambda_t^{k})x^{k} + \ldots \\
    &=&\lambda_t ( X -  b_t X^{k} + \ldots) +(a_{t}+b_t\lambda_t^{k})(X-b_{t}X^{k} + \ldots)^{k} + \ldots\\
    &=&\lambda_t X  +  (a_{t}+b_t(\lambda_t^{k}-\lambda_t)) X^{k}+ \ldots
\end{eqnarray*}
When $t\neq 0$, we have $|\lambda_t|\neq1$, so $\lambda_t\neq
\lambda_t^{k}$ for all $k$ with $2\leq k \leq q$. If $k$ is not congruent to $1$ modulo $q$, then 
$\lambda_{0}^{k}\neq \lambda_0$ as well, so we can set
\[
b_t = \frac{a_{t}}{\lambda_t-\lambda_t^{k}}
\]
and eliminate the term $a_{t}x^{k}$.  The
transformations $x+b_tx^{k}$ are injective on a uniform neighborhood
of $0$.

This proves that by successive coordinate transformations of the
form $X_t=x+b_{t}x^{k}$
we can
eliminate terms with powers that are not congruent to $1$ modulo
$q$, so the first term that cannot be eliminated in this way will
have a power of the form $a_{t}x^{\nu q + 1}$, for some integer
$\nu\geq 1$. The parabolic multiplicity of the fixed point
$\alpha_0$ is $1$, so $\nu=1$.

Thus the map takes the form
\begin{equation}\label{eq:pol3a}
\begin{array}{l}
    x_{1}= \lambda_t (x+a_{t}x^{q+1} + \bigO(x^{q+2}) )
\end{array}
\end{equation}
Of course, when $t\neq 0$, we could use the same 
map $X_t=x+b_{t}x^{q+1}$ to eliminate the term $a_tx^{q + 1}$,
however this would require shrinking the domain of injectivity
of $X_t$ to $0$ as $t\rightarrow 0$, as well as losing the
continuity of $X_t$ with respect to $t$ at $t=0$.

We can further reduce Equation \eqref{eq:pol3a} to 
\begin{equation}\label{eq:pol3b}
    x_{1}= \lambda_t (x+x^{q+1} + \bigO(x^{q+2})).
\end{equation}
by considering a linear map $X=A_tx$, where $A_t$ is a
constant chosen such that $A_t^{ q}=a_{t}$.

In Equation \eqref{eq:pol3b} we can eliminate all terms of the form
$a_t x^k$, with $q+1<k<2q+1$, using the same coordinate
transformation as before, $X_t=x+b_{t}x^{k}$, where
$b_t=\frac{a_{t}}{\lambda_t-\lambda_t^{k}}$. We
thus arrive at the normal form
\begin{equation*}
\widetilde{p_{t}}(x)=\lambda_{t} (x+x^{q+1}+C_t x^{2
q+1}+\bigO(x^{2q+2})).
\end{equation*}

When $t=0$ the regions $\Delta^+$ and $\Delta^-$ represent
attracting and respectively repelling sectors for the (normalized)
parabolic polynomial $\widetilde{p}_0$. The attractive sector
$\Delta^+$ belongs to the interior of the filled-in Julia set
$\mbox{int}(K_{\widetilde{p}_0})\cup \{0\}$ and all points in
$\Delta^+$ converge under forward iterations to the parabolic fixed
point $0$, which lies in the Julia set. When $t<0$, the sector $\Delta^+$ belongs to the basin of attraction of the attracting fixed point $0$. When $t>0$, the sector $\Delta^+$ belongs to
$\mbox{int}(K_{\widetilde{p}_t})\cup \{0\}$, because
$\Delta^+$ (with $0$ removed) is a trapping region for a
$q$-periodic attractive orbit \cite{BH}. We prove these facts directly for \He maps in
Propositions \ref{prop:sectors-t} and \ref{prop:trapping}, and the proofs apply also to the family of polynomials considered in this lemma.

To show that the derivative $\widetilde{p_t}'$ is expanding on
$\Delta^-$ when $t>0$, we perform the same computations as in the parabolic
case \cite{DH}. The choice of $\epsilon_{0}$ and $\epsilon_{1}$ guarantee that if $Re(x^{q})>\epsilon_{0}Im(x^{q})$, then $Re(x^{q})>\epsilon_{1}|x|^{q}$. Consider a constant $m$ so that $\big{|}\widetilde{p_t}'(x)-\lambda_t(1+(q+1)x^{q})\big{|}<m|x|^{2q}$ on $\D_{\rho}$. Using the normal form for $\widetilde{p_t}$, we get
\begin{eqnarray*}
|\widetilde{p_t} '(x)|&=&|\lambda_t| \big{|}1+(q +1)x^{ q}+\bigO(x^{2 q})\big{|}\geq |\lambda_t| \left(\big{|}1+(q +1)x^{ q}\big{|}-m|x|^{2q}\right)\\
&\geq& (1+t) \left(1 + (q+1)\epsilon_{1}|x|^{q} -
m|x|^{2q}\right)
> (1+t)\left(1 + (q+3/2)\epsilon_{1}|x|^{q}\right).
\end{eqnarray*}
for $|x|$ sufficiently small. Hence $|\widetilde{p_t}'(x)|>|\lambda_{t}|$ throughout $\Delta^{-}$, for all $t\in [0,\delta']$.
\qed

Let $\Delta^+_t=\phi_t^{-1}(\Delta^+)$ and
$\Delta^-_t=\phi_t^{-1}(\Delta^-)$. Recall that $\alpha_t=\phi_t^{-1}(0)$. 
By Proposition \ref{prop:normal-polynom-t}, the set $\Delta^+_t-\{\alpha_t\}$ belongs to 
the interior of the filled-in Julia set $K_{p_t}$. Moreover, when $t>0$,
the $q$-periodic attractive orbit of the polynomial $p_t$ is
contained in the sector $\Delta^+_t$. The Julia set $J_{p_t}$ near $\alpha_t$ is completely
contained in the repelling sectors:
\[
J_{p_t}\cap \D_{\rho'}(\alpha_{t}) = J_{p_t} \cap \Delta^-_t.
\]

When $t\in [-\delta',\delta']$, the polynomial $p_t$ has connected Julia
set; the critical point $0$ of $p_t$ is attracted to the
$q$-periodic orbit when $t>0$, respectively to the parabolic fixed
point when $t=0$, and to the attracting fixed point when $t<0$. So there exists a first iterate $n_t\in \N$ such
that $p_t^{\circ (n_t+1)}(0)\in \Delta^+_t$, otherwise said, there
exists a first iterate for which $p_t^{\circ n_t}(c_t)\in
\Delta^+_t$ and $p_t^{\circ n_t}(0)\notin \Delta^+_t$. The function
$t\rightarrow n_t$ is locally constant and we can assume without
loss of generality that when $\delta'$ is small, 
the number $n_t$ is the same for
all $t\in [-\delta',\delta']$. Therefore we can remove the dependence on $t$ and denote $n_t$ by $N$. 
Denote further by $p_t^{-\circ N}(\Delta^+_t)$ the
connected component of the $\small{N}^{th}$ preimage of the set $\Delta^+_t$
that contains the fixed point $\alpha_t$.

\subsection{A continuous family of bounded metrics}\label{sec:metrics1D}

For each value of the parameter $t$ we construct a
neighborhood $U_t$ of the Julia set $J_{p_t}$ and a metric $\mu_t$
on $U_t$ with respect to which the polynomial $p_t$ is expanding.
The family $(U_t, \mu_t)$ will be continuous with respect to the
parameter $t$. 

\begin{figure}[htb]
\begin{center}
\includegraphics[scale=1.1]{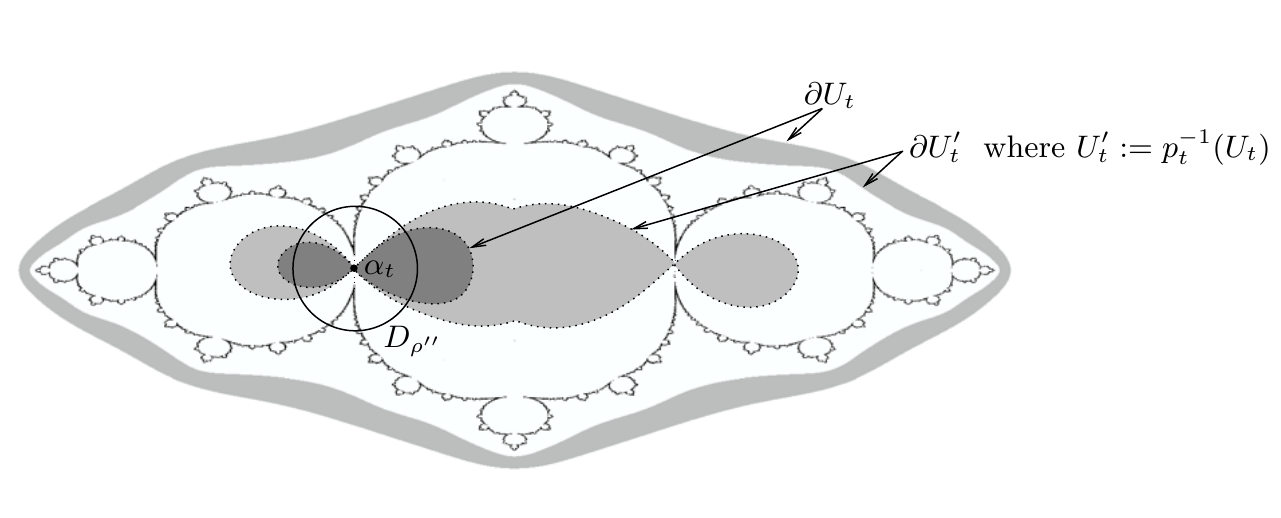}
\end{center}
\caption{The polynomial $p_{t}$ has a fixed point at $\alpha_{t}$. The corresponding neighborhoods $U_{t}$ and $U_{t}'$ are also shown.}
\label{fig:Ut}
\end{figure}

\noindent The outer boundary of the set $U_t$ is an
equipotential of the Julia set $J_{p_t}$. The inner boundary is
$\partial p_t^{-\circ N}(\Delta^+_t)$, where $N$ is defined above. Formally, choose $R>2$ and set
\begin{eqnarray}\label{def: U_t}
U_t = \C-\overline{p_t^{- \circ N}(\Delta^+_t)}-\{z\in \C-K_{p_t} : |\Psi_{p_{t}}^{-1}(z)|\geq R\}.
\end{eqnarray}
Let $U'_t=p_t^{-1}(U_t)$. The set $U'_t$ is contained in $U_t$ by
construction, and we can put on $U'_t$ the Poincar\'e metric of
$U_t$. The map $p_t:U'_t\rightarrow U_t$ is a covering map, hence
expanding:
\[
|(z,\xi)|_{U_t}<|(z,\xi)|_{U_t'}=|(p(z),p'(z)\xi)|_{U_t} \ \
\mbox{for}\ z\in U'_t\ \mbox{and}\ \xi\in T_zU'_t.
\]
However $U'_t$ is not relatively compact in $U_t$ because $\partial
U'_t \cap
\partial U_t = \{\alpha_t\}$, so there is no constant of uniform
expansion. On the repelling sectors $\Delta^-_t$, one can define a
metric $\mu_{\Delta^-_t}$ as the pull-back of the Euclidean metric
from the normalizing coordinates $\Delta^-$.
\[
|(z,\xi)|_{\Delta^-_t}:=|\phi_t'(z)\xi| \ \ \mbox{for}\ z\in \Delta^-_t\
\mbox{and}\ \xi\in T_z\Delta^-_t,
\]
where the latter length is the modulus of the complex number
$\phi_t'(z)\xi$.

\begin{defn}\label{defn: metric-inf-poly}
Let $\mu_t=\inf (\mu_{U_t}, M \mu_{\Delta_t^-})$, where $M$ is a
positive real number and
\[
M>\sup\left\{\frac{2\mu_{U_t}(z,\xi)}{\mu_{\Delta_t^-}(p_t(z),p_t'(z)\xi)}\
:\ z\in p_t^{-1}(\Delta_t^-), z\notin \Delta_t^- \mbox{ and }
t\in[0,\delta'] \right\}.
\]
\end{defn}

By choosing $M$ sufficiently large, one can assure that on the
boundary of $V_t$ the infimum is attained by the Poincar\'e metric
$\mu_{U_t}$. So the metric $\mu_t$ is continuous on $U_t'$. Note also
that the Poincar\'e metric is infinite at $\alpha_t$ while the Euclidean metric is bounded; therefore there exists a neighborhood of 
$\alpha_t$, uniform with respect to $t$, for which the infimum in
Definition \ref{defn: metric-inf-poly} is attained by the Euclidean
metric $M \mu_{\Delta_t^-}$.

\begin{lemma}\label{lemma: Poincare-dominated-by-Euclidian}
The family of metrics $\mu_t$ depends continuously on the parameter
$t$ and it is dominated above and below by the Euclidean metric.
There exist $m_1>0$ and $m_2>0$ such that
\[
m_1 |x-y|<d_{\mu_t}(x,y)< m_2 |x-y|, \mbox{ for any } x,y\in U_t'.
\]
\end{lemma}
\proof By construction, the neighborhood $U_t$ and the repelling
sectors $\Delta_t^{-}$ depend continuously on $t$. Let $\rho_t$
denote the density function of the Poincar\'e metric on $U_t$,
$\mu_{U_t}(z,dz)=\rho_t(z)|dz|$.  The map $\rho_t$ is
positive, $C^{\infty}$-smooth on $U_t'$ and continuous with respect
to $t$. Hence $\mu_{U_t}$ is bounded below by the
Euclidean metric on $U_t'$.  The metric $\mu_{U_t}$
on this set is also bounded above on $U_t'$, except on a small neighborhood
of the fixed point $\alpha_t\in\partial U_t$.

The metric $\ds\mu_{\Delta^-_t}$ is the pull-back of the Euclidean
metric by a holomorphic injective map $\phi_t$, continuous with
respect to $t$. We have $\mu_{\Delta^-_t}(z,
dz)=|\phi_t'(z)||dz|$, where $|\phi_t'(z)|>0$ is bounded above and
below on $\Delta^-_t$.
Therefore the infimum metric $\mu_t$ is bounded above and below with
respect to the Euclidean metric. \qed

\begin{lemma}\label{lemma: vector-1dim} The polynomial $p_t$ is
strictly expanding with respect to the metric $\mu_t$ on the set
$U_t'$ when $t>0$.
\end{lemma}
\proof Let $z,z'\in U_t'$ and $\xi \in T_z U_t',\  \xi'\in T_{z'}
U_t'$ such that $z'=p_t(z)$ and $\xi'=p_t'(z)\xi$.

We will show that for each $t>0$ there exists a constant $k_t>1$
such that
\[
\mu_t(z',\xi')> k_t\cdot \mu_t(z, \xi).
\]
There are four cases to consider:
\begin{enumerate}
  \item [(a)] $\mu_t(z,\xi)=\mu_{U_t}(z,\xi)$ and
  $\mu_t(z',\xi')=\mu_{U_t}(z',\xi')$.

  \noindent This happens only if both $z$
  and $z'$ are outside a small neighborhood of the
  point $\alpha_t$. Outside this neighborhood, the set $U_t'$ is compactly
  contained in $U_t$, so $p_t$ expands strongly with respect
  to the Poincar\'e metric. For all $t\in [0,\delta']$ there exists $\kappa_t>1$  such that
  \begin{equation}\label{eq: c_t-polynom}
  \mu_{U_t}(z',\xi')>\kappa_t\cdot \mu_{U_t}(z,\xi).
  \end{equation}
  The constant $\kappa_t$ depends only on the distance between the boundaries $\partial
  U_t$ and $\partial U_t'$ outside a disk of fixed size around the fixed point $\alpha_t$, so 
  $\inf\limits_{t\in [0,\delta']} \kappa_t>1$.

  \item [(b)] $\mu_t(z,\xi)=M\mu_{\Delta_t^-}(z,\xi)$ and
  $ \mu_t(z',\xi')=M\mu_{\Delta_t^-}(z',\xi')$.

  \noindent The normalized polynomial $\widetilde{p}_t$ expands with respect
  to the Euclidean metric, so by Proposition \ref{prop:normal-polynom-t} we have
  \[
  \mu_{\Delta_t^-}(z',\xi')>(1+t)\cdot\mu_{\Delta_t^-}(z,\xi).
  \]
Notice that the constant of expansion is $1$ if and only if
$t\rightarrow 0$ and $\phi_t(z)\rightarrow 0$ (that is, $z$ approaches the parabolic fixed point
$\alpha_0$).

  \item [(c)] $\mu_t(z,\xi)=M\mu_{\Delta_t^-}(z,\xi)$ and $\mu_t(z',\xi')=\mu_{U_t}(z',\xi')$.

\noindent Similar to case (a), the point $z'$ cannot be too close to the fixed
point $\alpha_t$, so
\[
\mu_{U_t}(z',\xi')>\kappa_t\cdot\mu_{U_t}(z,\xi)\geq \kappa_t\cdot
M\mu_{\Delta_t^-}(z,\xi).
\]

  \item [(d)] $\mu_t(z,\xi)=\mu_{U_t}(z,\xi)$ and $
  \mu_t(z',\xi')=M\mu_{\Delta_t^-}(z',\xi')$.

\noindent There are two sub-cases to consider
  \begin{enumerate}
    \item [(i)] If $z\in p_t^{-1}(\Delta_t^-)\cap \Delta_t^-$, then
    \begin{eqnarray*}
M\mu_{\Delta_t^-}(z',\xi')>(1+t)\cdot
M\mu_{\Delta_t^-}(z,\xi)>(1+t)\cdot\mu_{U_t}(z,\xi).
    \end{eqnarray*}
    \item [(ii)] If $z\in p_t^{-1}(\Delta_t^-)$ but $z\notin \Delta_t^-$, then
    the conclusion follows from the choice of the constant $M$, as shown below:
    \[
2\mu_{U_t}(z,\xi)=\frac{2\mu_{U_t}(z,\xi)}{\mu_{\Delta_t^-}(z',\xi')}\cdot
\mu_{\Delta_t^-}(z',\xi')<M\mu_{\Delta_t^-}(z',\xi').
    \]
  \end{enumerate}
\end{enumerate}
Set $k_t:=\min\left(1+t,\kappa_t\right)$. From estimates (a),(b),(c) and (d) we
can easily see that
\[ \mu_t(z',\xi')> k_t\cdot \mu_t(z, \xi).
\]
We get uniform expansion when $t>0$ because $k_t$
is strictly greater than $1$. \qed

The metric $\mu_t$ induces a natural path metric on $U_t'$. If
$\eta:[0,1]\rightarrow U_t'$ is a rectifiable path, then its
length with respect to the metric $\mu_t$ is given by the formula
\begin{equation}\label{eq:l-pol}
\ell_{\mu_t}(\eta)=\int_0^1 \mu_t \left(\eta(s),\eta'(s)\right) ds.
\end{equation}
The distance between two points $x$ and $y$ from $U_t'$ with
respect to the metric $\mu_t$ is
\begin{equation}\label{eq:d-pol}
d_{\mu_t}(x,y)=\inf \ell_{\mu_t}(\eta),
\end{equation}
where the infimum is taken after all rectifiable paths
$\eta:[0,1]\rightarrow U_t'$ with $\eta(0)=x$ and
$\eta(1)=y$.

\subsection{Contraction in the space of functions}\label{sec:contraction}

For each value of the parameter $t$, we will construct a sequence of equipotentials in the complement of the filled Julia set $K_{p_t}$ and show that they converge to the Julia set $J_{p_{t}}$, uniformly with respect to $t$.
In our setting, the filled Julia set $K_{p_t}$ is connected. Let $\Psi_{p_t}: \C-\overline{\D} \rightarrow \C-K_{p_t}$ be the inverse B\"ottcher isomorphism of the
polynomial $p_t$ as in \eqref{eq:psi-p} and let $\gamma_{t}:\s^{1}\rightarrow J_{p_{t}}$ be the
Carath\'eodory loop of $p_{t}$ as in \eqref{eq:gamma-p} (i.e. the continuous extension of $\Psi_{p_t}$ to the boundary). We write $\gamma_{t}$ instead of $\gamma_{p_{t}}$ to simplify notations. By the definition of the isomorphism $\Psi_{p_t}$ we have $\Psi_{p_t}(z^2)=p_t(\Psi_{p_t}(z))$, for $z\in   \C-\overline{\D}$.

Let $R>2$ be a fixed constant, chosen as in Equation \eqref{def: U_t}. For each $t\in [0,\delta']$, consider the space
of functions
\begin{eqnarray*}
\mathcal{F}_{t} &=& \left\{\gamma_{t,n}:\s^1\rightarrow U_t'\ : \
\gamma_{t,n}(s)=\Psi_{p_t}\left(R^{1/2^{n+1}} e^{2\pi i s}\right),\ n\in \N\right\}.
\end{eqnarray*}

For each $t$, the space $\mathcal{F}_t$ is just a sequence of parametrized equipotentials $\{\gamma_{t,n}\}_{n\geq 0}$, corresponding to the polynomial $p_t$.  The Green function $G_{p_t}$ of the polynomial $p_t$  is continuous with respect to $t$ and $z$.
Therefore each map $(t,s)\mapsto\gamma_{t,n}(s)$ is continuous with respect
to $t$ and $s$.
The polynomial $p_t$ maps each equipotential $\gamma_{t,n}$ to the equipotential $\gamma_{t,n-1}$ by a two-to-one covering map.  We can select a branch of $p_t^{-1}$ by using the inverse B\"ottcher isomorphism and setting
\[
p_t^{-1}\left( \Psi_{p_t} \left(R^{1/2^{n}} e^{2\pi
i\, (2s)} \right)\right) :=\Psi_{p_t} \left(R^{1/2^{n+1}} e^{2\pi
i s}\right)\ \mbox{for } s\in \s^1 \mbox{ and }  n\geq 1 .
\]

\noindent Therefore, the space $\mathcal{F}_{t}$ comes with a natural operator
$p_t^{-1}:\mathcal{F}_{t}\rightarrow \mathcal{F}_{t}$, given by the
rule
\begin{equation}\label{def: p-1}
p_t^{-1}(\gamma_{t,n-1}(2s))=\gamma_{t,n}(s), \ s\in \s^1,\ n\geq1.
\end{equation}
Endow the function space $\mathcal{F}_t$ with the supremum
metric
\[
d_{\mu_t}(\gamma_{t,n}, \gamma_{t,k})=\sup\limits_{s\in
\s^1}d_{\mu_t}\left( \gamma_{t,n}(s), \gamma_{t,k}(s)\right)
\]
and let $\mathcal{\overline{F}}_{t}$ be the completion of
$\mathcal{F}_{t}$ with respect to the supremum metric $d_{\mu_t}$. Notice also that the metric $d_{\mu_t}$ is bounded, by Lemma \ref{lemma: Poincare-dominated-by-Euclidian}.

\begin{thm}[Browder \cite{Br},\cite{KS}]\label{thm: Browder} Let $(X,d)$ be a complete metric space and suppose $f:X\rightarrow X$ satisfies
\[
	d(f(x),f(y))<h(d(x,y))\ \ \ \mbox{for all}\ \ x,y\in X,
\]
where $h:[0,\infty)\rightarrow [0,\infty)$ is increasing, continuous from the right, and  $h(s)<s$ for all $s>0$.  Then $f$ has a unique fixed point $x^{*}$ and $f^{n}(x)\rightarrow x^{*}$ for each $x\in X$.
\end{thm}
\begin{defn}
We will call a function $h$ that verifies the hypothesis of Theorem
\ref{thm: Browder} a \emph{Browder function}.
\end{defn}

\begin{remarka}\label{remark: Browder} Assume that the space $X$ from Theorem \ref{thm: Browder} is bounded.
The rate of convergence to the fixed point is controlled by the
function $h$, namely if we choose $L>0$ such that
$L-h(L)>\mbox{diam}(X)$, then the following estimate holds
\[
d(f^{\circ n}(x), x^*)<h^{\circ n}(L),\ \ \mbox{for any}\ x\in X,\
n\in \N.
\]
\end{remarka}
We apply Browder's Theorem to the complete metric space $(\mathcal{\overline{F}}_{t}, d_{\mu_t})$ with the operator $p_t^{-1}$, to show that  for each $t\geq 0$, the sequence $\{\gamma_{t,n}\}_{n\geq 0}$ converges uniformly as $n\rightarrow \infty$ to a continuous function $\gamma_t$. As a consequence of the same Theorem \ref{thm: Browder} and Remark \ref{remark: Browder} , we obtain the continuity of the map $t\rightarrow \gamma_t$ with respect to the parameter $t$. 

For each $t \in [0,\delta']$, consider the function $h_t:  [0,\infty)\rightarrow [0,\infty)$ given by
\begin{equation*}
h_t(s) := \sup\left\{d_{\mu_t}(x, y)\ : x,y\in p_{t}^{-1}(U_t')\
\mbox{and}\  d_{\mu_t}(p_t(x),p_t(y))\leq s\right\}.
\end{equation*}

Clearly $h_t$ is increasing (i.e. $s_1<s_2\Rightarrow h_t(s_1) \leq h_t(s_2)$) and satisfies the inequality
\[
d_{\mu_t}(\gamma_{t,n+1},\gamma_{t,k+1})<h_t(d_{\mu_t}(\gamma_{t,n},\gamma_{t,k})),\
k,n\in\N.
\]
 The family of metrics $\mu_t$ is continuous with respect
to $t$, so the map $t\rightarrow h_t$ is continuous with respect to
$t\in[0,\delta']$. Moreover, by Lemma \ref{lemma:
vector-1dim}, when $t>0$ we have 
\begin{equation}\label{eq:ht}
h_t(s)<\frac{s}{k_t}\ \ \ \mbox{for}\ s>0.
\end{equation}

Inequality \eqref{eq:ht} implies that $p_t^{-1}$ is
strictly contracting with respect to the metric $d_{\mu_t}$ when $t>0$. Banach
Fixed Point Theorem assures that for each $t>0$, the operator $p_t^{-1}$ has a unique fixed point
$\gamma_{t}:\s^1\rightarrow \overline{U_t'}$, and the sequence
$\gamma_{t,n}$ converges to $\gamma_{t}$ as $n\rightarrow \infty$.

In \cite{DH} (and also in \cite{H}) it is shown that the function $h_0$
verifies the
hypothesis of Theorem \ref{thm: Browder}, hence the sequence
$\gamma_{0,n}$ converges to the unique fixed point of the operator
$p_0^{-1}:\mathcal{F}_{0}\rightarrow \mathcal{F}_{0}$, that is to a
continuous function $\gamma_{0}:\s^1\rightarrow \overline{U_0'}$.
The image of $\gamma_{0}$ is invariant under the parabolic polynomial $p_0$ and it
parametrizes its Julia set $J_{p_0}$.

Notice that the constant $k_t$ goes to $1$ when $t$ goes to $0$, so we haven't obtained any information yet about the continuity of the map $t\rightarrow \gamma_t$ with respect to $t$ when $t=0$.

To provide a unified approach to the hyperbolic and parabolic cases,
we define a new map $h:[0,\infty)\rightarrow [0,\infty)$, $h(s):= \sup\limits_{t\in
[0,\delta']} h_t(s)$.

\begin{lemma}\label{lemma: properties-h} The map $h$ is
increasing and $h(s)<s$ for all $s>0$. 
\end{lemma}
\proof When $t=0$, it is proven in \cite{DH} and \cite{H} that $h_0(s)<s$ for all
$s>0$. When $t>0$, Inequality \eqref{eq:ht} yields that $h_t(s)<s$ for all $s>0$. 
For a fixed $s\in \R^+$, the map $t\rightarrow h_t(s)$ is continuous
with respect to $t$, thus it attains its supremum on $[0,\delta']$,
so $h(s)<s$. For each $t$, the
function $h_t$ is increasing, by definition.  
It is obvious then that the function $h$ is also increasing.
\qed

The function $h$ is increasing, so $h(s+)=\lim_{\epsilon \rightarrow 0^{+}}h(s+\epsilon)$ is well defined. 
The function $h^+:s\mapsto h(s+)$ is right
continuous and the following lemma holds:

\begin{lemma}\label{lemma: Browder-function-h}
The function $h^+:[0,\infty)\rightarrow[0,\infty)$ is a Browder function, i.e. it is right
continuous, increasing, and $h^+(s)<s$ for all $s>0$.  
Moreover
\begin{equation*}
d_{\mu_t}(\gamma_{t,n+1},\gamma_{t,k+1})<h^+(d_{\mu_t}(\gamma_{t,n},\gamma_{t,k}))\
\ \ \mbox{for all}\ t\in[0,\delta']\ \mbox{and}\ k,n\in\N.
\end{equation*}
\end{lemma}
\proof The only non-trivial property to check is the fact that
$h(s+)<s$ for all $s>0$. By Lemma \ref{lemma: properties-h} we know
that $h(s)<s$ for all $s>0$, so $h(s+)\leq s$ for all $s>0$.

Suppose that $h(s+)=s$ for some $s>0$. Let $\epsilon_n\searrow 0$ be
a decreasing sequence of positive numbers such that
$h(s+)=\lim_{\epsilon_n\rightarrow 0}h(s+\epsilon_n)$. From the
definition of  $h$ we have
\[
h(s+\epsilon_n)=\sup\limits_{t\in [0,\delta']}h_t(s+\epsilon_n).
\]
From the definition of the supremum, for every $n>0$ there exists
$t_n\in [0,\delta']$ such that
\begin{equation}\label{eq:h}
h_{t_n}(s+\epsilon_n)>\sup_{t\in
[0,\delta']}h_t(s+\epsilon_n)-\epsilon_n=h(s+\epsilon_n)-\epsilon_n
\end{equation}
The function $h_{t_n}$ satisfies
\begin{equation}\label{eq:h2}
h_{t_n}(s+\epsilon_n)<(s+\epsilon_n)\cdot
\frac{1}{k_{t_n}}<s+\epsilon_n.
\end{equation}
The sequence $t_n$ is bounded, so after passing to a convergent
subsequence, we may assume that $t_n\rightarrow \tau$ for some
$\tau\in [0,\delta']$. Let us show that $\tau=0$. Assume that
$\tau\neq 0$. From inequalities \eqref{eq:h} and \eqref{eq:h2} we know that
\begin{equation*}
h(s+\epsilon_n)-\epsilon_n<h_{t_n}(s+\epsilon_n)<\frac{1}{k_{t_n}}\cdot(s+\epsilon_n).
\end{equation*}
Taking the limit as $n\rightarrow\infty$, we get
\[
s=h(s+)\leq \frac{1}{k_{\tau}}\cdot s.
\]
Since $s>0$ we get $k_{\tau}\leq 1$. Then $\tau=0$, otherwise we would have $k_{\tau}>1$.

Pick $m$ an integer, and let $n$ be any integer $n\geq m$.
The sequence $\epsilon_n$ is decreasing, so $\epsilon_n
<\epsilon_m$. The function $h_{t_n}$ is increasing, so
\begin{equation}\label{eq:h3}
h_{t_n}(s+\epsilon_n)\leq h_{t_n}(s+\epsilon_m),\ \ \mbox{for any}\
n\geq m.
\end{equation}
From the inequalities \eqref{eq:h} and \eqref{eq:h3} we obtain
\[
h(s+\epsilon_n)-\epsilon_n< h_{t_n}(s+\epsilon_m),\ \ \mbox{for
any}\ n\geq m.
\]
After passing to the limit as $n\rightarrow \infty$ and using the
continuity of $h_t$ with respect to $t$, we get
\[
s=h(s+)\leq h_0(s+\epsilon_{m}),\ \mbox{for every}\ m\in \N.
\]
Letting $m\rightarrow\infty$ we get $s\leq h_{0}(s+)$. 
In the parabolic case $h_0(s+)<s$ for every $s>0$. This yields $s=0$, which is a contradiction. 
\qed

\begin{remarka} Lemma \ref{lemma: Browder-function-h} is very
important, because it provides a {\it reduction of the hyperbolic case to the
parabolic case}, hence allowing a uniform treatment of both cases.
Lemma \ref{lemma: Browder-function-h} uses only  minimum information about the parabolic case, that is, the fact that $h_0(s+)<s$ for all $s>0$. Another remark is that in the one-dimensional setting, one could presumably prove that the maps $h_t$ are already right continuous, so $h_t(s)=h_t(s+)$. However, we will apply this lemma in higher dimensions (where the maps $h_t$ will not necessarily be  right continuous), so the existence of the Browder function $h^+$ from Lemma \ref{lemma: Browder-function-h} bypasses this problem and is central for the application of Browder's Fixed Point Theorem.
\end{remarka}

\begin{thm}\label{thm: UniformConvergence}
For each $t$, the sequence $\gamma_{t,n}$ converges to a fixed point
$\gamma_{t}:\s^1\rightarrow \overline{U_t'}$ of the operator
$p_t^{-1}$. The rate of convergence to the fixed point is uniform
in $t$.
\end{thm}
\proof By Lemma \ref{lemma: Browder-function-h}, we can use the same Browder function $h^+$ for any parameter $t\in[0,\delta']$ and apply Theorem \ref{thm: Browder} to show the existence of a unique fixed point  $\gamma_t$ for the operator $p_t^{-1}$. Each map $\gamma_t:\s^1\rightarrow \overline{U_t'}$ is continuous.

We show that $\gamma_{t,n}$ converges to $\gamma_{t}$ uniformly
with respect to $t\in [0,\delta']$. 
By Lemmas \ref{lemma: Browder-function-h} and \ref{lemma: Poincare-dominated-by-Euclidian} and Remark
\ref{remark: Browder} there exists $L>0$ such that
\begin{equation*}\label{ineq: unif}
m_1\|\gamma_{t,n}- \gamma_{t}\|<d_{\mu_t}\left( \gamma_{t,n},
\gamma_{t} \right)\leq (h^+)^{\circ n}(L)\searrow 0, \mbox{ for any }
t\in [0,\delta'],
\end{equation*}
so the rate of convergence of $\gamma_{t,n}$ to the fixed point
$\gamma_{t}$ is bounded by the rate at which the sequence
$(h^+)^{\circ n}(L)$ decreases to $0$. Therefore the sequence of
functions $t\rightarrow\gamma_{t,n}$ converges as $n\rightarrow
\infty$ to $\gamma_{t}$, uniformly with respect to $t$.
\qed

In Theorem \ref{thm: UniformConvergence}, we have
constructed a sequence of functions $\gamma_{t,n}$, continuous with
respect to $t$, and proved that it converges uniformly as
$n\rightarrow\infty$ to $\gamma_{t}$. Hence the limit function
$\gamma_{t}$ is continuous with respect to $t$ on $[0,\delta]$.
The continuity of the Julia sets $J_{p_t}$ from Theorem \ref{thm:cont-pol} follows
immediately as $J_{p_t}=\mbox{Im}(\gamma_{t})$.

\section{Continuity of Julia sets for H\'{e}non maps}\label{ch:Henon}

\noindent Inspired by the one-dimensional setting from Chapter \ref{ch:Continuity1}, we now turn back to
dynamics in two complex variables and prove  a continuity result in Theorem \ref{thm:continuity} for
\He maps with a semi-parabolic fixed point. We will consider the \He
map written in the form
\[
H_{c,a}\hvec{x}{y}=\hvec{p(x)+ay}{ax}, \ \ \mbox{where} \ \ \ p(x)=x^2+c.
\]
When $a\neq 0$, this map is a biholomorphism of
constant Jacobian $-a^2$, whose inverse is
\[
 H_{c,a}^{-1} \hvec{x}{y} = \hvec{y/a}{(x-p(y/a))/a}.
\]

As in \cite{HOV1}, for $r>0$ large enough, the dynamical space $\C^2$ can be divided
into three regions: the bidisk $ \D_{r}\times\D_{r}=\{(x,y): |x|\leq r, |y|\leq r\}$,
\begin{equation}\label{eq:V+}
V^+=\{(x,y) : |x|\geq \max(|y|,r) \}\ \
\mbox{and}\ \ V^-=\{(x,y) : |y|\geq \max(|x|, r)\}.
\end{equation}
The escaping sets $U^{\pm}$ can be
described in terms of $V^{\pm}$ as follows: 
$U^+=\bigcup_{k\geq 0} H^{-\circ k}(V^+)$ and $U^-=\bigcup_{k\geq 0} H^{\circ k}(V^-)$.
By taking their complements in $\C^{2}$ we obtain $K^{+}=\C^{2}-U^{+}$, the set of points that do not escape to infinity in forward time, and 
$K^{-}=\C^{2}-U^{-}$,  the set of points that do not escape to infinity in backward time. The Julia set $J^{+}$ is the common boundary of $K^{+}$ and $U^{+}$. Similarly $J^{-}$ is the common boundary of $K^{-}$ and $U^{-}$. In fact, in the dissipative case, $J^{-}=K^{-}$ (see \cite{FM}). The sets $J=J^{+}\cap J^{-}$ and $K=K^{+}\cap K^{-}$ are contained in $\D_{r}\times\D_{r}$.

\begin{defn}\label{def:fixedpoint} Let $\fq$ be a fixed point of $H$ and $\lambda$ and $\nu$ be the two eigenvalues of $DH_{\fq}$. The fixed point $\fq$ is called:
\begin{itemize}
\item[(a)] {\it hyperbolic} if 
$|\nu|<1$ and $|\lambda|>1$;
\item[(b)] {\it semi-parabolic} if $|\nu|<1$ and 
 $\lambda=e^{2\pi i p/q}$;
 \item[(c)] {\it attracting} if 
$|\nu|<1$ and $|\lambda|<1$.
\end{itemize}
\end{defn}

Let $\lambda= e^{2 \pi i p/q}$ be a root of unity of order $q$ and set $\lambda_{t}:=(1+t)\lambda$.
In Section \ref{sec:Cont} we considered the family $p_t$ of polynomials, $p_{t}(x)=x^{2}+c_{t}$, with a fixed point $\alpha_t=\lambda_t/2$ of multiplier $\lambda_{t}$. The exact formula for the coefficient $c_{t}$ is given by Proposition \ref{prop:Plambda} below. 
In Theorem \ref{thm:cont-pol} we showed that the Julia sets $J_{p_t}$ converge to the Julia set $J_{p_0}$ as $t\rightarrow0$. 

In 1-D, the multiplier of a fixed point of a quadratic polynomial
uniquely identifies the polynomial. The following proposition
provides a description of the parameter space of \He maps for which
one eigenvalue of the fixed point is known.

\begin{prop}\label{prop:Plambda}
The set $\mathcal{P}_{\lambda_{t}}$ of parameters $(c,a)\in\C^{2}$ for which the \He map
$H_{c,a}$ has a fixed point with one eigenvalue $\lambda_{t}$ is a
curve of equation 
\begin{equation}\label{eq:ct}
c=
(1-a^{2})\left(\frac{\lambda_{t}}{2}-\frac{a^{2}}{2\lambda_{t}}\right)-
\left(\frac{\lambda_{t}}{2}-\frac{a^{2}}{2\lambda_{t}}\right)^{2}.
\end{equation}
\end{prop}

\noindent \textbf{Notations and conventions.} The
curve $\mathcal{P}_{\lambda_{t}}$ has degree $4$ in the variable $a$, and degree $2$ as a function of the Jacobian, which is $-a^2$. For this reason, we will sometimes call the curves $\mathcal{P}_{\lambda_{t}}$ complex parabolas.
When $t=0$, the curve
$\mathcal{P}_{\lambda}$ contains the \He maps that have a
semi-parabolic fixed point with one eigenvalue $\lambda$, a root of
unity.

For the rest of the paper, we denote by $c_t(a)$ the right hand side of Equation \eqref{eq:ct}. The \He map is completely determined by the choice of  $a$ and $t$,
so we will use $H_{a,t}$ in place of $H_{c_t(a),a}$ when there is
no danger of confusion. We write
\begin{equation}\label{eq: Henon-form}
H_{a,t}\hvec{x}{y}=\hvec{x^{2}+c_t(a)+ay}{ax}=\hvec{x^2+c_t+a^2w+ay}{ax},
\end{equation}
where the residual term $w$ is bounded and depends only on $a$ and $\lambda_{t}$, 
\begin{equation*}
w:=\frac{-1+\lambda_t-\lambda^2_t}{2\lambda_t}+\frac{a^2}{2\lambda_t}\left(1-\frac{1}{2\lambda_t}\right).
\end{equation*}

The \He map is also determined by the eigenvalues $\lambda$ and $\nu$ at a fixed point and we will sometimes write $H_{\lambda, \nu}$ in place of $H_{c,a}$ to stress this dependency. The  formula for $H_{\lambda,\nu}$ is the following: 
\[
H_{\lambda,\nu}(x,y) = \left(x^{2}+ (\lambda+\nu)(2+2\lambda\nu-\lambda-\nu)/4\pm i\sqrt{\lambda\nu}\,y,  \pm i\sqrt{\lambda\nu}\,x\right).
\] 
It may seem that there are two choices, but they are in fact conjugated by the affine change of variables $(x,y)\mapsto(x,-y)$.

A simple analysis shows that any constant $r>3$ works in the definition of sets $V^{\pm}$ from \eqref{eq:V+} for the whole family of \He maps $H_{a,t}$ for $|a|$ and $|t|$ small. From now on, we assume that $r>3$ is a fixed constant. Moreover, we assume that $|t|<1/(2q)$ and $|a|<1/2$ as minimal requirements and we will specify other restrictions when necessary. 

Let $\fq_{a,t}$ denote the fixed point of $H_{a,t}$ which has one eigenvalue $\lambda_t$. Suppose $|a|$ and $|t|$ are sufficiently small. 
We will see that the following bifurcation occurs:
\begin{itemize}
\item[(a)] if $t=0$ then $H_{a,0}$ has a semi-parabolic fixed point $\fq_{a,0}$
of multiplicity $q+1$.
\item[(b)] if $t>0$ then $H_{a,t}$ has a hyperbolic fixed point $\fq_{a,t}$ and
a $q$-periodic attractive orbit.
\item[(c)] if $t<0$ then $H_{a,t}$ has an attracting fixed point $\fq_{a,t}$ and
a $q$-periodic hyperbolic cycle.
\end{itemize}

In the degenerate case when $a=0$, the fixed point $\fq_{0,t}$ is just $(\alpha_t,0)$,
where $\alpha_t$ is the fixed point of the polynomial $p_t$. The \He maps becomes $H_{0,t}\hvec{x}{y}=\hvec{p_t(x)}{0}
$. The Julia set of the \He map $H_{0,t}$ is just the Julia set of the
polynomial $p_t$, so we write $J_{0,t}=J_{p_t}$. Moreover $J^+_{0,t}=J_{p_t}\times \C$.

\subsection{Local Dynamics -- Perturbed Normal Forms}\label{sec:localdynamics}

In this section we give a normal form for perturbations of semi-parabolic germs with semi-parabolic multiplicity $1$. This provides an analogue of Proposition \ref{prop:normal-polynom-t} for the two dimensional setting. Proposition \ref{prop:perturbed-form} is more general. Theorem \ref{thm:UniformNF} is specialized to the case of \He maps $H_{a,t}$ that come from perturbations of the family of polynomials $p_t$ discussed in the previous chapter. 

\begin{defn} Let $H$ be a holomorphic germ of $(\C^2,\fq)$ whose eigenvalues at the fixed point $\fq$ are $\lambda$ and $\nu$, with $0<|\nu|<1$ and $|\nu|<|\lambda|$. The strong stable manifold of  the fixed point $\fq$ corresponding to the eigenvalue $\nu$ is  
\begin{equation}\label{eq:Wss}
W^{ss}(\fq)=\{ z\in \C^2 : \lim\limits_{n\rightarrow \infty}|\nu|^{-n}\mbox{dist}(H^{\circ n}(z), \fq)={\rm const.}\}.
\end{equation}
\end{defn}

We refer to \cite{S} and \cite{MNTU} for a consistent treatment of strong stable manifolds. For the \He map $H_{a,t}$, the strong stable manifold $W^{ss}(\fq_{a,t})$ contains points that get attracted to the fixed point $\fq_{a,t}$ at an exponential rate $|\nu_{a,t}|^n$. When $\fq_{a,t}$ is hyperbolic or semi-parabolic, $W^{ss}(\fq_{a,t})$ lives in $J^+$. When the fixed point is attracting, or semi-Siegel, $W^{ss}(\fq_{a,t})$ belongs to the interior of $K^+$.
Note also that when $a=0$, $W^{ss}(\fq_{0,t})=\fq_{0,t}\times\C$.

We will denote by $W^{ss}_{loc}(\fq_{a,t})$ the local strong stable manifold of $\fq_{a,t}$ relative to the polydisk $\D_r\times \D_r$ that is, the connected component of $W^{ss}(\fq_{a,t})\cap \D_r\times \D_r$ that contains the point $\fq_{a,t}$.

Let $H$ be a semi-parabolic germ of transformation of $(\C^{2},0)$, with an isolated fixed point at $0$ with eigenvalues $|\nu|<1$ and  $\lambda=e^{2\pi i p/q}$. The multiplicity of $0$ as a solution of the equation $H^{\circ q}(x)=x$ is a number congruent to $1$ modulo $q$. Suppose therefore that $x=0$ is a fixed point of $H^{\circ q}$ of multiplicity $mq+1$; we call $m$ the \textit{semi-parabolic multiplicity} of $H$. 

\begin{prop}\label{prop:perturbed-form}
Let $\{H_{t}\}_{|t|<\delta'}$ be an analytic family of germs of diffeomorphisms of $(\C^{2},0)$ whose eigenvalues at $0$ are $\lambda_{t}=(1+t)\lambda$ and $\nu_{t}$, with $|\nu_{t}|<\min(1, |\lambda_t|)$ and $|\nu_{t}||\lambda_{t}|^{2 q}<1$. If the semi-parabolic multiplicity of $H_{0}$ is 1 then there exist local coordinates $(x,y)$ in which $H_{t}$ has the form $H_{t}(x,y)=(x_{1},y_{1})$, with
\begin{equation}\label{eq:NF3}
\left\{\begin{array}{l}
    x_{1}= \lambda_{t} (x+x^{ q+1} + Cx^{2 q+1}+ a_{2 q+2}(y)x^{2 q +2}+\ldots )\\
    y_{1}= \nu_{t} y + xh(x,y)
\end{array}\right.
\end{equation}
where $C$ is a constant depending on $t$, and $a_{i}(\cdot)$ and $h(\cdot,\cdot)$ are germs of holomorphic functions from $(\C,0)$ to $\C$, respectively from $(\C^{2},0)$ to $\C$, such that $a_{1}(0)=\lambda_{t}$ and $h(0,0)=0$. The coordinate transformations depend smoothly on $t$.
\end{prop}
\proof The case  $\lambda=1$ and $t=0$ was proved by Ueda 
\cite{U} and Hakim  \cite{Ha}. In \cite[Proposition~3.3]{RT}, this was stated for any primitive root of unity
$\lambda=e^{2\pi i p/q}$ and $t=0$. 

By straightening the local strong stable manifold of the fixed point $0$ we can assume that $H_t$ is written in the form:
\begin{equation}\label{eq:UniformNF1}
\left\{\begin{array}{l}
    x_{1}= a_{1}(y)x+a_{2}(y)x^{2}+\ldots \\
    y_{1}=\nu_t y + xh(x,y)
\end{array}\right.,
\end{equation}
where $a_{j}(\cdot)$ and $h(\cdot,\cdot)$ are holomorphic functions with $a_{1}(0)=\lambda_t$ and $h(0,0)=0$.

One can make a holomorphic change of coordinates to make the first $2q+1$ coefficients of the power series in the first coordinate constants. We proceed in three steps. In the first two steps we show that there exist local coordinates $(x,y)$ in which the map $H_t$ has the form:
\begin{equation}\label{eq:NF2}
\left\{\begin{array}{l}
    x_{1}= \lambda_t x+a_{2}x^{2}+\ldots + a_{2q+1}x^{2q+1} + a_{2q+2}(y)x^{2q+2}+\ldots \\
    y_{1}= \nu_t y + xh(x,y)
\end{array}\right.
\end{equation}
where $a_{2},\ldots, a_{2q+1}$ are constants. In the third step we show how to eliminate the terms $a_{k}x^{k}$ for $2\leq k\leq 2q+1$, $k$ not congruent to $1$ modulo $q$, and obtain Equation \eqref{eq:NF3}.

\noindent $(1)$ Reduction to $a_{1}(y)=\lambda_t$. Consider as in \cite{Ha} and \cite{U} a coordinate transformation
\begin{equation*}
\left\{\begin{array}{l}
    X=u(y)x \\
    Y=y
\end{array}\right. \ \ \ \mbox{with inverse}\ \ \
\left\{\begin{array}{l}
    x=X/u(Y) \\
    y=Y
\end{array}\right.
\end{equation*}
where $u$ is a germ of analytic functions from $(\C,0)$ to $\C$ with $u(0)=\lambda_t$.  We need to find $u$ such that
\begin{eqnarray*}
X_{1}&=& u(y_{1})x_{1}=u(\nu_t y +x h(x,y))\left(a_{1}(y)x+a_{2}(y)x^{2}+\ldots\right)\\
    &=& u(\nu_t Y + X/u(Y)h(X/u(Y) ,Y) )\left(a_{1}(Y)X/u(Y) +a_{2}(Y)(X/u(Y))^{2}+\ldots\right)\\
    &=& \frac{u(\nu_t Y)a_{1}(Y)}{u(Y)}X + \bigO(X^{2}) = \lambda_t X + \bigO(X^{2}).
\end{eqnarray*}
Let $b_1(Y)=a_1(Y)/\lambda_t$. The map $u$ satisfies the equation $u(Y) = u(\nu_t Y)b_1(Y)$. We successively substitute $\nu_t Y$ instead of $Y$ in this equation and obtain the unique solution
\begin{equation*}\label{eq:u(Y)}
u(Y) = \prod\limits_{n=0}^{\infty} b_1(\nu_{t}^{n}Y).
\end{equation*}
This product converges in a neighborhood of $0$ since $|\nu_t|<1$ and $b_{1}(Y) =1+\bigO(Y)$.

\smallskip
\noindent $(2)$ Reduction to $a_{k}(y)$ constants for $ 2\leq k\leq 2q+1$. We proceed by induction on $k$. The base case $k=1$ was discussed above. Suppose that $k\geq 2$ and that  there exist local coordinates $(x,y)$ in which $H_t$ has the form
\begin{equation*}
\left\{\begin{array}{l}
    x_{1}= \lambda_t x+a_{2}x^{2}+\ldots + a_{k-1}x^{k-1} + a_{k}(y)x^{k}+\ldots \\
    y_{1}= \nu_t y + xh(x,y),
\end{array}\right.
\end{equation*}
with $a_{2}, \ldots , a_{k-1}$  constant. We would like to find local coordinates so that $a_{k}(y)$ is also constant.  Consider the transformation
\begin{equation*}
\left\{\begin{array}{l}
    X=x+v(y)x^{k} \\
    Y=y
\end{array}\right. \ \ \ \mbox{with inverse}\ \ \
\left\{\begin{array}{l}
    x=X-v(Y)X^{k} + \ldots \\
    y=Y
\end{array}\right.
\end{equation*}
where $v$ is a germ of analytic functions from $(\C,0)$ to $\C$ with $v(0)=0$.  Using the coordinates given by this transformation we get
\begin{eqnarray*}
X_{1}&=& x_{1} + v(y_{1})x_{1}^{k}\\
    &=& \lambda_t x+a_{2}x^{2}+\ldots + a_{k-1}x^{k-1} + \left(a_{k}(y)+\lambda_t^k v(\nu_t y)\right)x^{k}+\bigO(x^{k+1})\\
    &=& \lambda_t X +  \ldots + a_{k-1}X^{k-1}+\left(a_{k}(Y) + \lambda_t^k v(\nu_t Y) - \lambda_t v(Y)\right)X^{k} + \bigO(X^{k+1})
\end{eqnarray*}
We need $v$ such that the coefficient of $X^{k}$ is constant. This gives the functional equation $\lambda_t v(Y)-\lambda_t^k v(\nu_t Y)=a_{k}(Y)-a_{k}(0)$. 
We successively substitute $\nu_t Y$ instead of $Y$ in this equation and obtain
\begin{equation*}\label{eq:v(Y)}
\lambda_t v(Y) = \sum\limits_{n=0}^{\infty} ( a_{k}\left( \nu_{t}^{n}Y)-a_{k}(0) \right)\lambda_t ^{n(k-1)}.
\end{equation*}
The series converges in a neighborhood of $0$  if $|\nu_t||\lambda_t|^{k-1}<1$. This is clearly achieved when $|\lambda_t|\leq 1$ since $|\nu_t|<1$. If $|\lambda_t|>1$ then $|\nu_t||\lambda_t|^{k-1}\leq |\nu_t||\lambda_t|^{2q}<1$ and the later inequality is true by hypothesis. Therefore  $H_t$ can be written as in Equation \eqref{eq:NF2}.

\smallskip
\noindent $(3)$ Suppose $2\leq k\leq 2q+1$ is not congruent to $1$ modulo $q$.
Assume by induction on $k$ that $H_t$ can be written as 
\begin{equation*}
\left\{\begin{array}{l}
    x_{1}= \lambda_t x+a_{k}x^{k}+\ldots + a_{2q+1}x^{2q+1} + a_{2q+2}(y)x^{2q+2}+\ldots \\
    y_{1}= \nu_t y + xh(x,y).
\end{array}\right.
\end{equation*}
Let $b = \frac{a_{k}}{\lambda_t-\lambda_t^{k}}$ and 
consider the coordinate transformation
\begin{equation}\label{eq: transf}
\left\{\begin{array}{l}
    X=x+bx^{k} \\
    Y=y
\end{array}\right. \ \ \ \mbox{with inverse}\ \ \
\left\{\begin{array}{l}
    x=X-bX^{k}+\ldots \\
    y=Y
\end{array}\right.
\end{equation}

In the new coordinate system, we get
\begin{eqnarray*}
X_{1}&=&x_{1}+bx_{1}^{k} = ( \lambda_t x+a_{k}x^{k}+\ldots )+ b(\lambda_t x+a_{k}x^{k}+\ldots )^k\\
    &=& \lambda_t x +(a_{k}+b\lambda_t^{k})x^{k} + \ldots\\
    &=&\lambda_t ( X -  b X^{k} + \ldots) +(a_{k}+b\lambda_t^{k})(X-bX^{k} + \ldots)^{k} + \ldots\\
    &=&\lambda_t X  +  (a_{k}+b(\lambda_t^{k}-\lambda_t)) X^{k}+ \ldots
\end{eqnarray*}
By this procedure, the term containing $X^k$ has been eliminated. The first monomial that cannot be eliminated in this way will be $a_{m q+1}X^{m q+1}$ for some integer $m$. If we assume that the parabolic multiplicity of the semi-parabolic germ is $1$, then $a_{q+1}\neq 0$ for $t=0$ (hence also for small $t$) and $a_{q+1}X^{q+1}$ will be the first term that we cannot eliminate by the above procedure. We can further reduce the normal form to $a_{q+1}=1$ by considering a linear transformation of the form $X=Ax, Y=y$, where $A$ is a constant such that $A^{ q}=a_{q+1}$.
We can therefore assume that the \He map can be written as 
\begin{equation*}
\left\{\begin{array}{l}
    x_{1}= \lambda_t (x+x^{q+1}+a_{k}x^{k}+\ldots + a_{2q+1}x^{2q+1} + a_{2q+2}(y)x^{2q+2}+\ldots )\\
    y_{1}= \mu y + xh(x,y).
\end{array}\right.
\end{equation*}
By repeating the coordinate transformations \eqref{eq: transf}, we can eliminate all monomials $a_k x^{k}$ with $q+1<k<2q+1$.
By abuse of notation we still denote by $a_k$ the term $\lambda_t a_k$. In the new coordinate system, we get
\begin{eqnarray*}
X_{1}&=&x_{1}+bx_{1}^{k} = ( \lambda_t (x+x^{q+1})+a_{k}x^{k}+\ldots )+ b(\lambda_t (x+x^{q+1})+a_{k}x^{k}+\ldots )^k\\
    &=& \lambda_t (x+x^{q+1}) +(a_{k}+b\lambda_t^{k})x^{k} + \ldots\\
    &=&\lambda_t ( X -  b X^{k} + (X -  b X^{k})^{q+1}+\ldots) +(a_{k}+b\lambda_t^{k})(X-bX^{k} + \ldots)^{k} + \ldots\\
    &=&\lambda_t (X+X^{q+1})  +  (a_{k}+b(\lambda_t^{k}-\lambda_t)) X^{k}+ \ldots
\end{eqnarray*}
therefore the term containing $X^k$ has been eliminated.
\qed

This following theorem is a generalization of \cite[Theorem~6.2]{RT}. 

\begin{thm}\label{thm:UniformNF}
Let $r>3$ be a fixed constant. There exist $\delta, \delta'>0$ such that for
any $(c,a)\in \mathcal{P}_{\lambda_{t}}$ with $|a|<\delta$ and $|t|<\delta'$ 
there exists a coordinate transformation $\phi_{a,t}$
from a tubular neighborhood $B=\D_{\rho'}(\alpha_{t})\times\D_{r}$
of the local strong stable manifold of the fixed
point $\fq_{a,t}$
\[
    \phi_{a,t}:B\rightarrow \D_{\rho}\times \D_{r+\bigO(|a|)}
\]
in which $H_{a,t}$ has the form
$\widetilde{H}_{a,t}(x,y)=(x_{1},y_{1})$, with
\begin{equation}\label{eq:UniformNF3}
\left\{\begin{array}{l}
    x_{1}= \lambda_t (x+x^{q+1} + C_{a,t}x^{2q+1}+ a_{2q+2}(y)x^{2q +2}+\ldots )\\
    y_{1}= \nu_{a,t} y + xh_{a,t}(x,y)
\end{array}\right.
\end{equation}
and $C_{a,t}$ is a constant (depending on $a$ and $t$) and
$xh_{a,t}(x,y)=\bigO(a)$. Moreover the transformations $\phi_{a,t}$ are analytic in $a$ and $t$, and
\[
    \lim\limits_{a\rightarrow 0}\phi_{a,t}=(\phi_t(x),y),
\]
uniformly with respect to $t$. The map
 $\phi_t:\D_{\rho'}(\alpha_t)\rightarrow
\D_{\rho}$
  is the change of coordinates from Lemma \ref{prop:normal-polynom-t} for the polynomial $p_t$ with a fixed point $\alpha_{t}$ of multiplier $\lambda_t$ and
\[
    \phi_t \circ p_t\circ \phi_t^{-1} (x)=\lambda_t (x+x^{q+1}+C_t x^{2
  q+1}+\bigO(x^{2q+2})).
\]
\end{thm}
\proof  
We choose $\delta$ small enough so that the local strong stable manifold $W^{ss}_{loc}(\fq_{a,t})$ has no foldings inside $\D_r\times \D_r$ and it can therefore be straightened using a holomorphic change of coordinates. The first part of the proof follows directly from Proposition  \ref{prop:perturbed-form}. 
We need to verify the conditions imposed on the eigenvalues at the fixed point.  We have $|\lambda_t|\in (1-\delta',1+\delta')$ and $|\nu_{a,t}|=|a|^2/|\lambda_t|\in [0,\delta^2/(1-\delta'))$. The bounds $\delta$ and $\delta'$ are chosen small enough so that $\delta\ll 1- \delta'$. It follows that $|\nu_{a,t}|\ll |\lambda_t|$ and $|\nu_{a,t}|<1$. Then $|\lambda_t||\nu_{a,t}|^{2q}< (1+\delta')\delta^{4q}(1-\delta')^{-2q}<1$. This inequality is not so restrictive. For example, it is verified for $\delta'<1/(2q)$ and $\delta<1/2$. The convergence of the coordinate transformation $\phi_{a,t}$ as $a\rightarrow 0$ follows immediately by comparing the coordinate transformations done in Proposition \ref{prop:perturbed-form} to those done in Proposition \ref{prop:normal-polynom-t}.
\qed

It is also worth mentioning that the change of coordinates function
$\phi_{a,t}$ from Theorem \ref{thm:UniformNF} maps horizontal curves to horizontal curves, that is 
\begin{equation}\label{eq:horcurves}
\phi_{a,t}(\D_{\rho'}(\alpha_t)\times \{y_1\})\subset \C\times
\{y_2\},
\end{equation} 
which will be useful later on. 

\subsection{Attracting and repelling sectors}\label{subsec:horizontal} 

In this section we continue the analysis of the local dynamics of holomorphic germs of diffeomorphisms of $(\C^2,0)$ with a semi-parabolic fixed point and their nearby perturbations.
Consider the set 
\[
\Delta_{R}= \left\{ x\in \C\ :\
\left(Re(x^{q})+\frac{1}{2R}\right)^{2} +
\left(|Im(x^{q})|-\frac{1}{2R}\right)^{2} <
\frac{1}{2R^{2}}\right\}.
\]
in the complex plane. There are $q$ connected components of $\Delta_{R}$, which we denote
$\Delta_{R,j}$, for $1\leq j\leq q$. Define $\Pet_{R,r} = \Delta_{R}\times \D_{r}$ and let $\Pet_{R,r,j} = \Delta_{R,j}\times \D_{r}$ be the connected components of $\Pet_{R,r}$.

\begin{prop}\label{prop:sectors-t}
For $R$ large enough and $r$ small enough there exists a positive
number $\delta'$ such that for all $t\in (-\delta',\delta')$
\[
H_{a,t}(\overline{\Pet}_{R,r,j})\subset \Pet_{R,r,j+p}\cup
\{0\}\times \D_{r}\ \ \ \ \mbox{for} \ 1 \leq j \leq  q.
\]
In particular $H_{a,t}(\overline{\Pet}_{R,r})\subset \Pet_{R,r}\cup
\{0\}\times \D_{r}$.
\end{prop}
\proof 
Assume that $R$ is large enough and $r$ is small enough so that the map
$H_{a,t}$ is well defined and has the expansion from Theorem
\ref{thm:UniformNF}:
\begin{equation*}
\left\{\begin{array}{l}
    x_{1}= \lambda_{t} (x+x^{q+1} + C_{a,t}x^{2q+1}+ a_{2q+2}(y)x^{2q +2}+\ldots )\\
    y_{1}= \nu_{a,t} y + xh(x,y)
\end{array}\right. .
\end{equation*}
Define the region $ U_{R}:=\left\{ X\in \C \ |\ R/q-Re(X)<|Im(X)|\right\}$ and set 
$W_{R,r}:=U_{R}\times \D_{r}$. 

Suppose $(x,y)\in \Pet_{R,r,j}$. The 
transformation $X=-1/(qx^{q})$, $Y=y$
maps each petal $\Pet_{R,r,j}$ to $W_{R,r}$. Thus $X\in U_{R}$ and $|Y|<r$. Let
$\widehat{H}(X,Y)=(X_{1},Y_{1})$ be the corresponding map in the new coordinates:
\begin{eqnarray*}
X_{1}&=&-\frac{1}{qx_{1}^{q}}= \frac{X}{\lambda_{t}^{q}\left(1+x^{q} + C_{a,t}x^{2q}+ a_{2q+2}(y)x^{2q +1}+\ldots \right)^{q}}\\
    &=&\frac{X}{\lambda_{t}^{q}}\left(1-q(x^{q} + C_{a,t}x^{2q}+\ldots) + \frac{q(q+1)}{2}x^{2q}+\ldots \right)\\
    &=&\frac{1}{\lambda_{t}^{q}}\left(X+1+\frac{A}{X}+\bigO_{Y}\left(\frac{1}{|X|^{1+1/q}}\right)\right),\  \mbox{where}\ A := \frac{1}{q}\left(\frac{q+1}{2}-C_{a,t}\right);\\
Y_{1} &=&\nu_{a,t} y + xh(x,y) = \nu_{a,t} Y +
\bigO_{Y}\left(\frac{1}{|X|^{1/q}}\right).
\end{eqnarray*}
The notation $\bigO_{Y}\left(|X|^{\alpha}\right)$
represents a
holomorphic function of $(X,Y)$ in $W_{R,r}$ which is bounded by $K |X|^{\alpha}$
for some constant $K$.

One can easily check that $|X|>\frac{R}{q\sqrt{2}}$ throughout the region $U_{R}$. Clearly $|\lambda_{t}|^{q}>1/2$ for small $|t|$.
There exist constants $K', K''$ and $K_{1}=2K'q\sqrt{2}$, $K_{2}=K''(q\sqrt{2})^{1/q}$ such that
\begin{eqnarray*}
\bigg{|}X_{1}-\frac{1}{\lambda_{t}^{q}}(X+1)\bigg{|} &\leq& \frac{K'}{|\lambda_{t}|^{q}|X|}<\frac{K_{1}}{R}\\
|Y_{1}-\nu_{a,t} Y| &\leq& \frac{K''}{|X|^{1/q}}<\frac{K_{2}}{R^{1/q}}.
\end{eqnarray*}
Choose $R$ large enough and $r$ small enough so that
\begin{equation}\label{eq:K1K2}
\left\{\begin{array}{ll}
    \ds \frac{K_{1}}{R}<\frac{1}{4}  \\
    \vspace{-0.35cm}\\
   \ds  \frac{K_{2}}{R^{1/q}}<(1-|\nu_{a,t}|)r \ .
\end{array}\right.
\end{equation}
The second condition immediately gives 
\[
    |Y_{1}|\leq |Y_{1}-\nu_{a,t} Y|+|\nu_{a,t}| |Y|<\frac{K_{2}}{R^{1/q}}+|\nu_{a,t}| r <r.
\]
The first condition of \eqref{eq:K1K2}  implies that
\begin{equation}\label{eq:1/4}
\bigg{|}X_{1}-\frac{1}{\lambda_{t}^{q}}(X+1)\bigg{|}<\frac{1}{4}.
\end{equation}
In our case $\lambda_{t}^{q} = (1\pm t)^{q}$ is a real positive number. Hence inequality \eqref{eq:1/4} yields the following estimates:
\begin{eqnarray*}
Re(X_{1}) &>& \frac{1}{(1\pm t)^{q}}Re(X)+\frac{1}{(1\pm t)^{q}}-\frac{1}{4}\\
|Im(X_{1})| &>& \frac{1}{(1\pm t)^{q}}|Im(X)|-\frac{1}{4}.
\end{eqnarray*}
Using these estimates and the fact that $X\in U_{R}$, we get 
\[
|Im(X_{1})|>\frac{R/q+1}{(1\pm t)^{q}}-\frac{1}{2}-Re(X_{1})>R/q-Re(X_{1}),
\]
provided that $|t|$ is small enough so that
\begin{equation}\label{eq:t-small}
\frac{R+q}{R+q/2}>(1\pm t)^{q}.
\end{equation}
The constant $\delta'>0$ is chosen so that this inequality holds for all $0\leq t<\delta'$. It follows that $\widehat{H}(\overline{W}_{R,r})\subset W_{R,r}$. 
\qed

Since the dynamics is different, we will treat the cases $t<0$ and $t>0$ separately. We 
begin with the latter case, so suppose $R,r$ and $t>0$ are as in Proposition \ref{prop:sectors-t}. 
We  show that points in $\Pet_{R,r}$ are attracted by an
attractive orbit of period $q$ under iterations by $H_{a,t}$. Each
region $\Pet_{R,r,j}$ contains a point of this orbit.  The fixed
point $0$ is hyperbolic. 

\begin{figure}[htb]
\begin{center}
\includegraphics[scale=0.37]{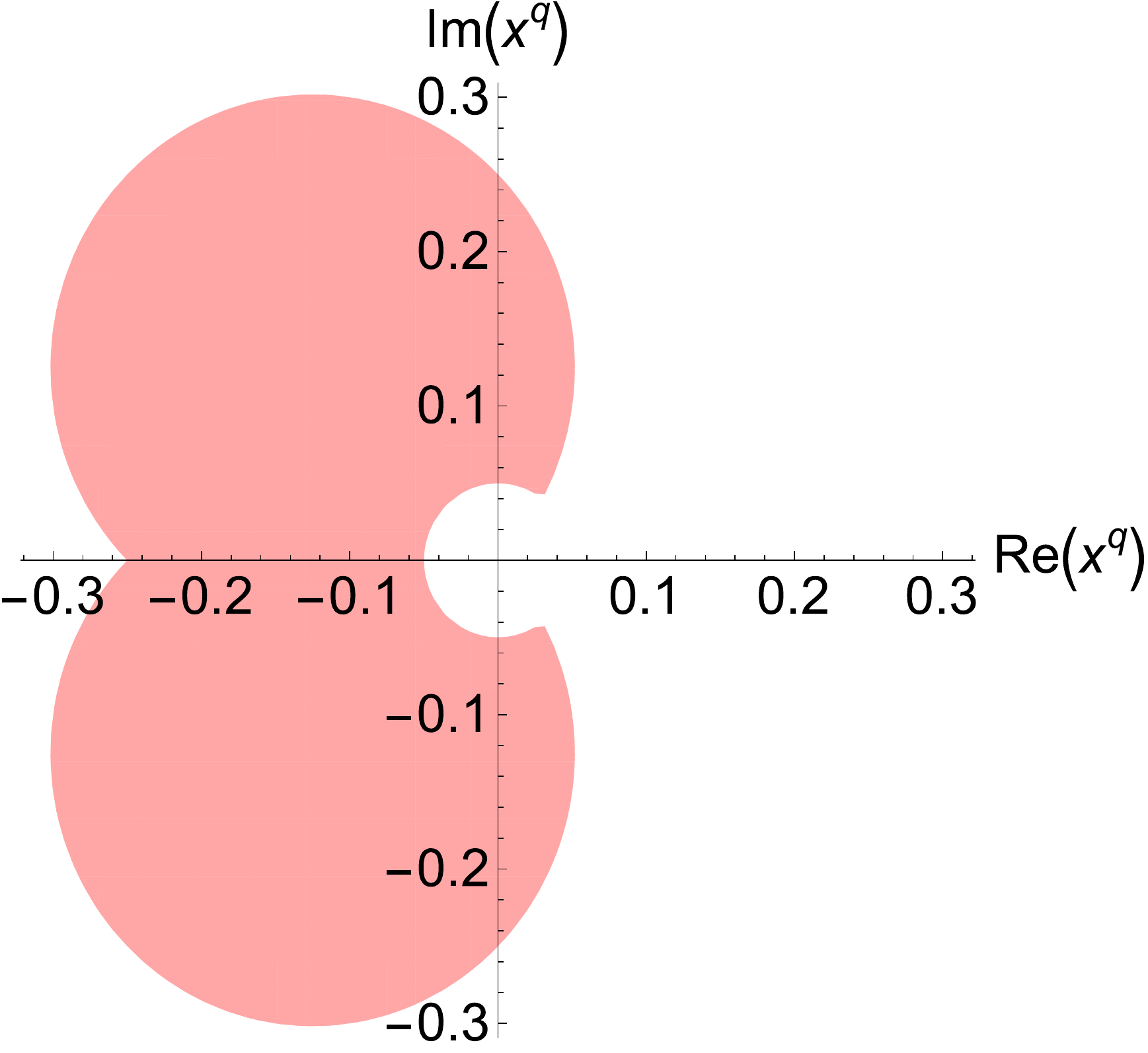}
\end{center}
\caption{(case $t>0$) The image of $\mathcal{D}_{R,r,t}$ under the map $x\mapsto x^{q}$ at a height $y={\rm const}$. 
A small disk of radius $\rho_{t}$ is removed around the origin.}
\label{fig:trapping0}
\end{figure}

Suppose $\rho_{t}>0$ is a small enough radius such that 
$\rho_{t}\leq \sqrt[q]{t/(2q)}$.  The number $\rho_{t}$ is just a local variable which will be used in the proposition below. 
Define
 \[
\mathcal{D}_{R,r,t}= \Pet_{R,r}-\left\{ (x,y)\in \C^{2}\ :\
|x|^{q}\leq \rho_{t}, |y|<r \right\}
\]
and let  $\mathcal{D}_{R,r,t,j}$ with $1\leq j \leq q$, be the
connected components of $\mathcal{D}_{R,r,t}$ (see Figure \ref{fig:trapping0}).

\begin{prop}[{\bf Trapping regions -- $t$ positive}]\label{prop:trapping}
For $R$ large enough and $r$ small enough there exists a positive
number $\delta'$ such that for all $t\in (0,\delta')$
\[
H_{a,t}( \overline{\mathcal{D}}_{R,r,t,j} ) \subset
\mathcal{D}_{R,r,t,j+p}\ \ \ \ \mbox{for} \ 1 \leq j \leq  q.
\]
In particular $H_{a,t}(\overline{\mathcal{D}}_{R,r,t})\subset
\mathcal{D}_{R,r,t}$ and all points of $\mathcal{D}_{R,r,t}$ are
attracted to an attractive orbit of period $q$ under iterations by
$H_{a,t}$.
\end{prop}
\proof We make a change of variables $X=-1/(qx^{q})$, $Y=y$ and
analyze the situation at infinity. This transformation maps each
$\mathcal{D}_{R,r,t,j}$ to a region $W_{R,r,t}:=U_{R,t}\times
\D_{r}$, where
\[
U_{R,t}:=\left\{ X\in \C \ :\ R/q-Re(X)<|Im(X)|\ \mbox{and}\
|X|<\frac{1}{q \rho_{t}^{q}}\right\}.
\]
From Equation \eqref{eq:1/4} from the proof of the previous
proposition we have
\[
|X_{1}|-\frac{1}{(1+t)^{q}}|X|-\frac{1}{(1+t)^{q}}\leq
\bigg{|}X_{1}-\frac{1}{(1+t)^{q}}(X+1)\bigg{|}<\frac{1}{4}
\]
which gives
\[
|X_{1}|<\frac{1}{q(1+t)^{q}\rho_{t}^{q}}+\frac{1}{(1+t)^{q}}+\frac{1}{4}<\frac{1}{q\rho_{t}^{q}}.
\]
The last inequality holds because
\[
q\rho_{t}^{q}<\frac{t}{2}<\frac{(1+t)^{q}-1}{(1+t)^{q}/4+1},
\]
based on our assumption on $t$ and our choice of $\rho_{t}$. 
With this choice we showed that
$\widehat{H}(\overline{W}_{R,r,t})\subset W_{R,r,t}$ and all points
of $W_{R,r,t}$ are attracted to an attractive orbit under iterations
by $\widehat{H}$. The existence of this orbit follows immediately
since we have a nested intersection of compact sets.
\qed

By choosing a smaller $\rho_{t}$ as necessary we can show that
all points in $\Pet_{R,r}$ are attracted by the $q$-periodic
attractive orbit under forward iterations by $H_{a,t}$.  Moreover,
every point that is attracted to this orbit must eventually land in
the interior of one of the regions $\Pet_{R,r,j}$ for $1\leq j\leq q$.

Let $\epsilon_{0}=\tan(2\pi/9)$. In order to simplify notation, define $\rho$ such that $\rho^{q}:=\frac{1-\epsilon_{0}}{R\sqrt{1+\epsilon_{0}^{2}}}$. The number $\rho$ measures the distance between the origin and one of the points of intersection of the lines $Re(x^{q})= \epsilon_{0}|Im(x^{q})|$ with the boundary of $\Delta_{R}$. 

Define the {\it
attractive} sectors
\begin{equation}\label{eq:W+t}
    \Delta^{+} := \left\{ x\in \C\ :\ Re(x^{q})\leq \epsilon_{0}|Im(x^{q})| \
     \mbox{and}\ |x^{q}|< \rho^{q} \right\},
\end{equation}
and the {\it repelling} sectors
\begin{equation}\label{eq:W-t}
\Delta^{-}:= \left\{ x\in \C\ :\ Re(x^{q})> \epsilon_{0}|Im(x^{q})|\
\mbox{and} \ |x^{q}|< \rho^{q} \right\}.
\end{equation}
Let $W^{+}:=\Delta^{+}\times \D_{r}\subset\Pet_{R,r}$ and $W^{-}:=\Delta^{-}\times \D_{r}$.

We will call $W^{-}$ {\it repelling} because as we will see, the \He
map expands horizontally when the Jacobian is small enough. We will
call $W^{+}$ {\it attractive} because points in $W^{+}$ are
attracted to the $q$-periodic attractive orbit as we have shown
above. There are $q$ components of $W^{\pm}$ which we denote
$W^{\pm}_{j}$ for $1\leq j\leq q$.

\begin{figure}[htb]
\begin{center}
\includegraphics[scale=0.37]{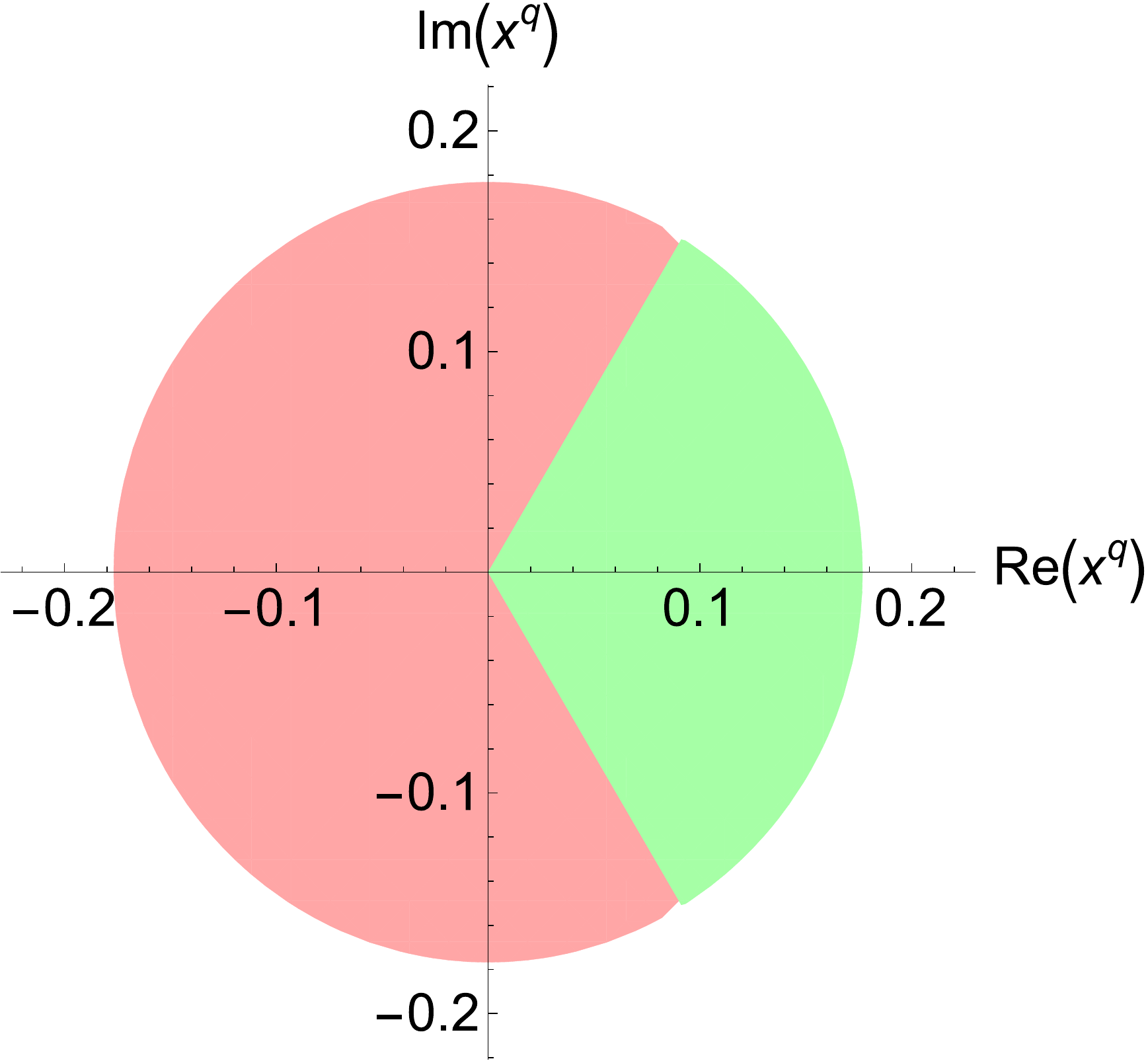}
\end{center}
\caption{(case $t\geq 0$) Repelling and attractive sectors near the semi-parabolic/hyperbolic fixed point. The attracting sector $\Delta^{+}$ is shown in red and the repelling sector $\Delta^{-}$ is shown in green. The angle of the green sector is $5\pi/9$.}
\label{fig:trapping}
\end{figure}

In the regions $\Delta^{-}$ and $W^{-}$ we have
\begin{equation}\label{eq:eps1}
    Re(x^{q})> \epsilon_{1}|x|^{q},\ \ \ \mbox{where}\ \  \epsilon_{1} := \frac{\epsilon_{0}}{\sqrt{1+\epsilon_{0}^{2}}}>\frac{3}{5}.
\end{equation}
The constants $\epsilon_{0}$ and $\epsilon_{1}$ are chosen such that
the image of $\Delta^{-}$ under $x\mapsto x^q$ has an angle opening
of $5\pi/9$ (see Figure \ref{fig:trapping}).

The definition of the sectors $W^{\pm}$ for $t>0$ is the same as in the case $t=0$ \cite[Section~4]{RT}. However,  when $t<0$, we need to modify the definition of the repelling sector $W^{-}$ so that we have a good horizontal expansion for the \He map. Suppose therefore that $t<0$. 

\begin{remarka}
When $\lambda=1$  the parametric paths described by $\lambda_{t}=1\pm t$ are in fact the same, so we can assume that $q\geq 2$. 
\end{remarka}

\noindent {\bf Assumption on $t$.} Suppose $t$ is sufficiently small so that $|t|<\frac{1}{24q+12}$. A restriction on $t$ of this form is needed for the local dynamics, but the choice for this bound will become clear later on. Let 
\begin{equation}\label{eq:Rt}
R_{t}:=\frac{|t|}{(q+1/3)\epsilon_{1}}
\end{equation} 
be fixed from now on. The constant $\epsilon_{1}>3/5$ is the same as in Equation \eqref{eq:eps1}. Suppose further that $|t|$ is small enough so that $R_{t}<1/(9R)$, where $R$ is as in Proposition \ref{prop:sectors-t}.

Define 
 \[
\mathcal{D}_{r,t}= \left\{ (x,y)\in \C^{2}\ :\
|x|^{q}\leq R_{t}, |y|<r \right\}.
\]

\begin{prop}[{\bf Attracting region -- $t$ negative}]\label{prop:trapping2origin}
For $R$ large enough and $r$ small enough there exists a positive
number $\delta'$ such that $H_{a,t}(\overline{\mathcal{D}}_{r,t})\subset
\mathcal{D}_{r,t}$ and all points of $\mathcal{D}_{r,t}$ are
attracted to the origin under iterations by $H_{a,t}$ for all $t\in (-\delta',0)$.
\end{prop}
\proof We make the change of variables $X=-1/(qx^{q})$, $Y=y$ and
analyze the situation at infinity. Let $\widehat{H}$ be the map written in these coordinates. This transformation maps 
$\mathcal{D}_{r,t}$ to the region
\[
W_{r,t}:=\left\{ X\in \C \ :\ |X|\geq \frac{1}{qR_{t}}\right\}\times\D_{r}.
\]
Let $(X_{1},Y_{1})=\widehat{H}(X,Y)$ for $(X,Y)\in W_{r,t}$. Note that $\lambda_{t}=1-|t|$, as $t$ is negative.
Similar to Equation \eqref{eq:1/4} we have that 
\[
\bigg{|}X_{1}-\frac{1}{(1-|t|)^{q}}X\bigg{|}<\frac{1}{4}.
\]
This gives
\[
\frac{1}{(1-|t|)^{q}}|X|-\frac{1}{(1-|t|)^{q}}-|X_{1}|\leq
\bigg{|}X_{1}-\frac{1}{(1-|t|)^{q}}(X+1)\bigg{|}<\frac{1}{4}
\]
and, after rearranging the terms, we want to obtain the following estimate 
\[
|X_{1}|>\frac{1}{(1-|t|)^{q}}|X|-\frac{1}{(1-|t|)^{q}}-\frac{1}{4}>|X|+\frac{1}{40}.
\]
The last inequality is equivalent to 
\[ |X|\left(\frac{1}{(1-|t|)^{q}}-1\right)>\frac{1}{(1-|t|)^{q}}+\frac{11}{40}.
\]
Note that $\frac{1}{(1-|t|)^{q}}-1>q|t|$ for small $|t|$. Using the fact that $|X|\geq(qR_{t})^{-1}$ and the particular choice of $R_{t}$ we get that
\[
|X|\left(\frac{1}{(1-|t|)^{q}}-1\right)>\left(q+1/3\right)\epsilon_{1}\geq \frac{7}{5}>\frac{1}{(1-|t|)^{q}}+\frac{11}{40},
\]
which is true whenever $(1-|t|)^{q}>8/9$. This condition is satisfied because, based on our assumption on $t$, we have
\[ (1-|t|)^{q}>\left(1-\frac{1}{24q+12}\right)^{q}\geq \left(\frac{59}{60}\right)^{2}>\frac{8}{9}.
\]
Note that the function $x\mapsto(1-1/(24x+12))^{x}$ is increasing on $[2,\infty)$ and that is why we can use the middle inequality. We have therefore shown that $|X_{1}|>|X|+1/40$. Let  $(X_{n},Y_{n})=\widehat{H}(X_{n-1},Y_{n-1})$ for some $(X_{0},Y_{0})\in W_{r,t}$. By induction we get that $|X_{n}|>|X_{0}|+n/40$. It follows that all points in $W_{r,t}$ are attracted to $(\infty,0)$. In order to prove that indeed $Y_{n}\rightarrow 0$ as $n\rightarrow\infty$ we need to do a similar analysis as in \cite{RT}; we leave the details to the reader.
\qed

\begin{figure}[htb]
\begin{center}
\includegraphics[scale=0.37]{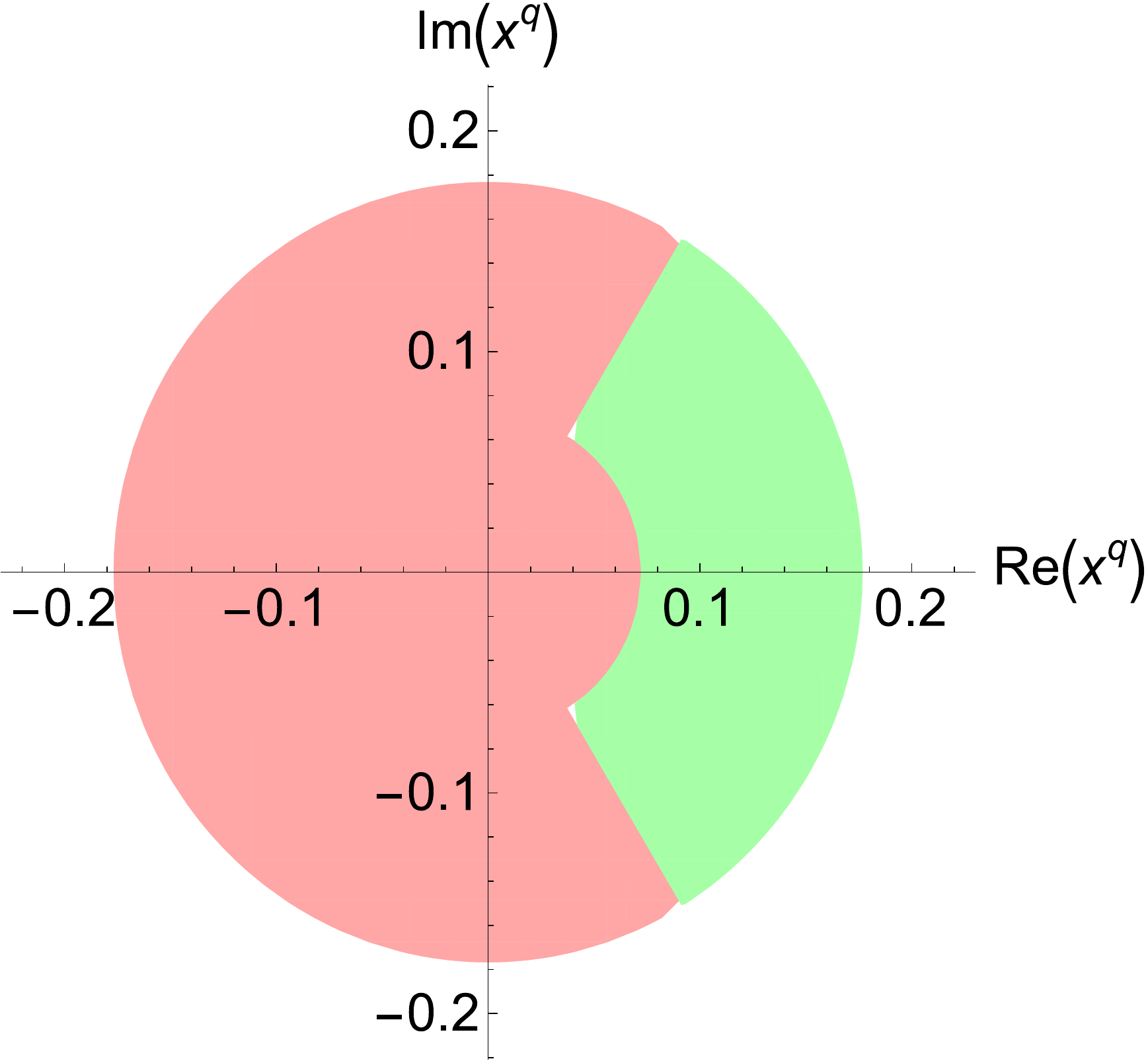}
\end{center}
\vspace{-0.15cm}
\caption{(case $t<0$) The repelling sector $\Delta^{-}_{R_{t}}$ is shown in green. The angle of the green sector is $5\pi/9$. The red region belongs to the basin of attraction of $0$. }
\label{fig:trapping2}
\vspace{-0.25cm}
\end{figure}

When $t<0$, we define the {\it repelling} sectors as follows
\begin{equation}\label{eq:W-t2}
\Delta^{-}_{R_{t}}:= \left\{ x\in \C\ :\ Re(x^{q})> \epsilon_{0}|Im(x^{q})|\
\mbox{and} \ R_{t}<|x^{q}|< \rho^{q} \right\},
\end{equation}
where $R_{t}$ is given in \eqref{eq:Rt}. The excluded region belongs to the basin of attraction of $0$ (see Figure \ref{fig:trapping2}). Set as before $W^{-}_{R_{t}}:=\Delta^{-}_{R_{t}}\times \D_{r}$ . The definition of the set $W^{+}$ is the same as in the case when $t$ is positive, i.e. $W^{+} = \Delta^{+}\times \D_{r}$, where $\Delta^{+}$ is given in Equation \eqref{eq:W+t}. Also, $\Delta^{-}_{R_{t}}$ is a subset of $\Delta^{-}$, defined in Equation \eqref{eq:W-t}. When $t\rightarrow 0^{-}$, the sets $\Delta^{-}_{R_{t}}$ converge to $\Delta^{-}$, so the definition of $W^{-}$ when $t=0$ is the same as in \cite{RT}. 

\smallskip
By choosing $t$ small enough so that $R_{t}<\frac{1}{9R}<\rho^{q}$ we made sure that the excluded region $\{x\in \C\ |\ |x|^{q}<R_{t}\}\cap\Delta_{R}$ is contained in $\Delta^{+}$.

\smallskip
Let $B=\D_{\rho'}(\alpha_{t})\times\D_{r}$ be the polydisc from Theorem \ref{thm:UniformNF} and $\phi_{a,t}$ the coordinate transformation defined on $B$. We define attractive and repelling sectors relative to $B$.

\begin{defn}\label{defn:WB} 
Let $W^+_B:=\phi_{a,t}^{-1}(W^+)$ for $t\geq0$ and $W^+_B:=\phi_{a,t}^{-1}(\D_{\rho}\times\D_{r}-W^{-}_{R_{t}})$ for $t<0$ be the attractive sectors inside $B$. Similarly, let $W^-_B:=\phi_{a,t}^{-1}(W^-)$ for $t\geq 0$ and $W^-_B:=\phi_{a,t}^{-1}(W^-_{R_{t}})$ for $t<0$ be the repelling sectors inside $B$. 
\end{defn}


\medskip
\begin{prop}[\textbf{Local dynamics}]\label{prop:localdyn}~\
\begin{itemize}
\item[a)] If $t\geq 0$  then the compact region $W^+_B$ satisfies $H_{a,t}(W_B^{+})\subset int(K_{a,t}^{+})\cup W^{ss}_{loc}(\fq_{a,t})$.
\item[b)] If $t<0$ then the compact region $W^+_B$ lies in the interior of $K_{a,t}^{+}$.
\end{itemize}
\end{prop}
\proof  Using the definitions from \ref{defn:WB}, the proof follows directly from Propositions \ref{prop:sectors-t} and \ref{prop:trapping2origin}.
\qed
\vspace{-0.15cm}

\subsection{Deforming the local semi-parabolic structure into a hyperbolic structure} \label{section:parabolic2hyperbolic}

In the parabolic case we have shown in \cite[Propositions~6.8, 9.2]{RT} that in the repelling sectors $W^-$ near a semi-parabolic fixed point, the \He map is weekly expanding in horizontal cones and strongly contracting in vertical cones, with respect to the Euclidean metric. We will reuse these cones and show that when $t$ is nonzero, the \He map is strongly expanding in horizontal cones and strongly contracting in vertical cones, and therefore has a local hyperbolic structure. In this section we only use the local normal form of the map, so all results are applicable to holomorphic germs of diffeomorhisms of $(\C^2,0)$ with a semi-parabolic fixed point at $0$.

\begin{defn}\label{def:vcones}
Define the \textbf{vertical cone} at a point $(x,y)$ from the set
$\D_{\rho}\times \D_r$ as
\[
\mathcal{C}^{v}_{(x,y)}=\left\{(\xi,\eta)\in T_{(x,y)}\D_{\rho}\times\D_r,\ |\xi|\leq |x|^{2q}|\eta|\right\}.
\]
Define the \textbf{horizontal cone} at a point $(x,y)$ from the set
$\D_{\rho}\times \D_r$ to be
\[
\mathcal{C}^{h}_{(x,y)}=\left\{(\xi,\eta)\in
T_{(x,y)}\D_{\rho}\times\D_r,\ |\xi|\geq |\eta|\right\}.
\]
\end{defn}
\noindent We consider the interior of a cone to be its topological interior together with the origin.

\smallskip
Consider the \He map $\widetilde{H}_{a,t}:\D_{\rho}\times \D_{r}\rightarrow \C^{2}$ written in the normal form given in Equation \eqref{eq:UniformNF3}. 
We write $\widetilde{H}_{a,t}$ whenever we want to stress the dependency on the parameters $a$ and $t$, but otherwise we simply write $\widetilde{H}$. We have
\[
\widetilde{H}_{a,t}\hvec{x}{y} = \hvec{\lambda_{t} (x+x^{q+1}+g_{a,t}(x,y))}{\, \nu_{a,t} y + xh_{a,t}(x,y)},
\]
where
\begin{eqnarray*}
    g_{a,t}(x,y)&=& C_{a,t}x^{2q+1}+ a_{2q+2}(y)x^{2q +2}+\ldots \\
    h_{a,t}(x,y)&=& b_{1} (y) +\ldots + b_{k}(y)x^{k}+\ldots
\end{eqnarray*}
and $g_{a,t}(x,y)=g_{0,t}(x)+\bigO(a)$ and $h_{a,t}(x,y)=\bigO(a)$. Here the term $\bigO(a)$ is in fact a holomorphic function in both $a$ and $t$.

When $a=0$, $\widetilde{H}_{0,t}(x,y)=(\widetilde{p}_{t}(x),0)$, where $\widetilde{p}_{t}(x)=\lambda_{t}(x+x^{q+1}+g_{0,t}(x))$ and
\[
g_{0,t}(x)=C_{0,t}x^{2q+1}+a_{2q+1}x^{2q+2}+\ldots .
\]
The function $g_{0,t}$ depends only on $x$ and $t$, hence $\partial_{y} g_{0,0}(x,y)\equiv 0$. For $|a|<\delta$ and $|t| <\delta'$ we
can assume that there exists a constant $M_{a,t}$ with $0<M_{a,t}<1$
such that 
\begin{equation}\label{eq:Mat}
\big{|}\partial_{y} g_{a,t}(x,y)\big{|}<M_{a,t}|x|^{2q+2}.
\end{equation}

When $a=0$ we also know that $xh_{0,t}(x,y)\equiv 0$.
Moreover, by the construction of the normalizing coordinates, we have
$xh_{a,t}(x,y)=\bigO(a)$. There exists a constant $N_{a,t}$, depending
on $a$,  with $0<N_{a,t}<1$ such that when $|a|<\delta$ the following
bounds hold
\begin{eqnarray}\label{eq:Nat}
\big{|}\partial_{x} ( xh_{a,t} )(x,y)\big{|} < N_{a,t} \ \ \mbox{and}\ \ \big{|}\partial_{y} ( xh_{a,t} ) (x,y)\big{|} < N_{a,t}.
\end{eqnarray}

Let $\partial_{x} g_{a,t}(x,y)=x^{2q}\beta_{a,t}(x,y)$, for some function $\beta_{a,t}$. As usual, $\partial_{x}$ denotes the partial derivative with respect to the variable $x$. 
Denote by $m$ the supremum of $|\beta_{a,t}(x,y)|$ on the set $W^{-}$, where the supremum is taken after all $|a|<\delta$ and $|t|<\delta'$. Thus
\begin{equation}\label{eq:m}
m:=\sup\limits_{\substack{(x,y)\in W^{-},\, |a|<\delta,\, |t|<\delta'}} |\beta_{a,t}(x,y)|
\end{equation}
and so $\big{|}\partial_{x} g_{a,t}(x,y)\big{|}< m|x|^{2q}$ for all $(x,y)\in W^{-}$. The repelling sectors $\Delta^{-}$ and $W^{-}= \Delta^{-}\times \D_{r}$ are defined in Equation \eqref{eq:W-t}.

By eventually reducing the radius $\rho>0$ from the definition of the set $\Delta^{-}$, we can assume that
\begin{equation}\label{eq:eps2}
|1+(q+1)x^{q}|-m|x^{2q}|>1+(q+2/3)\epsilon_{1}|x|^{q}>1,\
\mbox{for all}\ x\in \Delta^{-},
\end{equation}
where $\epsilon_{1}$ is given in Equation \eqref{eq:eps1}. Consider the polynomial
$\widetilde{p}_{t}$ as in Lemma \ref{prop:normal-polynom-t} with its corresponding repelling sector $\Delta^{-}_{R_{t}}$ 
(see Equations \eqref{eq:W-t} and \eqref{eq:W-t2}). The estimate above allows us to show that 
$|\widetilde{p}_{t}\! '(x)|>|\lambda_{t}|(1+(q+2/3)\epsilon_{1}|x|^{q})$ for all $x\in \Delta^{-}_{R_{t}}$.  
The polynomial $\widetilde{p}_{t}\! '$ is clearly expanding if $t$ is nonnegative since 
$|\lambda_{t}|=1+t\geq1$, but Lemma \ref{lem:expansion} shows that it is also expanding on $\Delta^{-}_{R_{t}}$ for negative $t$.

\begin{prop}[{\bf Vertical cones}]\label{prop:cones-Ve}
Consider $(x,y)$ and $(x_1,y_1)$ in the repelling sectors $W^{-}\subset \D_{\rho}\times \D_r$ (respectively in $W^{-}_{R_{t}}$ for $t<0$)
such that $\widetilde{H}(x,y)=(x_1,y_1)$. Then
\[
D\widetilde{H}^{-1}_{(x_1,y_1)}\left(\mathcal{C}^{v}_{(x_1,y_1)}\right)\subset
Int\ \mathcal{C}^{v}_{(x,y)}
\]
and $\big{\|} D\widetilde{H}^{-1}_{(x_{1},y_{1})}(\xi',\eta')\big{\|} \geq (|\nu_{a,t}|+3/2N_{a,t})^{-1}
 \|(\xi',\eta')\| $ for $(\xi',\eta')\in \mathcal{C}^{v}_{(x_1,y_1)}$.
\end{prop}
\proof Let $(\xi',\eta')\in \mathcal{C}^{v}_{(x_1,y_1)}$ with $(\xi',\eta')\neq (0,0)$, and set
$(\xi,\eta)=D\widetilde{H}^{-1}_{(x,y)}(\xi',\eta')$. We need to
show that $(\xi,\eta)\in \mathcal{C}^{v}_{(x,y)}$.  A direct computation gives
\[
D\widetilde{H}_{(x,y)}=\left(
                 \begin{array}{cc}
        \lambda_{t}(1+(q+1)x^q+\partial_{x} g_{a,t} (x,y)) & \lambda_{t}\partial_{y} g_{a,t} (x,y) \\
          \partial_{x} (xh_{a,t})(x,y) & \nu_{a,t}+x\partial_{y} h_{a,t}(x,y) \\
                 \end{array}
               \right)
\]
and so
\begin{eqnarray}
\xi' &=&\lambda_{t}\left(1+(q+1)x^q+\partial_{x} g_{a,t}(x,y)\right)\xi+\lambda_{t}\partial_{y} g_{a,t}(x,y)\eta\label{eq:xixi}\\
\eta' &=& \partial_{x} (x h_{a,t})(x,y)\xi+\left(\nu_{a,t}+\partial_{y} (xh_{a,t})(x,y)\right)\eta.\label{eq:etaeta}
\end{eqnarray}
Using the bounds from Equations \eqref{eq:Nat} and \eqref{eq:Mat} we get
\begin{eqnarray}\label{eq:Ve1}
|\xi'| &\geq& |\lambda_{t}|\left(|1+(q+1)x^q|-m|x|^{2q}\right)|\xi|-
|\lambda_{t}|M_{a,t}|x|^{2q+2}|\eta|\\
|\eta'|&\leq& N_{a,t}|\xi|+\left(|\nu_{a,t}|+N_{a,t}\right)|\eta|.\label{eq:Ve12}
\end{eqnarray}
Since $(\xi',\eta')$ belongs to the vertical cone at $(x_1,y_1)$, we
also know that
\begin{eqnarray*}
\ds
|\xi'|\leq|x_1|^{2q}|\eta'|&\leq&|\lambda_{t}|^{2q}|x|^{2q}|1+x^q+ g_{a,t}(x,y)/x|^{2q}|\eta'|\\
&\leq&|\lambda_{t}|^{2q}|x|^{2q}M_{1}^{2q}|\eta'|,
\end{eqnarray*}
where $M_{1}$ is the supremum of $|1+x^q+ g_{a,t}(x,y)/ x|$ on the
repelling sectors $W^{-}$ of the tubular neighborhood $\D_{\rho}\times \D_{r}$, that is
\begin{equation*}\label{eq:M1}
M_{1}:=\sup\limits_{\substack{(x,y)\in W^{-},\\ |a|<\delta,\, |t|<\delta'}}\big{|}1+x^q+g_{a,t}(x,y)/ x\big{|}.
\end{equation*}
Since $Re(x^{q})>\epsilon_{1}|x|^{q}$ on $W^{-}$ we can take $M_{1}>1$, but any constant $M_{1}>0$ would suffice. 
We have assumed that $|t|<1/(2q)$, so $|\lambda_{t}|^{2q}<3|\lambda_{t}|$ and 
\begin{equation}\label{eq:Ve2}
|\xi'|<3|\lambda_{t}||x|^{2q}M_{1}^{2q}|\eta'|.
\end{equation}

\noindent By combining estimates \eqref{eq:Ve1}, \eqref{eq:Ve12}, and \eqref{eq:Ve2}  we
get
\begin{eqnarray*}
&&|\lambda_{t}|\left(|1+(q+1)x^q|-m|x|^{2q}\right)|\xi|- |\lambda_{t}|M_{a,t}|x|^{2q+2}|\eta|\ \leq\ |\xi'|\ <\ \\\
&&\  <
3|\lambda_{t}||x|^{2q}M_{1}^{2q}|\eta'|\ \leq\  3|\lambda_{t}|M_{1}^{2q}N_{a,t}|x|^{2q}|\xi|+3|\lambda_{t}|M_{1}^{2q}(|\nu_{a,t}|+N_{a,t})|x|^{2q}|\eta|.
\end{eqnarray*}
After regrouping the terms, we write
\[
|\xi|<\frac{A_{2}}{A_{1}}|x|^{2q}|\eta|,
\]
where $A_{1}$ and $A_{2}$ are defined as follows
\begin{eqnarray*}
A_{1} &:=& |1+(q+1)x^q|-(m+3M_{1}^{2q}N_{a,t})|x|^{2q} \\
A_{2} &:=& 3M_{1}^{2q}(|\nu_{a,t}|+N_{a,t})+M_{a,t}|x|^2.
\end{eqnarray*}
 Since $x$ is chosen from the repelling sectors we have
$|1+(q+1)x^q|-m|x|^{2q}>1$.
The quantities $N_{a,t}$, $M_{a,t}$ and $\nu_{a,t}=-a^{2}/\lambda_{t}$  depend
on $a$ and on $t$, and they tend to $0$ as $a\rightarrow0$, uniformly
with respect to $t$. For $|a|$ and $|t|$ small we can therefore assume that
$A_{1}>2/3$ and $A_{2}<1/3$. Hence $(\xi,\eta)\in
\mathcal{C}^v_{(x,y)}$, so
\[
D\widetilde{H}^{-1}_{(x_1,y_1)}\left(\mathcal{C}^{v}_{(x_1,y_1)}\right)\subset
Int\ \mathcal{C}^{v}_{(x,y)}.
\]

We now show that inside the vertical cones the derivative
$D\widetilde{H}^{-1}$ is expanding with respect to the Euclidean
metric. We have
\begin{eqnarray*}
|\eta'|&\leq& N_{a,t}|\xi|+\left(|\nu_{a,t}|+N_{a,t}\right)|\eta| < N_{a,t}\frac{A_{2}}{A_{1}}|x|^{2q}|\eta|+\left(|\nu_{a,t}|+N_{a,t}\right)|\eta| \\
&<&\left(\frac{1}{2}N_{a,t}|x|^{2q}+|\nu_{a,t}|+N_{a,t}\right)|\eta|<\left(|\nu_{a,t}|+\frac{3}{2}N_{a,t}\right)|\eta|,
\end{eqnarray*}
provided that $|x|<1$ (which is already assumed since
$\rho<1$). By definition, as both $(\xi,\eta)$ and $(\xi',\eta')$
are taken from the vertical cones, we have
\begin{equation*}
\| (\xi,\eta)\| = \max(|\xi|,|\eta|) = |\eta| \ \ \ \mbox{and}\ \ \  \| (\xi',\eta')\| = \max(|\xi|',|\eta'|) = |\eta'|.
\end{equation*}

\noindent We obtain $\| (\xi,\eta)\|> (|\nu_{a,t}|+3/2N_{a,t})^{-1}\| (\xi',\eta')\|$.
\qed

For $|a|$ and $|t|$ sufficiently small the expansion factor $(|\nu_{a,t}|+3/2N_{a,t})^{-1}$ can be easily made larger than $1$. Hence
$D\widetilde{H}^{-1}$ expands in the vertical cones with a factor strictly greater than $1$.  

\begin{prop}[{\bf Horizontal cones}]\label{prop:cones-Ho}
Consider $(x,y)$ and $(x_1,y_1)$ in the repelling sectors
$W^{-}\subset \D_{\rho}\times \D_r$ (respectively in $W^{-}_{R_{t}}$ for $t<0$) such that $\widetilde{H}(x,y)=(x_1,y_1)$.
Then
\[
D\widetilde{H}_{(x,y)}\left(\mathcal{C}^{h}_{(x,y)}\right)\subset
Int\ \mathcal{C}^{h}_{(x_{1},y_{1})}
\]
and $\big{\|} D\widetilde{H}_{(x,y)}(\xi,\eta)\big{\|} \geq
|\lambda_{t}|\left(1+(q+1/2)\epsilon_{1}|x|^{q}\right)\|(\xi,\eta)\| $ \ for\ $(\xi,\eta)\in \mathcal{C}^{h}_{(x,y)}$.
\end{prop}
\proof Consider $(\xi,\eta)\in \mathcal{C}^{h}_{(x,y)}$, $(\xi,\eta)\neq (0,0)$, and let
$(\xi',\eta')=D\widetilde{H}_{(x,y)}(\xi,\eta)$. We first need to
show that $(\xi',\eta')\in \mathcal{C}^{h}_{(x_{1},y_{1})}$. Consider $\xi'$ and $\eta'$ written as in Equations \eqref{eq:xixi} and \eqref{eq:etaeta}, from 
the proof of the previous proposition. Since $(\xi,\eta)$ belongs to the horizontal cone at $(x,y)$, we know that $|\xi|\geq|\eta|$. As before, by using Equations \eqref{eq:Mat}, \eqref{eq:Nat}, and \eqref{eq:m}, we get the following estimates
\begin{eqnarray}
|\xi'| &\geq& |\lambda_{t}|\left(|1+(q+1)x^q|-m|x|^{2q}\right)|\xi|-
|\lambda_{t}|M_{a,t}|x|^{2q+2}|\eta| \nonumber \\
&\geq&   |\lambda_{t}|\left(|1+(q+1)x^q|-m|x|^{2q}-M_{a,t}|x|^{2q+2}\right)|\xi| \label{eq:Ho1}\\
|\eta'|&\leq& N_{a,t}|\xi|+\left(|\nu_{a,t}|+N_{a,t}\right)|\eta|\leq (2N_{a,t}+|\nu_{a,t}|)|\xi| \label{eq:Ho2} \nonumber
\end{eqnarray}

\noindent In the final analysis we obtain
\[
|\eta'|\leq\frac{B_{2}}{B_{1}}|\xi'|,
\]
where $B_{1}$ and $B_{2}$ are defined in the obvious way
\begin{eqnarray*}
B_{2} &:=& 2N_{a,t}+|\nu_{a,t}| \\
B_{1} &:=&
|\lambda_{t}|\left(|1+(q+1)x^q|-m|x|^{2q}-M_{a,t}|x|^{2q+2}\right).
\end{eqnarray*}

\noindent The bounds $N_{a,t}$, $M_{a,t}$ and $|\nu_{a,t}|$ tend to $0$ as
$a\rightarrow0$, uniformly with respect to $t$, so one can assume
that for $|a|$ and $|t|$ small enough we have $B_{2}<1/2$.  Moreover, using the bound from  Equation \eqref{eq:eps2}, we can assume that 
\begin{equation}\label{eq: expansion-h}
B_{1}>|\lambda_{t}|\left(1+(q+1/2)\epsilon_{1}|x|^{q}\right).
\end{equation}
Clearly, $B_{1}$ is bounded below by $|\lambda_{t}|>1-1/(2q)$.
In conclusion, we get  
\[
|\eta'|<\frac{q}{2q-1}|\xi'|,
\] 
which implies that $(\xi',\eta')\in Int\ \mathcal{C}^{h}_{(x',y')}$. The norm of the 
two vectors from the horizontal cones are $\|(\xi',\eta')\|=\max(|\xi'|,|\eta'|)=|\xi'|$ and 
$\|(\xi,\eta)\| = \max(|\xi|,|\eta|)=|\xi|$. We have already shown in Equation \eqref{eq:Ho1} that $|\xi'|\geq B_{1}|\xi|$. Together with the lower bound on $B_{1}$ from Equation \eqref{eq: expansion-h} this yields
\begin{eqnarray*}
\|(\xi',\eta')\|> |\lambda_{t}|\left(1+(q+1/2)\epsilon_{1})|x|^{q}\right)|\xi| = |\lambda_{t}|\left(1+(q+1/2)\epsilon_{1})|x|^{q}\right)\|(\xi,\eta)\|,
\end{eqnarray*}
which is what we needed to prove.
\qed

We now analyze the expansion factor $|\lambda_{t}|\left(1+(q+1/2)\epsilon_{1}|x|^q\right)$ in the horizontal cones from Proposition \ref{prop:cones-Ho}. If $t=0$, then $|\lambda_t|=1$ and the expansion factor reduces to $1+(q+1/2)\epsilon_{1}|x|^q$. In this case $D\widetilde{H}$ expands strictly, but not strongly in the horizontal cones. The expansion factor goes to $1$ when $x\rightarrow 0$, i.e. when we approach the local strong stable manifold of the semi-parabolic fixed point.

If $t$ is positive then $|\lambda_t|>1$ and $D\widetilde{H}$ expands strongly in the horizontal cones, by a factor of 
$(1+t)\left(1+(q+1/2)\epsilon_{1}|x|^q\right)\geq (1+t)$.
 If $t$ is negative, then $|\lambda_{t}|<1$ and we need to use the definition of the repelling sector $W^{-}_{R_{t}}$ to get a good expansion. The repelling sectors were carefully defined in Equation \eqref{eq:W-t2}, from the previous section. We use the fact that $|x|^{q}>R_{t}$ for the choice of $R_{t}$ from Equation \eqref{eq:Rt} to make the product $|\lambda_{t}|\cdot\left(1+(q+1/2)\epsilon_{1}|x|^q\right)$ strictly greater than $1$ throughout $\Delta^{-}_{R_{t}}$. So we use the particular choice of $R_{t}$ to make the second term dominate  $|\lambda_{t}|$, which is in fact smaller than 1. The following technical lemma deals with this situation. 

\begin{lemma}[\textbf{Expansion estimate}]\label{lem:expansion} If $t\in(-\delta',0)$, then 
\begin{equation*}\label{eq:eps12}
|\lambda_{t}|\left(1+(q+1/2)\epsilon_{1}|x|^{q}\right) > (1+\epsilon_2 |t|)\left(1+\frac{\epsilon_{1}}{16}|x|^{q}\right),
\end{equation*}
for all $x\in \Delta_{R_{t}}^{-}$, where $\epsilon_2 :=\frac{1}{16(q+1)}$. The inequality is also true for $t\in[0,\delta')$ and $x\in \Delta^-$.
\end{lemma}
\proof The proof is straightforward if $t$ is nonnegative. Suppose that $t$ is negative. Note that $|\lambda_{t}|=1-|t|$. 
We first show that for all $x\in \Delta_{R_{t}}^{-}$
\begin{equation*}\label{eq:18}
 |\lambda_{t}|\left(1+(q+1/2)\epsilon_{1}|x|^{q}\right) = (1-|t|)\left( 1+(q+1/2)\epsilon_{1}|x|^{q}\right)>1+\frac{\epsilon_1}{8}|x|^q,
\end{equation*}
which is equivalent to showing that $|x|^{q}\left((1-|t|)(q+1/2)\epsilon_{1}-\epsilon_{1}/8\right)>|t|$. 
On $\Delta^{-}_{R_{t}}$ we have that $|x|^{q}>R_{t}$ for $R_{t}=\frac{|t|}{(q+1/3)\epsilon_{1}}$, so 
\[
|x|^{q}\left((1-|t|)(q+1/2)\epsilon_{1}-\epsilon_{1}/8\right)>\frac{|t|}{(q+1/3)}\left((1-|t|)(q+1/2)-1/8\right)>|t|.
\]
This is verified for $|t|<\frac{1}{24q+12}$, which is one of the bounds already imposed on $t$. 

We then show by direct computation that 
\[ 
1+\frac{\epsilon_1}{8}|x|^q> (1+\epsilon_2 |t|)\left(1+\frac{\epsilon_{1}}{16}|x|^{q}\right),
\]
for all $x\in \Delta_{R_{t}}^{-}$ and some constant $\epsilon_2$. We take $\epsilon_2=\frac{1}{16(q+1)}$, but the choice is not optimal. The computational details are left to the reader. 
\qed

\subsection{Global analysis of the Julia set}\label{subsec:nbdV-t}

We would first like to show that the corresponding \He map is
hyperbolic on its Julia set $J_{a,t}$. We must show that the derivative of the \He map has appropriate contraction and
expansion in a family of vertical, respectively horizontal cones. We have already shown this to be true locally around the fixed point $\fq_{a,t}$ in Section \ref{section:parabolic2hyperbolic}.

If we look at the form of the \He map,
$H_{a,t}(x,y)=(p_t(x)+a^2w+ay,ax)$ given in Equation \eqref{eq: Henon-form},
we notice that the presence of the multiplicative factor $a$ in the second coordinate implies that the derivative of the \He map is strongly contracting in the ``vertical direction", or
equivalently, $DH^{-1}$ is expanding in the ``vertical direction". If
we analyze the first coordinate, we notice that the expanding
properties of $DH$ in the horizontal direction are closely related
to the expanding properties of the polynomial $p_t$ on a
neighborhood of its Julia set. We will construct a neighborhood $V$
of $J^{+}_{a,t}$ for the \He map $H_{a,t}$ inside a polydisk
$\D_{r}\times\D_{r}$, and put a metric on it with respect to which
the derivative of the \He map is expanding in the ``horizontal
direction".

For the \He map $H_{a,t}$, the construction of the neighborhood $V$
will be similar to the construction of the neighborhood $U_t'$ in
the polynomial case in Section \ref{sec:metrics1D} (see also \cite[Section~7]{RT} for the construction of the neighborhood in the semi-parabolic case).

Let $\alpha_{t}$ be a fixed point for the polynomial $p_t$. 
For $|a|<\delta$
consider the normalizing coordinates of the \He map $H_{a,t}$ on
the tubular neighborhood $B=\D_{\rho'}(\alpha_{t})\times \D_{r}$ as
defined in Theorem \ref{thm:UniformNF}. Let $\fq_{a,t}$ denote the
hyperbolic/semi-parabolic/attracting fixed point. 
Let $W^+_B$ and $W^-_B$ be the attractive, respectively repelling sectors inside $B$ from Definition \ref{defn:WB}. 
By Proposition \ref{prop:localdyn}, the set $W^+_B$ belongs to $int(K_{a,t}^{+})\cup W^{ss}_{loc}(\fq_{a,t})$.
The set $H_{a,t}^{-1}(B) \cap \D_{r}\times \D_{r}$ has two connected components, so let us denote by
\begin{equation}\label{eq:B'}
B':= \left(H_{a,t}^{-1}(B) - B\right) \cap \D_{r}\times \D_{r}
\end{equation}
the component disjoint from $B$. Let $W^{+}_{B'}$  be the preimage of the attractive sectors
$W^{+}_{B}$ in $B'$ that is, $W^{+}_{B'} := H^{-1}_{a,t}(W^{+}_{B})\cap B'$.

We start by defining a box neighborhood $U_t'\times \D_r$,
where $U_t'$ is constructed as in the one dimensional case (see
Definition \ref{def: U_t}). Recall that the set $U_t'$ was defined as
$U_t'=p_t^{-1}(U_t)$, where
\begin{equation*}
U_t: = \C - \overline{p_t^{-\circ N}(S_{att})}-\{ z\in \C-K_{p_{t}} :
|\Psi_{p_t}^{-1}(z)|\geq R\}.
\end{equation*}
We choose $R>2$ large enough so that the outer boundary of $U_t$ is in the set $V^{+}$ defined in \eqref{eq:V+} 
(and implicitly in the escaping set $U^{+}$).

The only difference will be that instead of removing the
attractive sectors $p_t^{-\circ N}(\Delta_t^+)$, we want to remove a
little bit less.
Construct  attractive sectors
$S_{att}\subset \Delta_t^+$ associated with the polynomial $p_t$ in
$\D_{\rho'}(\alpha_t)$, thin enough along the {\it attractive axes}
so that
\begin{equation*}\label{eq:Delta-}
\left( p_t^{ -\circ (N+1)}(S_{att})\cap A\right)\times
\D_{r}\subset W^{+}_{B}.
\end{equation*}
We denoted by $A$ the annulus between the
disk of radius $\rho'$ and the disk of radius $\rho''<\rho'/2$ centered at
$\alpha_{t}$, so  $A= \D_{\rho'}(\alpha_t) - \D_{\rho''}(\alpha_t)$.
Otherwise said, in the annular region $A$, we want the
small attractive sectors $p_t^{-\circ (N+1)}(S_{att})$ of the
polynomial $p_t$ to be compactly contained in the attractive sectors
$W^{+}_{B}$ of the \He map. As in Section \ref{sec:metrics1D}, when writing
$p_t^{-\circ N}(S_{att})$ we do not take into account all preimages of
$S_{att}$, but rather only the preimage of $S_{att}$ that is
contained in the immediate Fatou components of the
fixed point $\alpha_t$ and has $\alpha_t$ in the boundary.

Inside the tubes $B$ and $B'$, we forget all together about the
polynomial dynamics. So we take out the tubes completely and put
back in only the repelling sectors $W^{-}_{B}$ and  $W^{-}_{B'}$.
We can now finally define the set $V$ as
\begin{equation}\label{defn: V}
V := \left(U_t'\times \D_{r}-(B\cup B')\right)\cup
\left(W^{-}_{B}\cup W^{-}_{B'}\right).
\end{equation}

\begin{remarka}\label{remark:V}
Note that the set $B'$ (and consequently $W^{-}_{B'}$) is contained in the bigger set $U_t \times \D_r$ and its projection on the first coordinate is compactly contained in $U_t$ and bounded away from the critical point $0$ of the polynomial $p_t$. Denote by $B''$  the polydisk $\D_{\rho''}(\alpha_t)\times \D_r$.  When $|a|$ is small, the set $B''$ should be thought of as a small neighborhood
 of the local stable manifold
$W^{ss}_{loc}(\fq_{a,t})$ of the fixed point
$\fq_{a,t}$.  By the construction above, the set difference $W^-_B-B''$ is contained in $U_t'\times \D_r$. Hence $V-B''$ is a subset of $U_t\times\D_r$.
\end{remarka}

For $t\geq0$, let $\overline{V}$ denote the set $V$ together with the local stable manifold $W^{ss}_{loc}(\fq_{a,t})$  and together with $H^{-1}(W^{ss}_{loc}(\fq_{a,t}))\cap B'$. When $t<0$ there is no need to add the two stable manifolds as they belong to the interior of $K^{+}_{a,t}$. However, to preserve notation, we set $\overline{V}=V$ in this case. 

\begin{lemma}\label{lemma:nbdPolydisk} $J_{a,t}^{+}\cap \overline{V}= J_{a,t}^{+}\cap \D_r\times\D_r$.
Moreover, the Julia set $J_{a,t}$ is contained in $\overline{V}$ and
\begin{equation*}
    J_{a,t} = \bigcap\limits_{n\geq 0}H^{\circ n}(J^{+}_{a,t}\cap\overline{V}).
\end{equation*}
\end{lemma}
\proof The outer boundary of the set $V$ is an equipotential of the
polynomial $p_t$ cross $\D_{r}$, which belongs to the escaping set $U_{a,t}^{+}$.
From the tubular neighborhood $B$ of the local stable manifold we
have removed only the attractive sectors $W^{+}_{B}$, which are
contained in the interior of $K_{a,t}^{+}$ when $t<0$ and in the interior of $K_{a,t}^{+}$ union the local
stable manifold $W^{ss}_{loc}(\fq_{a,t})$ when $t\geq 0$.  From $B'$ we only removed
the attractive sectors $W^{+}_{B'}$ which are contained in the
interior of $K^{+}_{a,t}$ when $t<0$, and respectively in the
interior of $K^{+}_{a,t}$ union a preimage of the local stable manifold
$H^{-1}(W^{ss}_{loc}(\fq_{a,t})) \cap B'$ when $t\geq0$. Outside of $B\cup B'$, we
have removed a vertical tube $p_t^{-\circ (N+1)}(S_{att})\times \D_{r}$
which belongs to the interior of $K^{+}_{a,t}$. Therefore, when $t<0$, the set $J_{a,t}^{+}\cap \left(\D_{r}\times\D_{r}\right)$ is contained in the set $V$. When $t\geq 0$, in our
construction process of the neighborhood $V$, we have lost from
$J_{a,t}^{+}\cap \left(\D_{r}\times\D_{r}\right)$ only two local
stable manifolds. These local stable manifolds are no longer in $V$,
but they lie in the larger set $\overline{V}$,
\[
J_{a,t}^{+}\cap \left(\D_{r}\times\D_{r}\right)\  \subset
\overline{V}\ \ \mbox{when}\ t\geq 0.
\]
Any point in $J_{a,t}^+ \cap \overline{V}$ remains in $\overline{V}$
under forward iterates of the \He map, so the Julia set $J_{a,t}$ is
contained in $\overline{V}$. \qed

\subsection{Vertical and horizontal cones in the product metric}\label{subsec:hyper-region}

We  construct an invariant family of horizontal and vertical cones on the set $V$ defined in \eqref{defn: V},
such that the derivative of the \He map expands in the horizontal cones, and
contracts in the vertical cones.

In Section \ref{section:parabolic2hyperbolic} we have already constructed such an invariant family of cones in the repelling sectors 
of the fixed point $\fq_{a,t}$. These cones live only in a small neighborhood of $\fq_{a,t}$, where the map is conjugate to the normal form \eqref{prop:perturbed-form}.

In this section we define a family of cones on the set V away from a small neighborhood of $\fq_{a,t}$.  
At the end of this section, we show how to patch together these two types of cones, from Sections \ref{section:parabolic2hyperbolic} and \ref{subsec:hyper-region}, to get an invariant family on the entire set $V$.

The set $U_t\times \D_r$ comes equipped with the product metric
$\mu_{U_t}\times \mu_E$ of the Poincar\'{e} metric $\mu_{U_t}$ of the set $U_t$
and the regular Euclidean metric $\mu_E$ on the vertical disk $\D_r$. 
Tangent vectors $(\xi,\eta)$ from $T_{(x,y)}\C^2$ will be
measured with respect to the product metric
\[
\| (\xi,\eta)\| :=\max(\mu_{U_t}(x,\xi),|\eta|),
\]
where $|\eta|$ is the absolute value of the complex number $\eta$.

By Remark \ref{remark:V}, the set $V-B''$ is a subset of $U_t\times\D_r$, and we can endow $V-B''$ with the product metric that we have just constructed on the set $U_t\times\D_r$. 
Denote by $U_t''$ the projection of $V$ on the first coordinate, which is equal to $U_t' \cup pr_1(W^-_{B'})$. 
The set $U_t'' - \D_{\rho''}(\alpha_t)$ is compactly contained in $U_t$, so the Poicar\'{e} metric
$\mu_{U_t}$ is bounded above and below by the Euclidean metric on the set $U_t'' - \D_{\rho''}(\alpha_t)$, that is, there exist two positive constants
$m_1$ and $m_2$ such that
\begin{eqnarray}\label{eq:m1m2}
m_1 < \rho_{U_t}(x)<m_2,
\end{eqnarray}
for any $t\in [-\delta',\delta']$ and any $x\in
U_t''-\D_{\rho''}(\alpha_t)$. Therefore, the product metric is bounded on the set $V-B''$. If
we let $\rho_{U_t}$ be the density function of the metric
$\mu_{U_t}$,
\[
\mu_{U_t}(x,\xi)=\rho_{U_t}(x)|\xi|
\]
then $\rho_{U_t}$ is positive and $C^{\infty}$-smooth on
$U_t''-\D_{\rho''}(\alpha_t)$.

The sets $U_t'',\
|t|<\delta'$, avoid a neighborhood of fixed size of the critical
point of the polynomial $p_t$. Hence there exists a lower bound
$r_1>0$ such that 
\begin{equation}\label{eq:r1}
r_1<|p_t'(x)|, \mbox{ for any } x\in U_t''.
\end{equation}

\begin{defn}\label{def:vcones-poincare}
Let $\tau<1$. Let the \textbf{vertical cone} at a point $(x,y)$
from $V-B''$ be
\[
\mathcal{C}^{v}_{(x,y)}=\left\{(\xi,\eta)\in
T_{(x,y)}\C^2,\ \mu_{U_t}(x,\xi)\leq\tau\cdot
|\eta|\right\}.
\]
Define the \textbf{horizontal cone} at a point $(x,y)$ from the set
$U'_t\times \D_r$  to be
\[
\mathcal{C}^{h}_{(x,y)}=\left\{(\xi,\eta)\in
T_{(x,y)}\C^2,\ \mu_{U_t}(x,\xi)\geq|\eta|\right\}.
\]
\end{defn}

We will show that the vertical cones are invariant under
$DH_{a,t}^{-1}$ and that the horizontal cones are invariant under
$DH_{a,t}$.

\begin{prop}[{\bf Vertical cones}]\label{prop:cones-Ve-product-metric}
Consider $(x,y)$ and $(x',y')$ in $V-B''$ such that $H(x',y')=(x,y)$.
Then
\[
DH^{-1}_{(x,y)}\left(\mathcal{C}^{v}_{(x,y)}\right)\subset Int\
\mathcal{C}^{v}_{(x',y')}
\]
and $ \big{\|} DH^{-1}_{(x,y)}(\xi,\eta)\big{\|} \geq |a|^{-1} \|(\xi,\eta)\| $ 
for $(\xi,\eta)\in \mathcal{C}^{v}_{(x,y)}$.
\end{prop}
\proof
Let $(\xi,\eta)\in \mathcal{C}^{v}_{(x,y)}$  and 
$(\xi',\eta') = DH^{-1}_{(x,y)}(\xi,\eta)$.
From the formula of the inverse of the \He map
\[
\hvec{x'}{y'}=\hvec{\frac{y}{a}}{ \frac{x-p_t(y/a)-a^2w}{a} }
\]
we find $ \xi'=\frac{1}{a}\eta$ and $ \eta' = \frac{1}{a}\left(\xi-\frac{2x'}{a} \eta\right)$. Assume $(\xi,\eta)\neq(0,0)$, otherwise the proof is trivial. The vector
$(\xi,\eta)$ belongs to the vertical cone, so $ \mu_{U_t}(x,\xi)=\rho_{U_t}(x)|\xi|\leq \tau  |\eta|$. This
implies that 
\begin{equation}\label{eq: cone-v}
 |\xi|<\frac{\tau}{m_1} |\eta|.
\end{equation}
We can evaluate
\begin{equation}\label{eq: lh-1}
 \mu_{U_t}(x',\xi')=\rho_{U_t}(x') \frac{|\eta|}{|a|}\leq \frac{m_2}{|a|} |\eta|.
\end{equation}
Next, by using inequality \eqref{eq: cone-v}, we compute 
\begin{equation}\label{eq: rh-1}
 |\eta'|=\frac{1}{|a|} \left|\xi-\frac{2x'}{a}
\eta\right|>\frac{1}{|a|}  \left(\frac{|2x'|}{|a|} -
\frac{\tau}{m_1}\right) |\eta|.
\end{equation}
The point $x'$ belongs to $U_t''$, so $|2x'|>r_1$ by Equation \eqref{eq:r1}. 
Choose $|a|$ small so that $
\frac{r_1}{|a|}-\frac{\tau}{m_1} > \max\left(\frac{2 m_2}{\tau}, 1\right)$. Combining Equations \eqref{eq: lh-1} and
\eqref{eq: rh-1} gives 
$
 \mu_{U_t}(x',\xi') < \frac{\tau}{2} |\eta|.
$
Therefore $ DH^{-1}_{(x,y)}\left(\mathcal{C}^{v}_{(x,y)}\right)\subset Int\ \mathcal{C}^{v}_{(x',y')}$, which proves the cone invariance. Inequality \eqref{eq: rh-1} shows that
$DH^{-1}$ expands in the vertical cone as $ |\eta'|>|a|^{-1} |\eta|$.
By definition, since both $(\xi,\eta)$ and $(\xi',\eta')$ belong to
vertical cones, we have
\begin{equation*}
\| (\xi,\eta)\| = \max\left(\mu_{U_t}(x,\xi),|\eta|\right) = |\eta| \ \ \ \mbox{and}\ \ \
\| (\xi',\eta')\| = \max\left(\mu_{U_t}(x',\xi'),|\eta'|\right) = |\eta'|.
\end{equation*}
\noindent We therefore obtain $\| (\xi,\eta)\|> |a|^{-1}
\| (\xi',\eta')\|$, as claimed. 
\qed

\begin{remarka} When dealing with vertical cones, it is not really
necessary to measure the horizontal component of vectors with respect to the
Poincar\'e metric. Any bounded metric in the horizontal direction
would work, because we can always choose $|a|$ small to get the
invariance of the vertical cone field and the strong expansion of $DH^{-1}$ in
the vertical cones. The choice of the Poincar\'e metric is essential
however to show expansion of $DH$ in the horizontal cones.

The scalar $0<\tau<1$ in the definition of the vertical cone will
typically be chosen less than $\left(\rho/2\right)^{2q}$, so
that on a neighborhood of the boundary of $B$, the vertical cones
$\mathcal{C}^{v}_{(x,y)}$ from Definition \ref{def:vcones-poincare}
are contained in the pull-back by $D\phi_{a,t}$ of the vertical
cones from \ref{def:vcones} defined in the normalized coordinates.
\end{remarka}

\begin{prop}[{\bf Horizontal cones}]\label{prop:cones-Ho-product-metric}
Let $(x,y)$ and $(x',y')$ in
$V-B''$ such that
$H(x,y)=(x',y')$. Then we have
\[
DH_{(x,y)}\left(\mathcal{C}^{h}_{(x,y)}\right)\subset Int\
\mathcal{C}^{h}_{(x',y')}
\]
and $\big{\|} DH_{(x,y)}(\xi,\eta)\big{\|} \geq k \|(\xi,\eta)\| \ $ for $\ (\xi,\eta)\in \mathcal{C}^{h}_{(x,y)}$.

\end{prop}
\proof Let $(\xi,\eta)\in \mathcal{C}^{h}_{(x,y)}$ with  $(\xi,\eta)\neq (0,0)$ and let
$(\xi',\eta')=DH_{(x,y)}(\xi,\eta)$. We first need to show that
$(\xi',\eta')\in Int\ \mathcal{C}^{h}_{(x',y')}$. Since $DH_{(x,y)}(\xi,\eta)=(2x\xi+a\eta,a\xi)$, we find 
$ \xi'=2x\xi+a\eta$ and $\eta' = a\xi$.
The vector $(\xi,\eta)$ belongs to the horizontal cone at $(x,y)$, so
\begin{equation}\label{eq: cone-h}
|\eta|\leq\mu_{U_t}(x,\xi)=\rho_{U_t}(x)|\xi|<m_2|\xi|.
\end{equation}
Using $\hvec{x'}{y'}=\hvec{p_t(x)+a^2w+ay}{ax}$ we evaluate
\begin{equation}\label{eq: lh-2}
\mu_{U_t}(x',\xi')=\rho_{U_t}\left(p_t(x)+a^2w+ay\right)|2x\xi+a\eta|.
\end{equation}
Since $\rho_{U_t}$ is $C^{\infty}$-smooth, its derivative
$\rho'_{U_t}$ is also bounded on $U_t''-\D_{\rho''}(\alpha_t)$. There exists a constant $c>0$ (which just a local variable) such
that
\begin{eqnarray}\label{eq: density-smooth}
\frac{\left|\rho_{U_t}\left(p_t(x)+a^2w+ay\right)-\rho_{U_t}\left(p_t(x)\right)\right|}{|a|}
&\leq& |aw+y|\cdot \sup\limits_{} \rho'_{U_t}\cdot
\frac{\rho_{U_t}(p_t(x))}{\inf\limits_{} \rho_{U_t}} \nonumber
\\
&<& c\cdot \rho_{U_t}\left(p_t(x)\right).
\end{eqnarray}
\noindent The polynomial $p_t$ is expanding  with respect
to the Poincar\'e metric $\mu_{U_t}$. As in Lemma \eqref{eq:
c_t-polynom} part (a), there exists $\kappa_t>1$ with $\inf\limits_{t\in
[-\delta',\delta']}\kappa_t>1$ such that
\begin{equation}\label{eq: c_t-polynom2}
\ds \rho_{U_t}\left(p_t(x)\right)|p_t'(x)\xi|>\kappa_t\cdot \rho_{U_t}(x)
|\xi|,
\end{equation}
whenever $x,p_t(x)\in U_t''-\D_{\rho''}(\alpha_t)$. 
We
now turn back to relation \eqref{eq: lh-2}. Using \eqref{eq:
density-smooth}, \eqref{eq: c_t-polynom2}, \eqref{eq: cone-h} and \eqref{eq:r1} one
gets
\begin{eqnarray}\label{eq: rh-2}\ds
\mu_{U_t}(x',\xi') &>&(1-c|a|)\cdot
\rho_{U_t}\left(p_t(x)\right)|2x\xi|\cdot\frac{|2x\xi+a\eta|}{|2x\xi|}
\nonumber \\
&>& (1-c|a|)\cdot \kappa_t \cdot \rho_{U_t}(x)|\xi| \cdot
\left(1-|a|\frac{|\eta|}{|2x||\xi|}\right) \nonumber \\
&>& \kappa_t \cdot (1-c|a|)\cdot \left(1-|a|\frac{m_2}{r_1}\right)\cdot
\rho_{U_t}(x)|\xi|.
\end{eqnarray}
The constant $\kappa_t$ is bigger than $1$ for all $t\in [-\delta',\delta']$. We
write the dependence on $t$ to preserve notations from Lemma
\eqref{eq: c_t-polynom}, but we could drop the dependence on $t$ by
working with $\inf\limits_{t\in [-\delta',\delta']}\kappa_t$ which is also
strictly bigger than $1$. The factors $1-c|a|$ and
$1-|a|\frac{m_2}{r_1}$ are independent of $t$, and they can be made
arbitrarily close to $1$ by reducing $|a|$. In conclusion, for $|a|$
sufficiently small we can assume that
\begin{equation}\label{eq:k}
k := \inf_{t\in [-\delta',\delta']} \kappa_t \cdot (1-c|a|)\cdot
\left(1-|a|\frac{m_2}{r_1}\right) >1.
\end{equation}
From relation \eqref{eq: rh-2}, we obtain
\begin{equation*}
\mu_{U_t}(x',\xi') > k \cdot \rho_{U_t}(x)|\xi|=k\cdot
\mu_{U_t}(x,\xi),
\end{equation*}
which shows that $DH$ expands in the horizontal cones. Also from
 \eqref{eq: rh-2} we infer that
\begin{equation*}
|a|\cdot \mu_{U_t}(x',\xi') > k \cdot \rho_{U_t}(x)|a\xi|> k \cdot
m_1 \cdot |\eta'|,
\end{equation*}
which proves that $DH_{(x,y)}\left(\mathcal{C}^{h}_{(x,y)}\right)\subset
Int\ \mathcal{C}^{h}_{(x',y')}$, so the horizontal cones are invariant. \qed

On the set $V-B''$ we have one family of horizontal/vertical cones, $\mathcal{C}^{h}_{(x,y)}$ and $\mathcal{C}^{v}_{(x,y)}$, defined in \ref{def:vcones-poincare}. On $W^-_B$ we have another family of horizontal/vertical cones
\[
D\phi_{a,t}^{-1}|_{\phi_{a,t}(x,y)}(\mathcal{C}^{h}_{\phi_{a,t}(x,y)}) 
\ \mbox{ and } \ 
D\phi_{a,t}^{-1}|_{\phi_{a,t}(x,y)}(\mathcal{C}^{v}_{\phi_{a,t}(x,y)})
\]
defined in \ref{def:vcones} with respect to the Euclidean metric in the normalized coordinates given by $\phi_{a,t}$.
For those points $(x,y)\in W^-_B$ where both types of cones are defined, we take the horizontal/vertical cone to be their intersection.

\subsection{Combining infinitesimal metrics}\label{subsec: combined}

On the neighborhood $V$ defined in \ref{defn: V} we have given two
infinitesimal metrics. On the set $V-B''$, where $B''\subset B$ was defined in Remark \ref{remark:V},
we put the product of the
Poincar\'e metric $\mu_{U_t}$ with the Euclidean metric on $\D_r$,
\begin{equation}
\mu_P((x,y),(\xi,\eta)):= \max \left( \mu_{U_t}(x,\xi), |\eta| \right),
\end{equation}
where $(x,y)\in V-B''$ and $(\xi,\eta)\in
T_{(x,y)}V-B''$.

In the repelling sectors $W^-_B$ of the tubular neighborhood $B$ of
the local  strong stable manifold of the hyperbolic/semi-parabolic/attractive fixed
point (see Definition \ref{defn:WB}), we have the pull-back Euclidean metric from the
normalizing coordinates $\phi_{a,t}:W^-_B\rightarrow W^-\subset
\D_{\rho}\times\D_r$. Let
\begin{equation*}
\mu_B((x,y),(\xi,\eta)):= \max\left( |\widetilde{\xi}|,|\widetilde{\eta}|\right),
\end{equation*}
where
$(\widetilde{\xi},\widetilde{\eta})=D\phi_{a,t}\big|_{(x,y)}\left(\xi,\eta\right)$
and $\phi_{a,t}:B\rightarrow \D_{\rho}\times \D_r$ is the change of
coordinate function from Theorem \ref{thm:UniformNF}.

Just like in the polynomial case, we can define an infinitesimal pseudo-norm on
$V$,
\begin{equation}\label{defn: mu-metric}
\mu:=\inf\left(M \mu_B, \mu_P\right),
\end{equation}
where $M$ is a positive real number, chosen so that  the derivative of the \He
map is still expanding in the horizontal cones when we map from the repelling sectors $W^{-}_{B'}$ of $B'$ (see Definition \eqref{eq:B'})
into the repelling sectors $W^{-}_B$ of $B$. We take $M$ so that
\begin{equation}\label{defn: M}
M > \sup\limits_{\substack{(x,y)\in W^{-}_{B'}\\
(\xi,\eta)\in \mathcal{C}^h_{(x,y)} - \{(0,0)\} }}
\small{\frac{2 \mu_{P}\left((x,y),(\xi,\eta)\right)}{
\mu_{B}\left(H(x,y),DH_{(x,y)}(\xi,\eta)\right)}}.
\end{equation}
The supremum from Equation \eqref{defn: M} is bounded, because $\mu_P((x,y),(\xi,\eta))= \mu_{U_t}(x,\xi)$ for $(x,y)\in W^{-}_{B'}$, $(\xi,\eta)\in \mathcal{C}^h_{(x,y)}$, and  the Poincar\'e metric on $B'$ is bounded.

Let $(x,y)$ be any point in $V$. Let $(\xi,\eta)$, $(\xi',\eta')$ be any vectors from the two dimensional tangent space $T_{(x,y)}V$.

Notice that $\mu$ is homogeneous that is, $\mu((x,y), \alpha (\xi, \eta))= |\alpha| \mu((x,y), (\xi, \eta))$, for all complex numbers $\alpha$, as both $\mu_P$ and $\mu_B$ are homogeneous metrics.  It also satisfies the relation $\mu((x,y),(\xi,\eta))\geq 0$, with equality if and only if $(\xi,\eta)=(0,0)$. However, $\mu$ does not necessarily satisfy the triangle inequality $\mu((x,y), (\xi+\xi', \eta+\eta'))\leq \mu((x,y), (\xi, \eta))+\mu((x,y), (\xi', \eta'))$. Nonetheless, $\mu$ induces a regular path metric on horizontal curves between points in $V$ (see Definitions $\eqref{eq:l-pol}$ and $\eqref{eq:d-pol}$) by integration. If
$g:[0,1]\rightarrow V$ is a horizontal rectifiable path $g(s)=(g_1(s),y)$, then its
length with respect to $\mu$ is given by the formula
$
\ell_{\mu}(g)=\int_0^1 \mu \left(g(s),\left(g_1'(s),0\right)\right) ds.
$
The distance between two points $(x,y)$ and $(x',y)$ from $V$ with
respect to the induced metric $\mu$ is
\begin{equation}\label{eq:d-pol2henon}
d_{\mu}\left((x,y),(x',y)\right)=\inf \ell_{\mu}(g),
\end{equation}
where the infimum is taken after all horizontal rectifiable paths
$g:[0,1]\rightarrow V$ with $pr_2(g)=y$,  $g(0)=(x,y)$ and
$g(1)=(x',y)$.

Another way to combine the two metrics is by defining a true product metric, where the second coordinate is measured with respect with the Euclidean metric, and the first coordinate is an infimum of two metrics. Choose $(x,y)\in U_t'\times \D_r \cap W^-_B$ and nonzero $(\xi,\eta)\in
T_{(x,y)}U_t'\times \D_r \cap W^-_B$ as before. Using Remark \eqref{eq:horcurves}, let $(\skew{3.75}\widetilde{\widetilde{\xi}},0)=D\phi_{a,t}\big|_{(x,y)}\left(\xi,0\right)$ and define
\begin{equation*}
\mu'((x,y),(\xi,\eta)):=\max\left(\inf\left(\mu_{U_t}(x,\xi), M|\skew{3.75}\widetilde{\widetilde{\xi}}|\right), |\eta|\right).
\end{equation*}
The constant $M$ is greater than $\sup\left(2\mu_{U_t}\left(x,\xi\right)/ |\skew{3.75}\widetilde{\widetilde{\xi}}_1 |\right)$, where $(\xi_1,\eta_1)=DH_{(x,y)}(\xi,\eta)$ and the supremum is taken after all $(x,y)\in W^{-}_{B'}$ and nonzero vectors $(\xi,\eta)\in \mathcal{C}^h_{(x,y)}$.

If we let $\phi_{a,t}=(\phi_1, \phi_2)$, then $ \skew{3.75}\widetilde{\widetilde{\xi}}=\partial_x\phi_1(x,y) \xi$. Also $\mu_{U_t}(x,\xi)=\rho_{U_t}(x)\xi$, where $\rho_{U_t}$ is the density function of the Poincar\'e metric of $U_t$. In conclusion, if we define $m(x,y):= \inf(\rho_{U_t}(x), M|\partial_x\phi_1(x,y)|)$ we get
\begin{equation*}
\mu'((x,y),(\xi,\eta)):=\max\left( m(x,y)|\xi|, |\eta|\right).
\end{equation*}

With this definition it is easy to see that the triangle inequality is satisfied and $\mu'$ is an infinitesimal metric, i.e. a norm. Notice also that $\mu$ and $\mu'$ coincide when restricted to horizontal curves, and they induce the same horizontal path metric.

\begin{lemma}\label{lemma: kt}
Let $t\in[-\delta',\delta']$ and $|a|<\delta$. Let $k$ be chosen as in Proposition \ref{prop:cones-Ho-product-metric} and the constant $\epsilon_2$ as in Lemma \ref{lem:expansion}.
There exists a constant
$k_t\geq \min(1+\epsilon_2|t|, k)>1$, independent of $a$, such that
\begin{equation}\label{eq:kt}
\mu\left(H_{a,t}(x,y),DH_{a,t}\big|_{(x,y)}(\xi,\eta)\right)> k_t
\cdot \mu\left((x,y),(\xi,\eta)\right),
\end{equation}
for any $(x,y)\in V$ and any nonzero tangent vector $(\xi,\eta)$ in the
horizontal cone at $(x,y)$.

If $t=0$ and $|a|<\delta$, there exists $k_0(x,y)>1$ such that
\begin{equation}\label{eq:k0}
\mu\left(H_{a,t}(x,y),DH_{a,t}\big|_{(x,y)}(\xi,\eta)\right)> k_0(x,y)
\cdot \mu\left((x,y),(\xi,\eta)\right),
\end{equation}
and $k_0(x,y)$ goes to $1$ precisely when $(x,y)$ tends to the local
stable manifold $W^{ss}_{loc}(\fq_{a,0})$ of the semi-parabolic fixed
point. 

Moreover, the inequalities \eqref{eq:kt} and \eqref{eq:k0} hold true for $\mu'$ instead of $\mu$.
\end{lemma}
\proof The proof is identical to the proof of Lemma \ref{lemma:
vector-1dim}. We use the estimates in the horizontal cones from
Propositions \ref{prop:cones-Ho} and \ref{prop:cones-Ho-product-metric}.
The choice of the constant $M$ in Definition \ref{defn: M} is useful when dealing with
the analogue of case (d)(i), from Lemma \ref{lemma: vector-1dim}. The
case $t=0$ is given by \cite[Theorem~8.7]{RT}. 
\qed

\noindent \textbf{A larger region of hyperbolicity.} The following theorem is a classical result on dominated splitting. 
\begin{thm}[\cite{KH}]\label{thm:KH}
A compact $f$-invariant set $\Lambda$ is hyperbolic if there exists
$\beta_{0}<1<\beta_{1}$ such that for every $x\in \Lambda$ there is
a decomposition $T_{x}M=S_{x}\oplus T_{x}$, a family of horizontal
cones $\mathcal{C}^{h}_{x}\supset S_{x}$, and a family of vertical
cones $\mathcal{C}^{v}_{x}\supset T_{x}$ associated with that
decomposition such that
\[
Df_{x}\mathcal{C}^{h}_{x}\subset Int\ \mathcal{C}^{h}_{f(x)},\ \
Df_{x}^{-1}\mathcal{C}^{v}_{f(x)}\subset Int\ \mathcal{C}^{v}_{x},
\]
$\|Df_{x}\ \xi \|\geq \beta_{1}\|\xi\|$ for $\xi\in
\mathcal{C}^{h}_{x}$ and $\|Df_{x}^{-1}\ \xi\|\geq
\beta_{0}^{-1}\|\xi\|$ for $\xi\in \mathcal{C}^{v}_{x}$.
\end{thm}

We now have all the ingredients to prove the hyperbolicity part of Theorem \ref{thm:HypRegion}.

\begin{thm}[\textbf{Hyperbolicity}]\label{thm:HypRegion-b}
There exist $\delta,\delta'>0$ such that in the parametric region
\begin{eqnarray*}
\mathcal{HR}_{\delta, \delta'}=\left\{(c,a)\in \mathcal{P}_{\lambda_t}\ :  \ 0<|a|<\delta\ \ \mbox{and}\
 -\delta'<t<\delta',\  t\neq 0\right\}
\end{eqnarray*}
the Julia set $J_{c,a}$ is connected and the \He map $H_{c,a}$ is hyperbolic.
\end{thm}
\proof In Section \ref{subsec:hyper-region},  we built a family of horizontal and vertical cones,
invariant under $DH$, respectively under $DH^{-1}$, such that $DH$ expands with a factor of $\beta_1>1$
inside the horizontal cones, and $DH^{-1}$ expands with a factor of $1/\beta_0>1$ inside the vertical cones. The expansion is measured with respect to the metric $\mu'$ from Lemma \ref{lemma: kt}.
The proof follows from Propositions \ref{prop:cones-Ve},
\ref{prop:cones-Ho}, \ref{prop:cones-Ve-product-metric} and \ref{prop:cones-Ho-product-metric} and Lemma \ref{lem:expansion} by taking $\beta_{0}= \max(|a|, |\nu_{a,t}|+3/2N_{a,t})\ll 1$ and $\beta_{1}=k_t>1$ for $t\neq0$. The constant $k_t$ is given in Lemma \ref{lemma: kt}. We then apply Theorem \ref{thm:KH} for the set
$\Lambda:=\overline{V}$, which includes $J^+\cap\D_r\times\D_r$, by Lemma \ref{lemma:nbdPolydisk}. The set $V$ was constructed in Section \ref{subsec:nbdV-t}. 

The fact that this hyperbolic region is inside a
component of the connectedness locus follows from Corollary \ref{cor:connected}.
\qed

It is worth mentioning that 
$J=J^*$ (the closure of the saddle periodic points) throughout the parametric region defined by $(c,a)\in \mathcal{P}_{\lambda_t}$, $0<|a|<\delta$ and $-\delta'<t<\delta'$. This follows from Theorem \ref{thm:HypRegion-b} and \cite{BS1} for $t\neq 0$ and \cite{RT} for $t=0$.

\bigskip

\noindent\textbf{Period doubling.}
Let $\mathcal{P}^{2^{n}}_{-1}$ be the set of parameters $(c,a)\in\C^{2}$ for which the \He map $H_{c,a}$ has a cycle of period $2^{n}$ with one multiplier $\lambda=-1$. Theorem \ref{thm:HypRegion-b} can be generalized to show that there are regions of hyperbolicity to the left and to the right of the real curve $\mathcal{P}_{-1}^{2^{n}}\cap\R^{2}$ (see Figure \ref{fig:perspective}). Moreover, for each $n$, there is a region of hyperbolicity connecting $\mathcal{P}_{-1}^{2^{n}}\cap\R^{2}$ to $\mathcal{P}_{-1}^{2^{n+1}}\cap\R^{2}$ of ``vertical'' size (the size of the parameter $a$) $\delta_{n}>0$. Presumably $\delta_{n}\rightarrow 0$ as we approach the Feigenbaum parameter.

\subsection{The function space $\mathcal{F}$}\label{sec: FunctionSpaceHenon}

We will do the same construction as in Section \ref{sec:contraction}. Let $R$ be fixed as in Equation \eqref{def: U_t}. 
Recall that 
\[
\gamma_{t,0}:\s^1\rightarrow U_t',\ \
 \gamma_{t,0}(s)=\Psi_{p_t}\left(R^{1/2}e^{2\pi i s}\right)
\]
is the equipotential of $p_t$ that defines the outer boundary of the
neighborhood $U_t'$ constructed in Section \ref{sec:metrics1D}. Define $f_0:\s^1\times \D_r \rightarrow V$ as 
\[
f_0(s,z)=(\gamma_{t,0}(s),z).
\]
The image of $f_0$ is thus a solid torus contained in the escaping set $U^+$ that represents the outer boundary of the set $V$.

\begin{defn}\label{def:F}
Consider the space of functions:
\begin{equation*}
\mathcal{F}_{a,t} = \left\{ f_n:\s^{1}\times \D_{r} \rightarrow V :
        f_0(s,z)=(\gamma_{t,0}(s),z),\ f_n(s,z)=F_{a,t}\circ f_{n-1}(s,z)\
        \mbox{for}\ n\geq 1
        \right\}\!,
\end{equation*}
where the graph transform $F_{a,t}:\mathcal{F}_{a,t}\rightarrow\mathcal{F}_{a,t}$ is 
defined as
$ F_{a,t}(f)=\tilde{f}$, where the map $\tilde{f}$ is continuous with respect to $s$, holomorphic with respect to $z$, and $ \tilde{f}\big{|}_{s\times\D_{r}}$ is the reparametrization  $\tilde{f}(s,z)=(\varphi_{s}(z),z)$ of one of the two vertical-like disk components of 
\[
H_{a,t}^{-1}\left(f(2s\times\D_{r}) \right)\cap V
\]
as a graph of a function over the second coordinate, via the Inverse Function Theorem.
\end{defn}

\begin{remarka} The maps $f_{n}$ are essentially reparametrizations of the backward iterates of $f_{0}$ (inside $V$) under the \He map.  The picture to keep in mind is the following: The image of the map $f_n\in
\mathcal{F}_{a,t},\ n\geq 0$, is a solid torus $T_n$ contained in the escaping set $U^+$. In the $s$-coordinate, the \He map behaves like angle doubling, whereas in the vertical $z$-coordinate, it behaves like a strong contraction. Therefore, the \He map maps $T_{n+1}$ to another solid torus, wrapped around two times inside $T_{n}$. 
\end{remarka}

We say that a complex disk is vertical-like if any tangent vector to it belongs to the vertical cones \ref{def:vcones-poincare}  and \ref{def:vcones} if both of them are defined, or to the one that is defined. The invariance of vertical cones (Propositions \ref{prop:cones-Ve-product-metric}  and \ref{prop:cones-Ve}) and the fact that the \He map has degree $2$ imply that the preimage of a vertical-like complex disk contained in the set $V\cap U^+$ consists of two vertical-like complex disks. Assume by induction on $n\geq 0$, that we have $f_n(s,z)=(\varphi^{n}_s(z),z)$, where $f_n$ is injective, continuous with respect to $s$ and holomorphic with respect to $z$, and for any $s\in\s^1$, $L=f_n(2s\times \D_r)$ is a vertical-like disk in the escaping set $U^{+}$. Let us show how to construct $f_{n+1}$. The projection $\Delta$ of $L$ on the first coordinate is almost constant and bounded away from $0$, the critical point of $p_t$, so the preimage $H_{a,t}^{-1}(L)\cap V$ is a disjoint union of two vertical-like disks that we would like to first label as $s$ and $s+1/2$ and then parametrize as graphs over the second coordinate. As in the polynomial case, there are exactly two possible choices of labelings that would make the function $f_{n+1}$ continuous with respect to $s\in \s^1$. 
There are two holomorphic branches 
of the backward iterate of the polynomial $p_t$ defined on $\Delta$.  Let now $(\varphi^{n}_{2s}(z),z)$ be any point of $L$ such that $H_{a,t}^{-1}(\varphi^{n}_{2s}(z),z)\in V$. In particular, by analyzing the second coordinate of $H_{a,t}^{-1}(\varphi^{n}_{2s}(z),z)$, we see that the condition $|(\varphi^{n}_{2s}(z)-p(z/a)-a^{2}w)/a|<r$ must be satisfied, which means exactly that the first coordinate $z/a$ is $\bigO(a)$ close to one of the two preimages of  $\varphi^{n}_{2s}(z)$ under the polynomial $p_t$. The curves $f_{n+1}(s\times \D_r)$ and $f_{n+1}((s+1/2)\times \D_r)$ correspond to different choices of the branch of $p_t^{-1}$ (see also Section \ref{sec:contraction}, Equation \ref{def: p-1}). The Inverse Function Theorem can be used to write the two vertical-like disks as graphs of functions over the second coordinate. Thus $f_{n+1}(s,z)=(\varphi^{n+1}_s(z),z)$ where $\varphi^{n+1}_s$ is a holomorphic function, continuous with respect to $s$. The map $f_{n+1}$ is injective.

\begin{remarka} 
This procedure can be used to define external rays for the \He map. External rays are very useful tools, because they give combinatorial models for the Julia set. In \cite{BS6} and \cite{BS7} it was shown 
that external rays for polynomial diffeomorphisms of $\C^2$ can be defined when $J$ is connected. A priori we do not know that our family has connected $J$, but this will be shown to be true as a result of our construction, in Corollary \ref{cor:connected}. The construction of the space $\mathcal{F}$ and of the operator $F$ for $t=0$ is given in \cite[Section~11]{RT} and is identical for $t$ real and small. 
\end{remarka}

On the set $V$ we use the modified metric $d_{\mu}$ from Definition \ref{defn: mu-metric}. 
On the function space $\mathcal{F}_{a,t}$ we
consider the metric
\begin{eqnarray}\label{eq:d}
 d(f, g) = \sup_{s\in \s}\sup_{z\in \D_{r}} d_{\mu}\left(f(s,z),\ g(s,z)\right),
\end{eqnarray}
where 
$d_{\mu}\left(f(s,z),\ g(s,z)\right)$ is defined in \eqref{eq:d-pol2henon} as the infimum of the
length of horizontal rectifiable paths $\gamma:[0,1]\rightarrow V$
with $\gamma(0)=f(s,z)$ and $\gamma(1)=g(s,z)$. 
The length is
measured with respect to the metric $\mu$.

\begin{thm}\label{prop:strict} Suppose that $t\in [-\delta', \delta']$ and
$|a|<\delta$. If $t\neq 0$, then the operator $F_{a,t}:\mathcal{F}_{a,t}\rightarrow \mathcal{F}_{a,t}$
is a strong contraction, i.e. there exists a constant $K_t>1$, which
depends on $t$, such that
    \[
    d(F_{a,t}(f),F_{a,t}(g))<\frac{1}{K_t}\, d(f,g), \mbox{ for any } f,g\in
    \mathcal{F}_{a,t}.
    \]
\end{thm}
\proof We can use the expanding properties of the infinitesimal pseudo metric $\mu$ constructed in \eqref{defn: mu-metric} to show that the operator $F_{a,t}$ contracts distances between vertical-like disks with respect to the induced metric $d_{\mu}$.
We show that there exists $K_t>1$ such that
\begin{equation}\label{d:fibers}
d\left(F_{a,t}\circ f(s\times \D_r),F_{a,t}\circ g(s\times
\D_r)\right)<\frac{1}{K_t}d\left(f(2s\times \D_r),g(2s\times
\D_r)\right),
\end{equation}
for all $f,g\in \mathcal{F}_{a,t}$ and $s\in \s^{1}$.

We first discuss the strategy in the semi-parabolic case $t=0$, which is harder and treated in \cite{RT}. 
When $t=0$ we showed a similar inequality  in \cite[Proposition~11.9]{RT}: for $f,g\in \mathcal{F}_{a,0}$ and $s\in \s^1$ there exists a constant $0\leq C(f,g,s)<1$ such that
\begin{equation*}
d\left(F_{a,0}\circ f(s\times \D_r), F_{a,0}\circ g(s\times \D_r)\right)<C(f,g,s) d\left(f(2s\times \D_r),g(2s\times \D_r)\right).
\end{equation*}
The contraction factor $C(f,g,s) $ depends only on the distance from the fibers $f(2s\times \D_r)$, $g(2s\times\D_r)$ to $W^{ss}_{loc}(\fq_{a,0})$ and goes to $1$ precisely when these fibers approach $W^{ss}_{loc}(\fq_{a,0})$, the local stable manifold of the semi-parabolic fixed point $\fq_{a,0}$. 

We briefly explain how the factor $C(f,g,s) $ is obtained when $t=0$. 
When the disks $f(2s\times \D_r)$, $g(2s\times\D_r)$  are close to $W^{ss}_{loc}(\fq_{a,0})$, they become almost vertical (the vertical cones have angle opening of $\approx |x|^{2q}$, where $|x|$ measures the distance to the stable manifold $W^{ss}_{loc}(\fq_{a,0})$ and is close to $0$). Meanwhile, by \cite[Proposition~6.8]{RT}, the expansion factor of the derivative of the \He map in the horizontal direction is at least $(1+\frac{\epsilon_1}{16} |x|^{q})>1$. By using the fact that $|x|^{q}$ dominates $|x|^{2q}$ when $|x|$ is small, we showed in  \cite[Theorem~10.2]{RT}  that the factor $C(f,g,s)$ is strictly smaller than $1$ and goes to $1$ precisely when $x\rightarrow 0$.  When the disks $f(2s\times \D_r)$, $g(2s\times\D_r)$  do not belong to a small neighborhood of $W^{ss}_{loc}(\fq_{a,0})$ we proved in \cite[Proposition~11.9]{RT}  that in fact we have a strong contraction factor $C(f,g,s)<1/(k+\bigO(|a|))<1$.  The constant $k>1$ is defined in Equation \eqref{eq:k}.

We now return to the proof of inequality \eqref{d:fibers}. When $t\neq 0$, the same proof as outlined above works, but the computations are greatly simplified, because the expansion factor
from Lemma \ref{lemma: kt} is  at least $\min\left(1+\epsilon_2|t|, k\right)$, hence strictly greater than $1$. 

In a small neighborhood of $W^{ss}_{loc}(\fq_{a,t})$, the disks $f(2s\times \D_r)$, $g(2s\times\D_r)$ are almost vertical. Indeed, by Definition \ref{def:vcones} and the invariance of vertical cones from Proposition \ref{prop:cones-Ve}, the vertical cones in the repelling sectors $W^-_B$ have a narrow angle opening $\approx |x|^{2q}$ when $x$ is close to $0$.  By Proposition \ref{prop:cones-Ho} and Lemma \ref{lem:expansion}, the derivative $DH_{a,t}$ expands horizontally by a factor of $(1+\epsilon_2 |t|)\left(1+\frac{\epsilon_1}{16}|x|^{q}\right)$, so the operator $F_{a,t}$ contracts the distance between vertical-like disks by a factor of $(1+\epsilon_2 |t|) C(f,g,s)$. Away from the local  stable manifold $W^{ss}_{loc}(\fq_{a,t})$, we have the same strong contraction factor $1/(k+\bigO(|a|))$ as in the case $t=0$. Let $K_t$ be $\min\left(1+\epsilon_2|t|, k+\bigO(|a|)\right)$.
 In conclusion, when $t\neq 0$, we have a strong contraction factor $1/K_t$, strictly less than $1$.
\qed

As in the polynomial case, we can reduce the hyperbolic
case $t\in [-\delta', \delta']$ to the semi-parabolic case by considering
\begin{equation*}
h:[0,\infty)\rightarrow [0,\infty),\ \ h(s):=\sup\limits_{t\in
[-\delta',\delta']}h_t(s),
\end{equation*}
where 
\begin{eqnarray*}
h_t(s)&:=&\sup\limits_{|a|\leq \delta}\sup \left\{d(F_{a,t}\circ f(\theta\times \D_r),F_{a,t}\circ g(\theta\times
\D_r))\ :\ f,g\in \mathcal{F}_{a,t}\ \mbox{and}\ \theta\in \s^1\right.\\
&&\left.\hspace{1.85cm}\mbox{and}\ d\left(f(2\theta\times\D_r),g(2\theta\times\D_r)\right)\leq s \right\}.
\end{eqnarray*}

By definition, for each $t$, the function $h_{t}:[0,\infty)\rightarrow [0,\infty)$ is increasing and satisfies
    \[
    d(F_{a,t}(f),F_{a,t}(g))<h_{t}\left(d(f,g)\right), \mbox{ for any } f,g\in \mathcal{F}_{a,t}.
    \]
By Theorem \ref{prop:strict} we know  that for $t\neq0$, $h_{t}(s)<\frac{1}{K_{t}}s<s$ for all $s>0$. When $t=0$, we 
know a bit more: by \cite[Theorem~11.10]{RT}, $h_{0}(s+)<s$ for all $s>0$. We can therefore apply Lemma \ref{lemma: Browder-function-h} and conclude that the function $h^{+}:s\mapsto h(s+)$ is
a Browder function that works for all $t\in [-\delta',\delta']$. This
proves that the construction can be done uniformly with respect to
$t$ and $a$. Uniformity with respect to $|a|<\delta$ was already shown in
the semi-parabolic case $t=0$ \cite{RT}. Now Browder's Theorem \ref{thm: Browder} proves
the existence of a unique fixed point $f^{*}_{a,t}:\s^1\times\D_r
\rightarrow \overline{V}$ of the operator $F_{a,t}$.  
We summarize below some basic properties of the fixed point $f^{*}_{a,t}$, which are direct consequences of our construction.

\begin{prop}\label{prop:fixedpoint} The operator $F_{a,t}$ has a unique fixed point
$f^*_{a,t}:\s^1\times \D_r\rightarrow J_{a,t}^+\cap \overline{V}$,
$f^*_{a,t}(s,z)=(\varphi_{t,s}(z),z)$. The map $f_{a,t}^{*}$ is surjective, continuous
with respect to $t$ and $s$, and holomorphic with respect to $a$ and $z$.
\end{prop}

\subsection{Stability and continuity of $J$ and $J^{+}$}\label{subsec:JJ+}

The Julia sets $J$ and $J^+$ depend lower-semicontinuously on the parameters, and discontinuities can occur at a parameter for which the \He map has a semi-parabolic fixed point \cite{BSU}. In Theorem \ref{thm:continuity-2} we prove that in our family of complex \He maps $H_{a,t}$ the sets $J$ and $J^+$ depend continuously on the parameters as $t\rightarrow 0$. 

We begin by analyzing the properties of the fixed point $f^{*}_{a, t}$ in more detail.
The analysis is similar to \cite[Section~12]{RT}, but the role of the parameter $t$ is different.
By Proposition \ref{prop:fixedpoint}, $f^{*}_{a,t}(s,z) = (\varphi_{t,s}(z),z)$, where $\varphi_{t,s}(z)$
is continuous with respect to $s\in \s^{1}$ and analytic with
respect to $z\in \D_{r}$. The map $\varphi_{t,s}$ depends analytically on the
parameter $a$ as well, but we choose to disregard this to simplify notations. We will point out the
dependency on $a$ when needed. For each $t$, let $\sigma_{a,t}:\s^{1}\times
\D_{r}\rightarrow \s^{1}\times \D_{r}$ be defined by
\begin{equation}\label{eq:sigma}
\sigma_{a,t}(s,z) = \left( 2s, a\varphi_{t,s}(z) \right).
\end{equation}
By Proposition \ref{prop:Prop0} below, for sufficiently small $|t|$ and $|a|\neq 0$,
the map $\sigma_{a,t}$ is well-defined, open, and injective. 

As in \cite[Lemma~12.2]{RT}, for each $t$, the map $\varphi_{t,s}$ has the following expansion
\begin{equation}\label{eq:phiOa2}
\varphi_{t,s}(z)=\gamma_t(s)-\frac{az}{2\gamma_t(s)}+a^{2}\beta_{t}(s,z,a),
\end{equation}
where $\gamma_{t}:\s^{1}\rightarrow J_{p_{t}}$ is the Carath\'eodory loop associated to the polynomial $p_{t}$. Recall that $\gamma_t$ is continuous, surjective, and does not vanish on $\s^1$ since the critical point of the polynomial $p_t$ does not belong to the Julia set $J_{p_t}$.
The tail $\beta_{t}(s,z,a)$ is bounded with respect to $a$ and $t$.

\begin{prop}\label{prop:Prop0}
For sufficiently small
$ |t|<\delta'$ and $0<|a|<\delta$ the map $\sigma_{a,t}$ is open and
injective. Moreover $\sigma_{a,t}(\s^{1}\times \D_{r})\subset \s^{1}\times
\D_{|a|r'}$, for some $r'<r$.
\end{prop}
\proof Since the critical point of $p_t$ is far away from the Julia set $J_{p_t}$, there exists $\epsilon>0$ such that 
$|\gamma_t(s)-\gamma_t(s+1/2)|>\epsilon$ for all $s\in\s^{1}$ and $t\in [-\delta',\delta']$. The expansion from \eqref{eq:phiOa2} shows that there exist constants $M_t>0$ such that
$|\varphi_{t,s}(z)-\gamma_t(s)|<|a|M_t$ for all $s\in \s^{1}$ and
$z\in \D_{r}$. The constant $M:=\sup_{|t|<\delta'}M_t$ does not
depend on $t$.  Then for $|a|<\frac{\epsilon}{2M}$ the map $\sigma_t$
is injective. It is open because locally it is a homeomorphism.
\qed

The following theorem is a direct consequence of our construction so far.

\begin{thm}\label{thm:diagram}
\justifying Let $\lambda= e^{2 \pi i p/q}$ and
$\lambda_{t}=(1+t)\lambda$. 
There
exists $\delta,\delta'>0$  such that
\begin{itemize}
\item[\small$\bullet$] for all $-\delta'< t<\delta'$ and
\item[\small$\bullet$] for all parameters $(c,a)\in \mathcal{P}_{\lambda_{t}}$
with $0<|a|<\delta$
\end{itemize}
the diagram commutes:
\[
\diag{\s^{1}\times \D_{r}}{J^{+}\cap \D_{r}\times\D_{r}}{\s^{1}\times
\D_{r}}{J^{+}\cap \D_{r}\times\D_{r}}
{f^{*}_{a,t}}{\sigma_{a,t}}{H_{c,a}}{f^{*}_{a,t}}
\]
\end{thm}
\proof The existence of the fixed point $f^{*}_{a,t}$ has already
been established in the previous section \ref{sec: FunctionSpaceHenon}. By construction, 
we have that $H\circ
f^{*}_{a,t}(s\times \D_{r})$ is compactly contained in $f^{*}_{a,t}(2s\times
\D_{r})$. Thus we can write
\begin{eqnarray*}
H\circ
f^{*}_{a,t}(s,z)&=&\hvec{p_{t}(\varphi_{t,s}(z))+a^{2}w+az}{a\varphi_{t,s}(z)}\\ 
&=& \hvec{\varphi_{t, 2s}(a\varphi_{t,s}(z))}{a\varphi_{t,s}(z)} = f^{*}_{a,t}\circ
\sigma_{a,t} (s,z).
\end{eqnarray*}
The last equality follows from the fact that 
\begin{equation*} 
f^{*}_{a,t}\circ \sigma_{a,t}(s,z) = f^{*}_{a,t}(2s,a\varphi_{t,s}(z)) =
\left(\varphi_{t,2s}(a\varphi_{t,s}(z)),a\varphi_{t,s}(z)\right).
\end{equation*}
Therefore $f^{*}_{a,t}$ semiconjugates $H$ on $J^{+}\cap \overline{V}$ to $\sigma_{a,t}$ on
$\s^{1}\times \D_{r}$, as claimed. The fact that $J^{+}\cap \overline{V} = J^{+}\cap \D_{r}\times\D_{r}$ follows from Lemma \ref{lemma:nbdPolydisk}.
\qed

\begin{figure}[htb]\label{fig:param}
\begin{center}
\mbox{\subfigure[$t=0.25$]{
\begin{overpic}[scale=0.33]{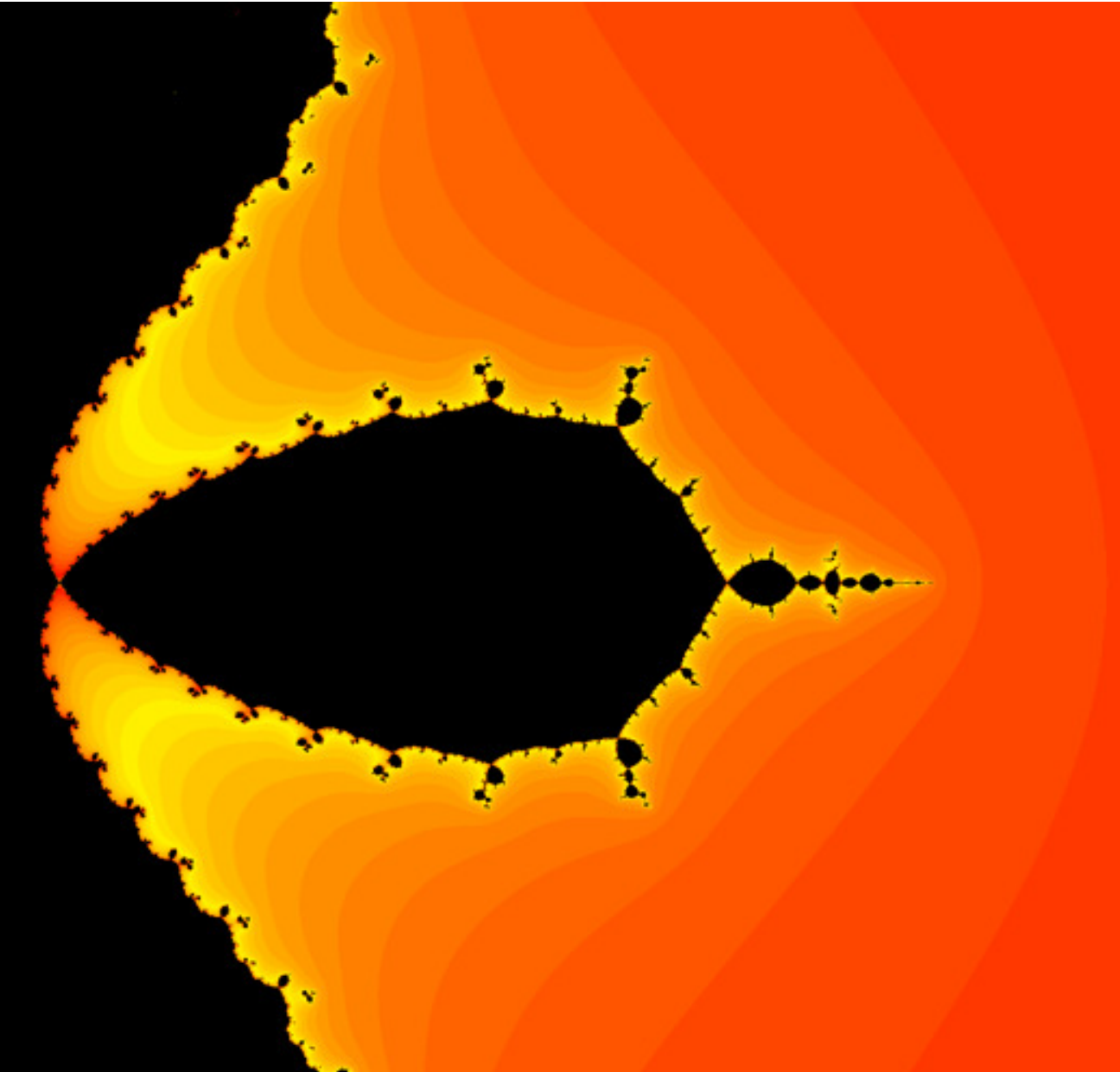}
 \put (41,43) {\color{white} {\tiny $\bullet$} {\scriptsize $0$} }
\end{overpic}}
\quad
\subfigure[$t=0.1$]{
\begin{overpic}[scale=0.33]{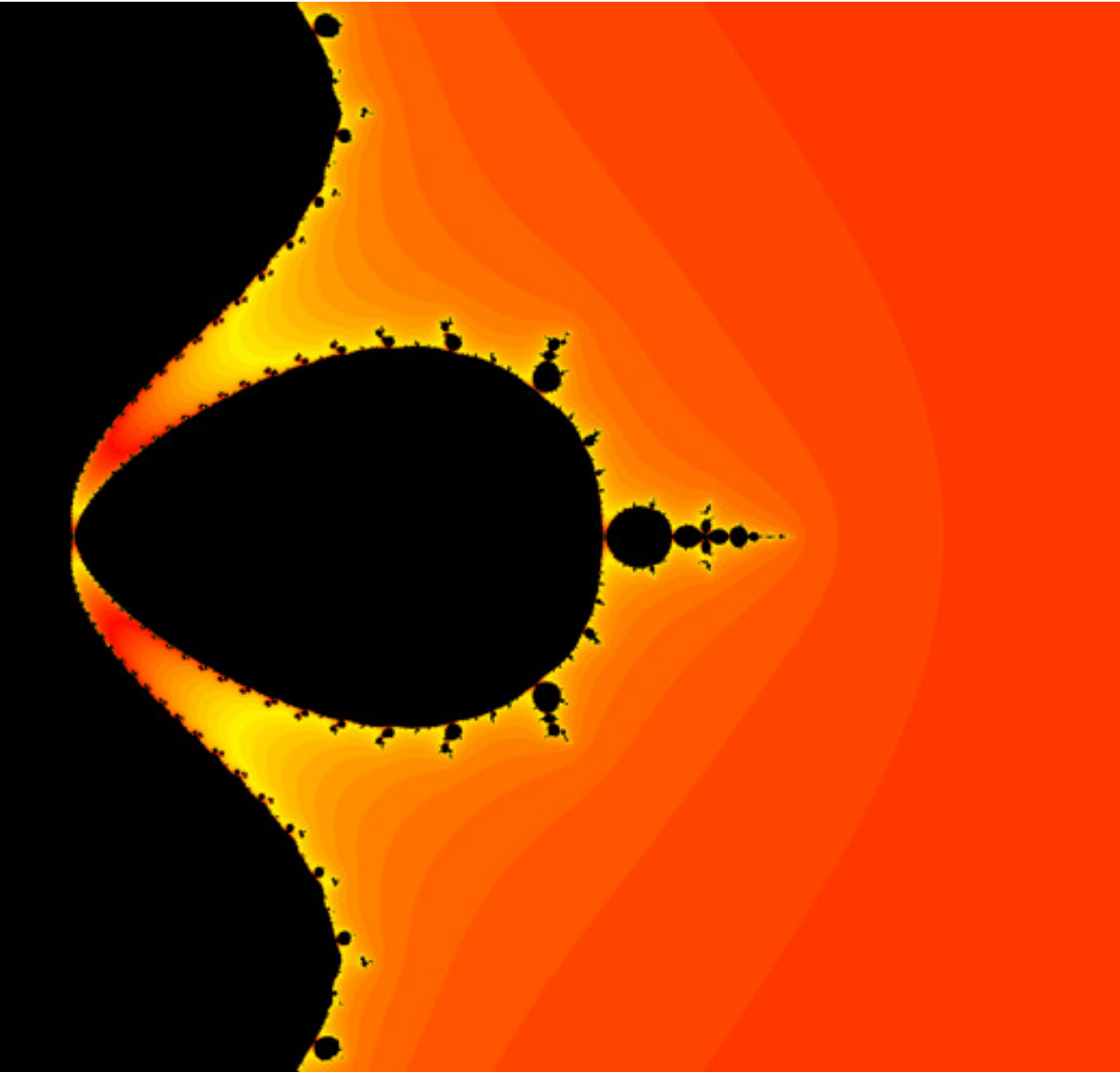} 
\put (35,46) {\color{white} {\tiny $\bullet$} {\scriptsize $0$} }
\end{overpic}}
\quad
\subfigure[$t=0.025$]{
\begin{overpic}[scale=0.33]{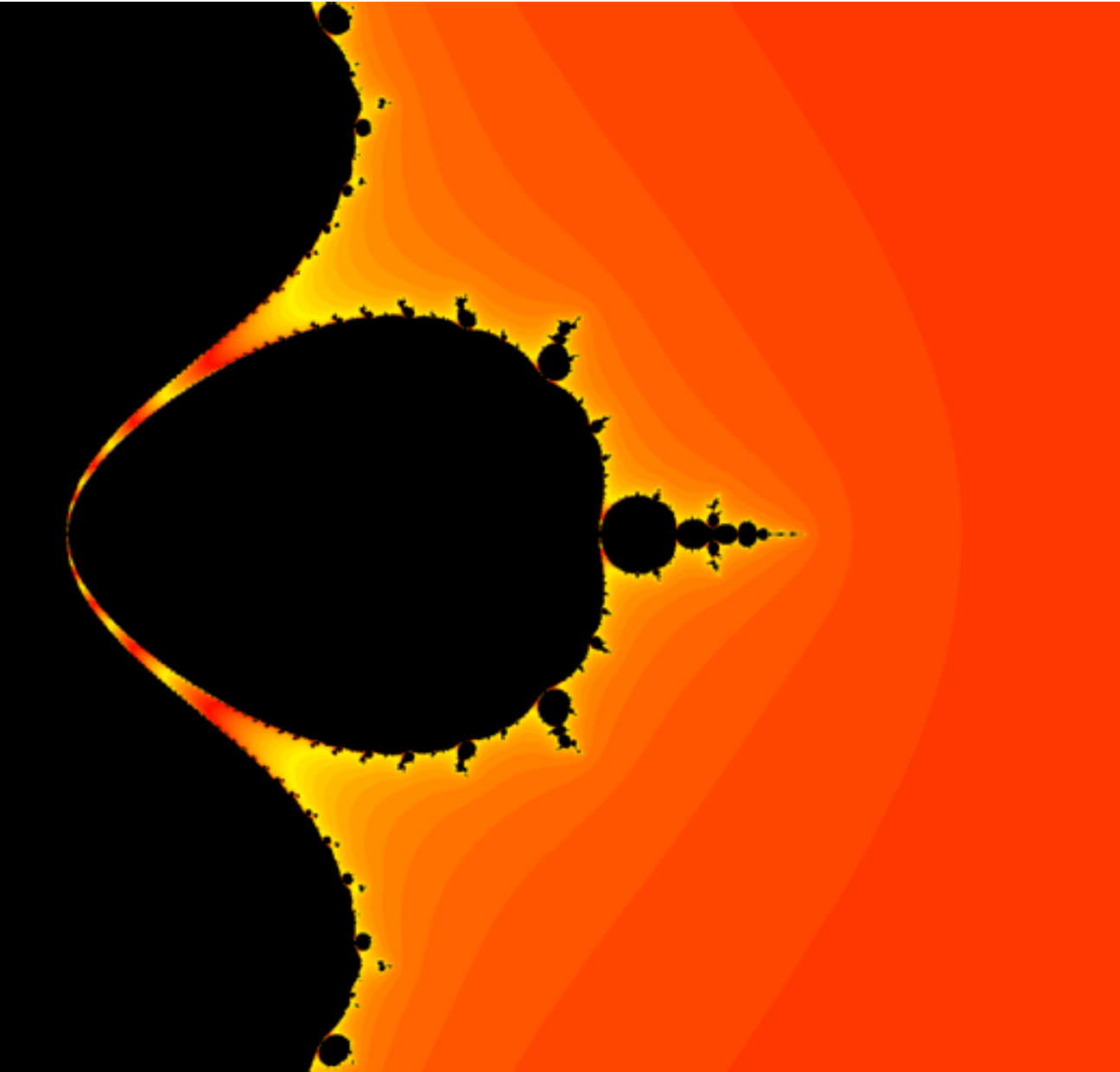} 
\put (35,46) {\color{white} {\tiny $\bullet$} {\scriptsize $0$} }
\end{overpic}}
}
\mbox{\subfigure[$t=0$]{
\begin{overpic}[scale=0.33]{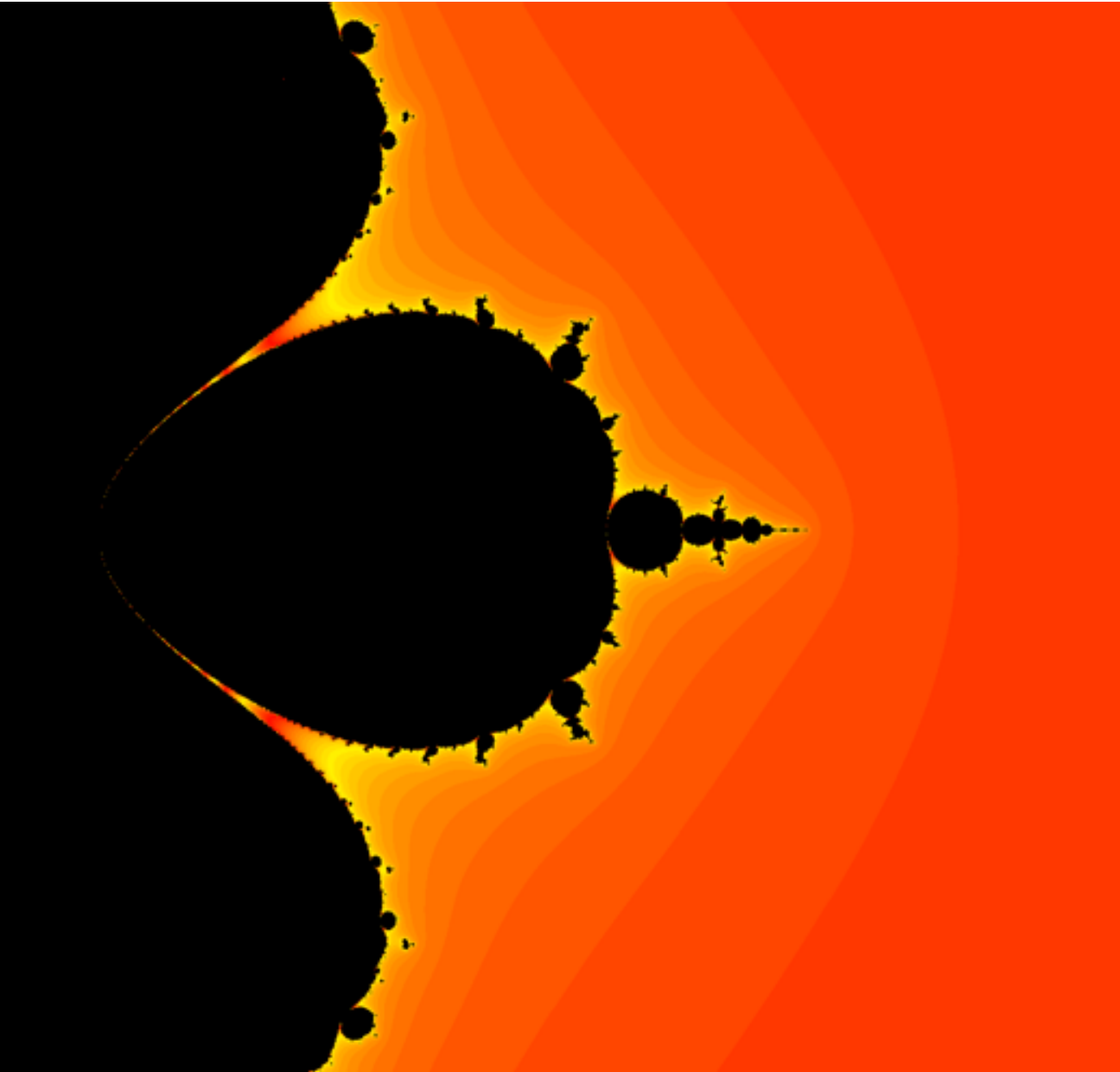}
 \put (37,47) {\color{white} {\tiny $\bullet$} {\scriptsize $0$} }
\end{overpic}
}
\quad
\hspace{-0.25cm}
\subfigure[$t=-0.025$]{
\begin{overpic}[scale=0.33]{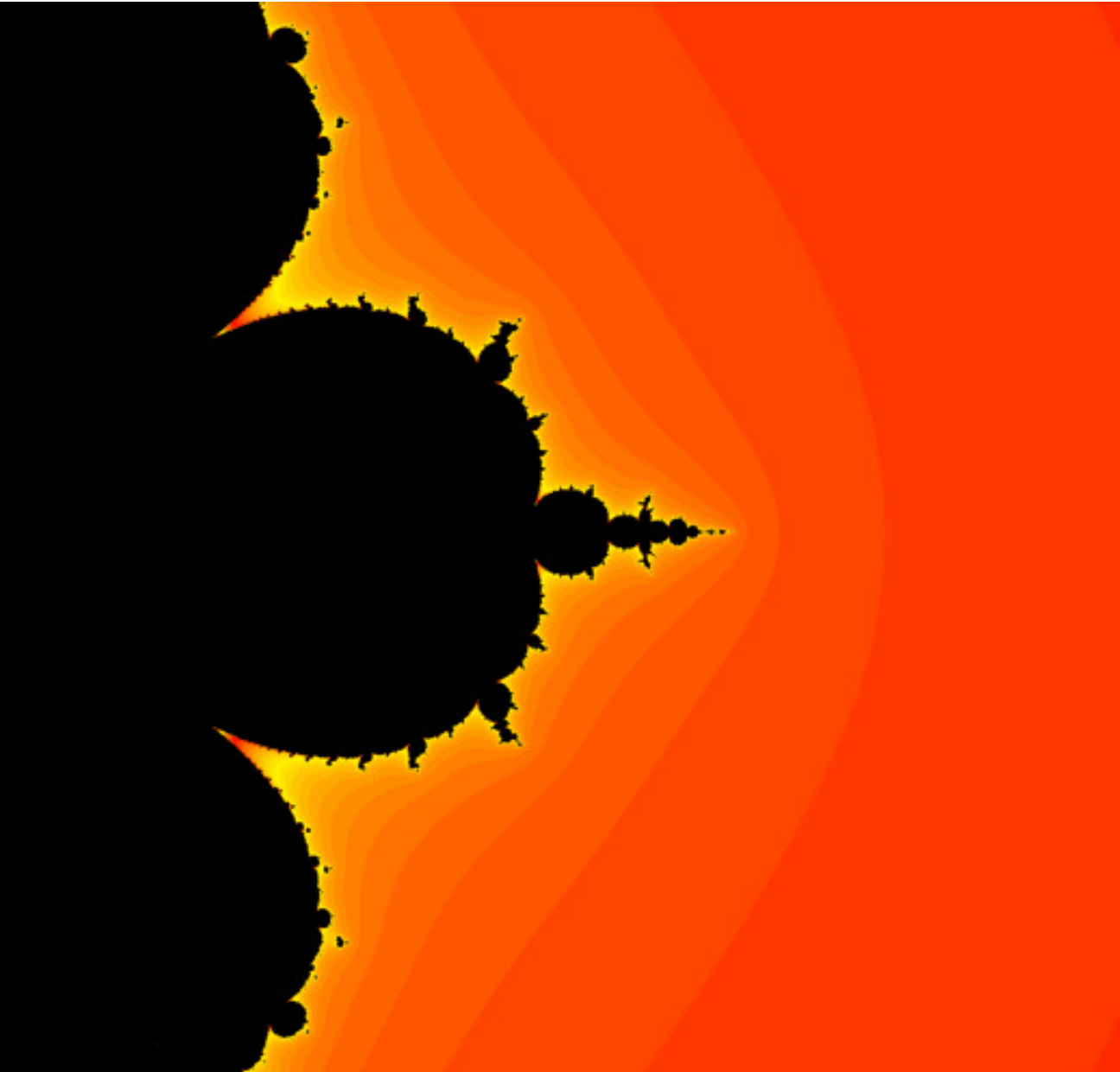} 
 \put (33,47) {\color{white} {\tiny $\bullet$} {\scriptsize $0$} }
\end{overpic}
}
\quad
\hspace{-0.25cm}
\subfigure[$t=-0.1$]{
\begin{overpic}[scale=0.33]{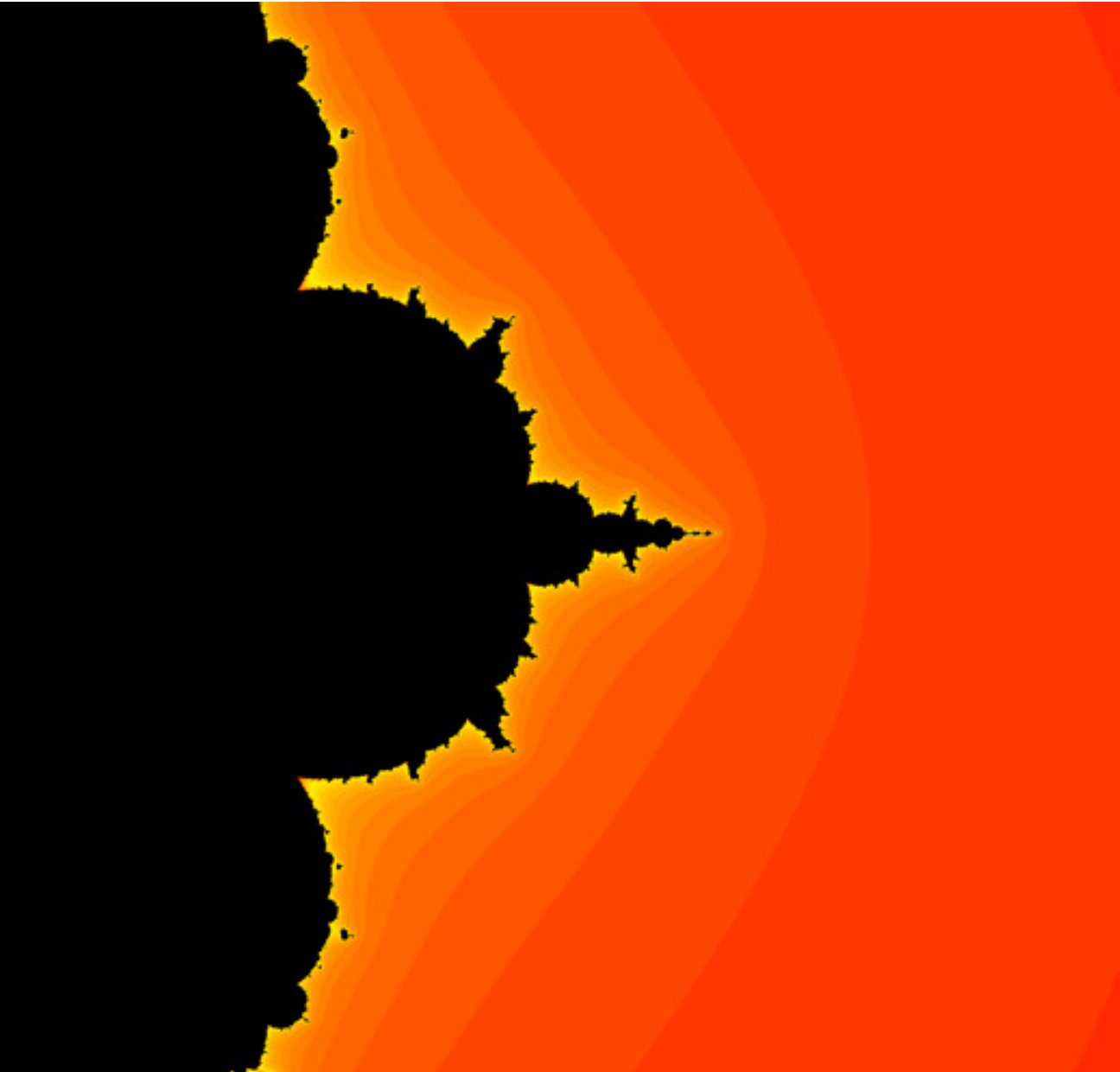} 
 \put (33,47) {\color{white} {\tiny $\bullet$} {\scriptsize $0$} }
\end{overpic}
}
}
\end{center}
\caption{Parameter plots inside the curves $\mathcal{P}_{(1+t)\lambda}$ for $\lambda=-1$ and several values of $t$. In each picture, the large region in the center contains the disk $|a|<\delta$. The black region represents a rough approximation of the set of parameters $(c,a)\in \mathcal{P}_{(1+t)\lambda}$ for which the Julia set $J_{c,a}$ is connected. Here the \He map is written in the standard form $H_{c,a}(x,y)=(x^{2}+c-ay,x)$. 
The pictures were generated using FractalStream.
}
\label{fig:parameter}
\end{figure}

\begin{cor}\label{cor:connected} The Julia set $J$ is connected.
\end{cor} 
\proof By Theorem \ref{thm:diagram} and  Lemma \ref{lemma:nbdPolydisk}, we get
\begin{equation}\label{eq:JJ}
J=f^{*}_{a,t}\left(\bigcap_{n\geq 0}\sigma_{a,t}^{\circ
n}\left(\s^{1}\times \D_{r}\right)\right).
\end{equation}
By Proposition \ref{prop:Prop0}, the intersection above is a nested intersection of connected, relatively compact sets, hence connected. Then $J$ is connected, since $f_{a,t}^{*}$ is continuous. See Figure \ref{fig:parameter} for parameter space pictures of the connectivity region. 
\qed

The \He map $H_{a,t}$ has a fixed point $\fq_{a,t}$ with eigenvalues $\lambda_{t}$ and $\nu_{t}$. The product of the eigenvalues equals the Jacobian of the map, so $|\lambda_{t}||\nu_{t}|= |a|^{2}$. We write $J_{(\lambda_{t},\nu_{t})}$ and $J^{+}_{(\lambda_{t},\nu_{t})}$ to denote the dependency of the Julia sets $J$ and $J^{+}$ on the eigenvalues, rather than on the parameters $a$ and $t$. 

\begin{thm}[\textbf{Continuity}]\label{thm:continuity-2}
There exists $\delta>0$ such that if $|\nu_t|<\delta$ and $\nu_t\rightarrow \nu$ as $t\rightarrow 0$, then 
the Julia sets $J$ and $J^+$ depend continuously on the parameters, i.e. 
\[
    J^{+}_{(\lambda_{t},\nu_t)}\rightarrow J_{(\lambda,\nu)}^{+}\ \ \ \ \ \mbox{and}\ \ \ \     J_{(\lambda_t,\nu_t)}\rightarrow J_{(\lambda,\nu)}
\]
in the Hausdorff topology.
\end{thm}
\proof By Theorem \ref{thm:diagram}, Lemma \ref{lemma:nbdPolydisk}, and Proposition \ref{prop:fixedpoint}
we know that
\begin{equation}\label{eq:aux1}
J^{+}_{(\lambda_{t},\nu_t)}\cap \D_{r}\times\D_{r}=f^{*}_{a,t}\left(\s^1\times \D_r\right)
\end{equation}
and $f^{*}_{a,t}$ is continuous with respect to $t$ and holomorphic in $a$. Therefore  
\[
J^{+}_{(\lambda_{t},\nu_t)}\cap \D_{r}\times\D_{r}\rightarrow J^{+}_{(\lambda,\nu)}\cap \D_{r}\times\D_{r}
\]
in the Hausdorff topology, as $t\rightarrow 0$. Let $H= H_{(\lambda_{t},\nu_{t})}$ be the \He map corresponding to a pair of eigenvalues $(\lambda_{t},\nu_{t})$. Clearly $H^{-1}$ is continuous with respect to $t$.  Let $n$ be a positive integer. Taking $H^{-\circ n}$ in Equation \eqref{eq:aux1} gives  
\[
H^{-\circ n}\left(J^{+}_{(\lambda_{t},\nu_t)}\cap \D_{r}\times\D_{r}\right)=H^{-\circ n}f^{*}_t\left(\s^1\times \D_r\right),
\]
which converge in the Hausdorff topology, as $t\rightarrow 0$. We have accounted for all of $J^{+}$, because globally the set $J^{+}$ is $\bigcup_{n\geq0}H^{-\circ n}(J^{+}\cap \D_{r}\times\D_{r})$.

The Julia set $J_{(\lambda_{t},\nu_{t})}$ can be written as in Equation \eqref{eq:JJ}. The maps $f^{*}_{a,t}$ and $\sigma_{a,t}$ are continuous in $a$ and $t$, so $J_{(\lambda_t,\nu_t)}$ converges to $J_{(\lambda,\nu)}$ in the Hausdorff topology. 
\qed

\begin{remarka}\label{rem:obs-radial} We have established a continuity result for real values of $t$, but the situation is much more general, similar to the one-dimensional case. If $t$ is real, then the local attractive/repelling sectors from Section \ref{subsec:horizontal} are ``straight'', as in the semi-parabolic case $t=0$. If we allow $t$ to be complex, then we need to adapt the computations from Sections \ref{subsec:horizontal} and \ref{subsec:hyper-region} to construct ``spiralling petals" for the \He map. 
\end{remarka}

Suppose $t$ is fixed. For each $s\in\s^1$, $f^{*}_{a,t}(s\times \D_r)$ is a vertical-like holomorphic disk. Any two such disks corresponding to distinct angles $s_1$ and $s_2$ are either disjoint or coincide (since they were obtained as a uniform limit of disjoint holomorphic disks $f_{n}(s_1\times \D_r)$ and $f_{n}(s_2\times \D_r)$). In the latter case, their parametrizing maps coincide, i.e. \[f^{*}_{a,t}(s_1,z)=(\varphi_{s_1}(z),z)=(\varphi_{s_2}(z),z)=f^*_{a,t}(s_2,z)\] for all $z\in \D_r$. The fixed point $f^{*}_{a,t}$ is holomorphic with respect to
$a$, so we can determine the equivalence classes of $f^{*}_{a,t}$ by letting
$a\rightarrow 0$. When $a=0$,  we have $J^+ \cap \overline{V} = J_{p_t}\times
\D_r$ so all the identifications are given by the polynomial $p_t$.
An application of Hurwitz's theorem (see \cite[Propositions~12.4-12.6]{RT}) gives
\begin{equation*}\label{eq: equiv-classes}
f^{*}_{a,t}(s_1,z_1)=f^*_{a,t}(s_2,z_2) \mbox{ if and only if }
\gamma_t(s_1)=\gamma_t(s_2) \mbox{ and } z_1=z_2.
\end{equation*}
\begin{defn}\label{def:equiv-rel}
We define an equivalence relation $\sim$ on $\s^{1}\times \D_{r}$  as follows: 
\begin{equation*}\label{eq: equiv-rel}
(s_{1},z)\sim(s_{2},z)\ \ \ \mbox{whenever}\ \  \ \gamma_t(s_{1})=\gamma_t(s_{2}). 
\end{equation*}
\end{defn}
We obtain $\varphi_{t,s_{1}}(z) = \varphi_{t,s_{2}}(z)$  iff $\gamma_t(s_{1}) = \gamma_t(s_{2})$. By Equation \eqref{eq:phiOa2} this also gives $\beta_{t}(s_{1},z,a,t) =\beta_{t}(s_{2},z,a,t)$ whenever $\gamma_t(s_{1}) = \gamma_t(s_{2})$. The relation $\sim$ is clearly closed. 
Moreover, since all polynomials $p_t$ have the same Thurston lamination for $t\in[0,\delta')$ hence the same combinatorial model (see \cite{Th}), the
equivalence relation $\sim$ does not depend on $t$ when $t\in [0,\delta')$. Thus, in 1-D, the polynomial $p_{t}$ acting on $J_{p_{t}}$ is conjugate to the parabolic polynomial $p_{0}$  on $J_{p_{0}}$, for all $t\in [0,\delta')$. Note that this is not true for all $t\in (-\delta',\delta')$. For example, if $\lambda =-1$, the Julia set of $p_{1-t}$ is a quasi-circle and the associated Thurston lamination is empty. However, the Julia set of $p_{1+t}$ is homeomorphic to the Julia set of $z\mapsto z^{2}-3/4$ (the ``fat Basilica'') and the corresponding lamination is non-empty.  The same situation is true in 2-D as we will show below.

\begin{thm}[\textbf{Stability}] The family of complex \He maps $\mathcal{P}_{\lambda_t} \ni (c,a) \to  H_{c,a}$   is a structurally stable family  on $J$ and $J^{+}$ for  $0<|a| < \delta$ and $0\leq t<\delta'$.
\end{thm}
\proof
In view of Equation \eqref{eq:phiOa2}, the map $\sigma_{a,t}:\s^{1}\times \D_{r}\rightarrow \s^{1}\times \D_{r}$ has the form 
\begin{equation}\label{eq:sigma-at}
\sigma_{a,t} (s,z)= \left(2s, a\gamma_t(s)-\frac{a^{2}z}{2\gamma_t(s)}+\bigO(a^3) \right).
\end{equation}
By Theorem \ref{thm:diagram}, the map $H_{a,t}$ on $J^{+}_{a,t}\cap \D_{r}\times\D_{r}$ is semiconjugate to  $\sigma_{a,t}$ on $\s^{1}\times \D_{r}$. For $0<|a|<\delta$ and $t\in[0,\delta')$ small enough, the maps $\sigma_{a,t}$ are conjugate to each other. The proof of this fact is the same as that of \cite[Lemmas~12.7, 12.8]{RT} stated below. 

\begin{lemma}[\cite{RT}]\label{lemma:conj} Suppose $0<|a|<\delta$ and $t=0$.
\begin{itemize}
\item[a)] The map $\sigma_{a,0}:\s^{1}\times \D_{r}\rightarrow \s^{1}\times \D_{r}$ is conjugate to the map $\sigma_{a,0}':\s^{1}\times \D_{r}\rightarrow \s^{1}\times \D_{r}$, defined by $\sigma'_{a,0}(s,z)= \left(2s, a\gamma_0(s)-\frac{a^{2}z}{2\gamma_0(s)}\right)$.
\item[b)] The maps $\sigma'_{a,0}$ are conjugate to $\sigma'_{\epsilon,0}$ for some $\epsilon>0$ independent of $a$. 
\end{itemize}
\end{lemma}

It is important to note that in each fiber $\{s\}\times\D_{r}$ the image of $\sigma_{a,t}$ consists of two disjoint disks.  This follows from Proposition \ref{prop:Prop0} as $\gamma_{t}(s)$ and $\gamma_{t}(s+1/2)$ are at least $\epsilon$-apart, for some $\epsilon>0$ independent of $a$ and $t$. 

Suppose $0<|a|<\delta$ and $t\in [0,\delta')$. The equivalence classes of $f_{a,t}^{*}$ are exactly the ones given by the equivalence relation $\sim$, in the sense that $(s_{1},z)\sim(s_{2},z)$ iff $f_{a,t}^{*}(s_{1},z) = f_{a,t}^{*}(s_{2},z)$. Moreover, by Definition \ref{def:equiv-rel} and the discussion following it, if $(s_{1},z) \sim (s_{2},z)$ then $\sigma_{a,t}(s_{1},z)\sim \sigma_{a,t}(s_{1},z)$. 
Hence, $f_{a,t}^{*}$ and $\sigma_{a,t}$ are well defined on the quotient $\s^{1}\times \D_{r}/_{\sim}$ and  $f_{a,t}^{*}:\s^{1}\times \D_{r}/_{\sim} \rightarrow J_{a,t}^{+}\cap \D_{r}\times\D_{r}$  is bijective. The equivalence relation
does not depend on $t$ or $a$. We get that $(H_{a,t}, J^{+}_{a,t}\cap \D_{r}\times\D_{r})$ is conjugate to $(\sigma_{a,t}, \s^{1}\times \D_{r}/_{\sim})$, which are conjugate to each other and to $(\sigma'_{\epsilon,0}, \s^{1}\times \D_{r}/_{\sim})$, for all $0<|a|<\delta$ and $t\in [0,\delta')$. Stability on $J$ and $J^{+}$ follows from these observations.
\qed

Using the same arguments as in the previous theorem, we also have that the family $\mathcal{P}_{\lambda_t} \ni (c,a) \to  H_{c,a}$ is a structurally stable family on $J$ and $J^{+}$ for  $0<|a| < \delta$ and $t\in (-\delta',0)$. However, this is not so surprising: by Theorem \ref{thm:HypRegion-b} this family of maps is hyperbolic and has connected Julia set $J$, so the family belongs to the same hyperbolic component of the \He connectedness locus. As in the polynomial case presented earlier, stability does not hold for all parameters $t\in (-\delta',\delta')$.  

Using the equivalence relation from Definition \ref{def:equiv-rel} we can identify the quotient
space $\s^{1}\times \D_{r}/_{\sim}$ with $J_{p_t}\times \D_{r}$ and
the map $\sigma_{a,t}:\s^{1}\times \D_{r}\rightarrow \s^{1}\times \D_{r}$ defined in Equation
\eqref{eq:sigma} with a similar map $\psi_{a,t}: J_{p_t}\times \D_{r}\rightarrow J_{p_t}\times \D_{r}$ of the form
\begin{equation}\label{eq:psi-at}
 \psi_{a,t}(\zeta,z)=\left(p_{t}(\zeta),a\zeta -
    \frac{a^{2}z}{2\zeta}+\bigO(a^{3})\right).
\end{equation}
The following theorem is a direct consequence of the construction above and provides concrete model maps for the \He family. 
The corollaries following the theorem are immediate consequences. 

\begin{thm}\label{thm:Parabolas-Pt}
\justifying Let $\lambda= e^{2 \pi i p/q}$ and
$\lambda_{t}=(1+t)\lambda$. Suppose $p_{t}(x)=x^{2}+c_{t}$ is a
polynomial with a fixed point of multiplier $\lambda_{t}$. There
exists $\delta,\delta'>0$  such that
\begin{itemize}
\item[\small$\bullet$] for all $t\in(-\delta',\delta')$ and
\item[\small$\bullet$] for all parameters $(c,a)\in \mathcal{P}_{\lambda_{t}}$ with $0<|a|<\delta$
\end{itemize}
there exists a homeomorphism
$     \Phi_{a,t}: J_{p_{t}}\times \D_{r} \rightarrow J^{+}\cap \D_{r}\times\D_{r} $
which makes the diagram
\[
\diag{J_{p_{t}}\times \D_{r}}{J^{+}\cap  \D_{r}\times\D_{r}}{J_{p_{t}}\times \D_{r}}{J^{+}\cap  \D_{r}\times\D_{r}}
{\Phi_{a,t}}{\psi_{t}}{H_{c,a}}{\Phi_{a,t}}
\]
commute, where
\begin{equation}\label{eq:psi-t}
    \psi_{t}(\zeta,z)=\left(p_{t}(\zeta),\epsilon\zeta -
    \frac{\epsilon^{2}z}{2\zeta}\right)\!,
\end{equation}
for some $\epsilon>0$ independent of $a$ and $t$.
\end{thm}
\proof Most of the work has already been done. The idea of the proof is the same as in \cite[Theorem~1.1]{RT}. As in \cite[Lemma~12.7]{RT} we can construct a homeomorphism $h_{a,t}:J_{p_{t}}\times \D_{r}\rightarrow J_{p_{t}}\times \D_{r}$ conjugating the map $\psi_{a,t}$ from Equation \eqref{eq:psi-at} to the map $\psi_{t}$ from Equation \eqref{eq:psi-t}, for some $\epsilon>0$, independent of $a$ and $t$. The map $\Phi_{a,t}$ is just a composition between the homeomorphism $h_{a,t}$ and the map $f^{*}_{a,t}$ from Theorem \ref{thm:diagram}.
\qed

\begin{cor}\label{cor:J-t}
The Julia set $J$ is homeomorphic to $\bigcap_{n\geq 0}\psi_{t}^{\circ n}(J_{p_{t}}\times \D_{r})$.
\end{cor}

\begin{cor}\label{cor:globagJ+t}
Passing to the inductive limit we obtain a global model for the Julia set $J^{+}$. The map $\Phi_{a,t}$ extends naturally
to a homeomorphism $\widecheck{\Phi}_{a,t}$ which makes the following
diagram
\[
\diag{\limind(J_{p_{t}}\times
\D_{r},\psi_{t})}{J^{+}}{\limind(J_{p_{t}}\times
\D_{r},\psi_{t})}{J^{+}}
{\widecheck{\Phi}_{a,t}}{\widecheck{\psi}_{t}}{H_{c,a}}{\widecheck{\Phi}_{a,t}}
\]
commute.
\end{cor}

\end{document}